\documentclass{amsart}
\usepackage[latin1]{inputenc}
\usepackage{amssymb}
\usepackage{amsmath}
\usepackage{enumerate}
\usepackage{epsfig}

\usepackage[all]{xy}

\newtheorem{theorem}{Theorem}[section]
\newtheorem{lemma}[theorem]{Lemma}
\newtheorem{proposition}[theorem]{Proposition}
\newtheorem{corollary}[theorem]{Corollary}

\newtheorem{definition}[theorem]{Definition}
\newtheorem{example}[theorem]{Example}
\newtheorem{construction}[theorem]{Construction}
\newtheorem{ass}[theorem]{Assumption}
\newtheorem{notation}[theorem]{Notation}

\newtheorem{rem}[theorem]{Remark}

\numberwithin{equation}{section}
\newcommand{\rk}{\mbox{rank}}
\newcommand{\NS}{\mbox{NS}}
\newcommand{\T}{\mbox{T}}
\newcommand{\Tr}{\mbox{Tr}}
\newcommand{\Red}{\mbox{Red}}
\newcommand{\MW}{\mbox{MW}}
\newcommand{\MWL}{\mbox{MWL}}

\newcommand{\la}{\leftarrow}
\newcommand{\ra}{\rightarrow}

\newcommand{\bprf}{{\it Proof.~}}

\newcommand{\eprf}{\hfill $\square$ \smallskip\par}
\newcommand{\erem}{\hfill $\square$}

\newcommand{\PP}{ \mathbb{P}}
\newcommand{\C }{ \mathbb{C}}
\newcommand{\R}{ \mathbb{R}}
\newcommand{\Z}{\mathbb{Z}}
\newcommand{\Q}{\mathbb{Q}}

\newcommand{\N}{\mathbb{N}}

\makeatletter
\def\blfootnote{\xdef\@thefnmark{}\@footnotetext}

\begin{document} 
\title{Van Geemen--Sarti involutions and elliptic fibrations on K3 surfaces double cover of $\mathbb{P}^2$}

\author{Paola Comparin and Alice Garbagnati
\address{Paola Comparin, Universit\'e de Poitiers, 
Laboratoire de Math\'ematiques et Applications, 
 T\'el\'eport 2 
Boulevard Marie et Pierre Curie
 BP 30179,
86962 Futuroscope Chasseneuil Cedex, France}
\email{paola.comparin@math.univ-poitiers.fr}
\address{Alice Garbagnati, Dipartimento di Matematica, Universit\`a di Milano,
  via Saldini 50, I-20133 Milano, Italy}
\email{alice.garbagnati@unimi.it}
\urladdr{http://sites.google.com/site/alicegarbagnati/home}
}

\begin{abstract}
{In this paper we classify the elliptic fibrations on K3 surfaces which are the double cover of a blow up of $\mathbb{P}^2$ branched along rational curves and we give equations for many of these elliptic fibrations. Thus we obtain a classification of the van Geemen--Sarti involutions (which are symplectic involutions induced by a translation by a 2-torsion section on an elliptic fibration) on such a surface. Each van Geemen--Sarti involution induces a 2-isogeny between two K3 surfaces, which is described in this paper.}
\end{abstract}


\maketitle

\pagestyle{myheadings}
\markboth{Paola Comparin and Alice Garbagnati}{van Geemen--Sarti involutions on double covers of $\mathbb{P}^2$}
\section{Introduction}
An involution on a complex K3 surface is called symplectic if it acts trivially on the second cohomology group, otherwise it is called non--symplectic. The quotient of a K3 surface by a non--symplectic involution is a smooth surface (it is either an Enriques surface or a surface with Kodaira dimension $-\infty$). The quotient  of a K3 surface by a symplectic involution is a singular surface whose desingularization is a K3 surface. Thus a symplectic involution on a K3 surface induces a 2-isogeny between two K3 surfaces (a very explicit description of certain isogenies between K3 surfaces induced by symplectic involutions can be found, for example, in \cite{vGT}). If a K3 surface admits a symplectic involution, then its Picard number is at least 9 (cf. \cite{vGS}). Let $X$ be a K3 surface with a sufficiently large Picard number, then one can ask if $X$ admits at least one symplectic involution and, if yes, how many symplectic involutions it admits. To positively answer to the first question it suffices to show that there is a primitive embedding of a certain lattice, $E_8(-2)$, in the N\'eron--Severi group of the surface (cf. \cite{Nikulin symplectic}, \cite{vGS}). To answer to the second question one has to analyze the number of primitive embeddings  of $E_8(-2)$ in the N\'eron--Severi of the surface. In general this is a difficult problem. On the other hand it is of a certain interest to analyze different symplectic involutions on the same K3 surface, because each of them gives a different 2-isogeny with another K3 surface.
In this paper we restrict our attention to particular symplectic involutions, the so called van Geemen--Sarti involutions, which are the translation by a 2-torsion section on an elliptic fibration on the K3 surface. Our purpose is in fact to classify all the van Geemen--Sarti involutions and the associated 2-isogenies on certain K3 surfaces, $X_{(r,22-r,\delta)}$, which are double cover of a blow up of $\mathbb{P}^2$ branched along rational curves. Since the presence of a van Geemen--Sarti involution on a K3 surface assumes the presence of an elliptic fibration on it, our strategy consists in the classification of all the elliptic fibrations on the surfaces $X_{(r,22-r,\delta)}$ and then in the analysis of those admitting at least a  2-torsion section. We underline that the classification of the elliptic fibrations on a family of K3 surfaces is the object of study of several papers (cf. \cite{O},\cite{K}, \cite{SZ}), some of them very recent (cf.\ \cite{BL}, \cite{Ku}). To classify the elliptic fibrations on the K3 surfaces $X_{(r,22-r,\delta)}$, we use techniques which go back to \cite{O} and \cite{K}, where the elliptic fibrations on a K3 surface $X$ are classified if $X$ is respectively the Kummer surface of the product of two non--isogenous elliptic curves or the double cover of $\mathbb{P}^2$ branched along the union of 6 lines in general position. In both these papers one of the main properties of $X$ is that it admits a non--symplectic involution acting trivially on the N\'eron--Severi group. In particular in \cite{K} the author suggests that his results can be generalized to K3 surfaces with such an involution and other particular properties. Here we do exactly this generalization: the K3 surfaces considered are double covers of $\mathbb{P}^2$ branched over singular (and reducible) sextics. The choice of the branch sextic is not unique, and we use this property in order to give a geometric description and an equation for several elliptic fibrations on $X_{(r,22-r,\delta)}$. Indeed, we associate to certain elliptic fibrations on $X_{(r,22-r,\delta)}$ a pencil of rational curves in $\mathbb{P}^2$ whose general member meets the branch sextic in 4 smooth points. Thus the equation of the elliptic fibration can be deduced by the equation of the branching sextic and of the pencil. In several cases an appropriate choice of the branch sextic makes it possible to choose the pencil of rational curves to be a pencil of lines and thus to obtain a very easy and clear geometrical description of the elliptic fibration. Moreover, this technique allow us to show explicitly how certain elliptic fibrations on the surface $X_{(r,22-r,\delta)}$ specialize to elliptic fibrations on $X_{(r+1,21-r,\delta)}$ (this kind of problem was analyzed for a particular elliptic fibration on $X_{(16,6,1)}$ in \cite{CD:vGS involutions}). Neither in \cite{O} nor \cite{K} the equations of the classified elliptic fibrations are found. More recently the equations of the elliptic fibrations classified in \cite{O} were given in \cite{KS} and the equations of the elliptic fibrations classified in \cite{K} were given in \cite{U}.

In Section \ref{subsec: non symplectic involutions} we introduce the K3 surfaces $X_{(r,22-r,\delta)}$ and in Section \ref{section:double cover} we describe them as double covers of $\mathbb{P}^2$ branched over sextics. We denote by $\iota$ the cover involution. In Section \ref{section: possibilities} we give a list of the possible elliptic fibrations on the K3 surfaces $X_{(r,22-r,\delta)}$. From a geometrical point of view there are essentially two possibilities: $\iota$ acts only on the basis of the elliptic fibrations (i.e. it is induced by an involution of $\mathbb{P}^1$, base of the elliptic fibration, see Example \ref{ex: involution on the basis}) or $\iota$ acts only on the fibers of the elliptic fibration (as the hyperelliptic involution, see Example \ref{ex: hyperelliptic involution}). In the first case (analyzed in Section \ref{section:existence case 1}) the elliptic fibrations admit sections of infinite order (if $r\neq 20$), they are obtained from a rational elliptic fibration by a base change and they are not associated to a pencil of rational curves on $\mathbb{P}^2$ as described before. Moreover, there are few elliptic fibrations of this type admitting a van Geemen--Sarti involution. On the contrary, in the second case, the elliptic fibrations have no sections of infinite order and are always associated to a pencil of rational curves in $\mathbb{P}^2$: in Section \ref{sec: existence case 2)} we describe this pencil for each admissible elliptic fibration of this type. In this case there are several elliptic fibrations admitting a 2-torsion section and we observe that the choice of the reducible fibers of the fibrations does not determine the intersection properties of the torsion section: for example, in \cite{K} the author proves that there exists an elliptic fibration on $X_{(16,6,1)}$ with reducible fibers $I_4^*+6I_2$ and a 2-torsion section by giving explicitly the nef divisor $F$ defining this elliptic fibration. Thus the intersection properties of the 2-torsion section are determined. In \cite{CD:vGS involutions} an elliptic fibration on $X_{(16,6,1)}$ with the same reducible fibers is considered: the 2-torsion section does not meet the reducible fibers in the same components as the elliptic fibration considered in \cite{K}. Thus these two elliptic fibrations on $X_{(16,6,1)}$ define two different van Geemen--Sarti involutions. In Section \ref{sec: existence case 2)} we give a complete classification of the van Geemen--Sarti involutions considering also this type of situation, i.e. we do not classify only the type of reducible fibers of an elliptic fibration but also the intersections of the torsion sections with such reducible fibers. The main ingredient to describe the intersection properties of the torsion section with the reducible fiber is the height formula.\\
In Section \ref{section: van Geemen--Sarti involutions} we return to one of our original problems: describe the 2-isogenies between K3 surfaces which are induced by the quotient map of a van Geemen--Sarti involution. Thus we describe the quotients of the elliptic fibrations classified in the previous sections: we recover and describe in a slightly different way some known 2-isogenies (for example the ones analyzed in \cite{vGT} and in \cite{CD:vGS involutions}) and we obtain many other isogenies. It is of a certain interest that there exist some isogenies among surfaces $X_{(r,22-r,\delta)}$ and surfaces $X_{(r,a,\delta)}$, where $a\neq 22-r$. The K3 surfaces $X_{(r,a,\delta)}$ are surfaces admitting a non--symplectic involution acting trivially on the N\'eron--Severi group. They are the object of study for several reasons, one of them is that they are involved in the Borcea--Voisin construction, which allows one to obtain Calabi--Yau 3-folds from these K3 surfaces. If $a\neq 22-r$ and $(r,a,\delta)\neq (14,6,0)$ both the family of K3 surfaces $X_{(r,a,\delta)}$ and the associated Calabi--Yau 3-folds admit a known mirror family. This is not the case if $a=22-r$, thus we are considering an isogeny between K3 surfaces (and Calabi--Yau 3-folds) with a known mirror family and K3 surfaces (and Calabi--Yau 3-folds) for which the mirror is not known.\\

{\bf Acknowledgments }{\it We warmly thank Bert van Geemen and Alessandra Sarti for interesting discussions and useful suggestions.}

\section{Background materials}\label{section: backgroupnd materials}
\begin{definition} Let $X$ be a compact complex surface. It is called  K3 surface if its canonical bundle is trivial and $h^{1,0}(X)=0$.\\
The space $H^{2,0}(X)$ is 1-dimensional and we choose a basis $\omega_X$, which we call the period of $X$.
 \end{definition} We will always assume that $X$ is projective.\\
Let $\alpha$ be an automorphism of $X$. Then $\alpha^*$ preserves the Hodge structure of $H^2(X)$. Thus $\alpha(\omega_X)=\lambda_\alpha(\omega_X)$, with $\lambda_\alpha\in\C^*$. Moreover, if $\alpha$ has finite order $m$ then $\lambda_{\alpha}^m=1$.
\begin{definition} Let $X$ be a K3 surface. An involution $\iota\in Aut (X)$ is called symplectic if $\lambda_\iota=1$, non--symplectic if $\lambda_\iota=-1$.\end{definition}
\subsection{Lattices}
A lattice $(L,b)$ is a free $\Z$-module of finite rank with a $\Z$-valued symmetric bilinear form $b:L\times L\longrightarrow \Z$. When it is clear that the bilinear form is $b$, we denote the lattice $(L,b)$ only by $L$. A lattice $L$ is said to be even if the quadratic form associated to $b$ takes only even values, otherwise it is called odd. The discriminant $d(L)$ of $L$ is the determinant of a matrix which represent $b$ with respect to a certain basis, $L$ is said to be unimodular if $d(L)=\pm 1$. If $L$ is non-degenerate, i.e. $d(L)\not= 0$, then the signature of $L$ is the signature of the $\R$-linear extension of $b$ to $L\otimes\R$. A lattice with signature $(1, \rk(L)-1)$ is called hyperbolic.\\
The dual of the lattice $L$ is $L^{\vee}=Hom_{\Z}(L,\Z)=\{v\in L\otimes_{\mathbb{Z}}\mathbb{Q}~|~b(v,x)\in\mathbb{Z}~\mbox{for~all}~x\in L\}$. There is a natural embedding of $L$ in $L^{\vee}$ via $c\mapsto b(c,-)$. The quotient $L^{\vee}/L$ is called the discriminant group of $L$ and is denoted by $A_L$. The minimal number of generators of the group $A_L$ is called length of $L$ and is indicated by $l(L)$.\\
The bilinear form $b$ on $L$ induces a bilinear form $L^{\vee}\times L^{\vee}\ra \Q$ on $L^{\vee}$ and thus a bilinear form $b_{A_L}:A_L\times A_L \rightarrow \Q/\Z$ defined over $A_L=L^{\vee}/L$. The bilinear form $b_{A_L}$ is called discriminant form of $L$. \\
A lattice $L$ is called 2-elementary if $A_L\simeq (\Z/2\Z)^a$. A 2-elementary lattice is uniquely determined by three invariants: $r$, which is the rank of $L$, $a$ which is the length of $L$ and $\delta$ which takes values either $0$, if the discriminant form takes value in $\Z$, or $1$ otherwise (cf.\ \cite{Nikulin non symplectic}).\\
We denote by $U$ the unique (up to isometries) even unimodular hyperbolic lattice of rank 2. It is associated to the matrix
$\left(\begin{array}{cc}
0&1\\
1&0\\
\end{array}\right).
$ We denote by $E_8$ the unique even unimodular positive
definite lattice of rank eight (it is associated to the Dynkin diagram $E_8$).\\
Let $X$ be a K3 surface. Then $H^2(X,\Z)$ with the cup product is an even unimodular lattice with signature $(3,19)$. Up to isometries there exists only one lattice with these properties, which is $U\oplus U\oplus U\oplus E_8(-1)\oplus E_8(-1)$. This lattice is denoted by $\Lambda_{K3}$.\\
The N\'eron--Severi group of $X$, defined as $\NS(X):=H^{1,1}(X)\cap H^2(X,\Z)$, is a sublattice of $\Lambda_{K3}$. Since we consider algebraic K3 surfaces, the N\'eron--Severi group is an hyperbolic lattice, by the Hodge index theorem. We denote by $\rho(X)$ the Picard number of $X$, i.e.\ the rank of its N\'eron--Severi group. The transcendental lattice $\T_X$ is the orthogonal of $\NS(X)$ in $H^2(X,\Z)$, thus its signature is $(2,20-\rho(X))$
  
\subsection{Elliptic fibrations}\label{subsec: elliptic fibration}
We recall some known facts on elliptic fibrations, which are very important in the following. A good and recent reference for these results is \cite{SS}.

\begin{definition} Let $S$ be a surface and $D$ a curve. An elliptic fibration $\mathcal{E}:S\ra D$ on the surface $S$ is a surjective morphism such that the generic fiber is a smooth curve of genus one and such that a section $s:D\ra S$ is given. We call this section the zero section.
\end{definition}

Every elliptic fibration can be regarded as an elliptic curve over the function field of the basis. We will assume $D\simeq \mathbb{P}^1$ with homogeneous coordinate $(\tau:\sigma)$. Under this assumption each elliptic fibration admits a minimal Weierstrass equation of the form
$$y^2=x^3+A(\tau:\sigma)x+B(\tau:\sigma),\ \ \ A(\tau:\sigma),\ B(\tau:\sigma)\in \C[\tau:\sigma]_{hom},\ \ \deg A(\tau:\sigma)=4m,\ \deg B(\tau:\sigma)=6m
$$
for a certain $m\in \N_{>0}$ and where there exists no polynomials $C(\tau:\sigma)$ such that $C(\tau:\sigma)^4|A(\tau:\sigma)$ and $C(\tau:\sigma)^6|B(\tau:\sigma)$. If $S$ is a K3 surface, then $m=2$, i.e., $\deg(A(\tau:\sigma))=8$, $\deg(B(\tau:\sigma))=12$.
We use the notation $A(\tau)$ and $B(\tau)$ to indicate the polynomials $A(\tau:1)$ and $B(\tau:1)$.\\
There are a finite number of singular fibers, which are the fibers over the points $\overline{\tau}\in\mathbb{P}^1$ where $\Delta(\overline\tau)=-16(4A^3(\overline\tau)+27B^2(\overline\tau))$ is zero. For each singular fiber $F_{\overline{\tau}}$ of the fibration we will denote by $\delta(F_{\overline{\tau}})$ the multiplicity of zero of $\Delta$ in $\overline{\tau}$.
The possible singular fibers on an elliptic fibration are described by Kodaira: their components are rational curves. One of these components meets the zero section, the other ones meet each other with  the configuration of a Dynkin diagram.  For each singular fibers $F$ we denote by $r(F)$ the number of components of $F$ and by $d(F)$ the discriminant of the lattice associated to the Dynkin diagram of $F$. A simple component of a fiber is a component with multiplicity 1.\\
In the following table we describe the singular fibers of an elliptic fibration. We denote by $\Theta_0$ the component of a fiber meeting the zero section. The first column contains the name of the reducible fiber according to Kodaira classification, the second the Dynkin diagram associated to the fiber, the third column contains the description of the intersections among the components of the fibers, the last column components which are simple.
\begin{eqnarray*}
\tiny{
\begin{array}{|c|c|c|c|}
\hline
II&&\mbox{one cuspidal rational curve }\Theta_0&\Theta_0\\
\hline
III&A_1&\mbox{two tangent rational curves }\Theta_0,\ \Theta_1&\Theta_0,\Theta_1\\
\hline
IV&A_2&\mbox{three lines meeting in a point }\Theta_0,\ \Theta_1,\ \Theta_2&\Theta_i, i=0,1,2\\
\hline
I_1&&\mbox{one nodal rational curve }\Theta_0&\Theta_0\\
\hline
I_2&A_1&\tiny{\begin{array}{ccc}\Theta_0&=&\Theta_1\end{array}}&\Theta_0,\Theta_1\\
\hline
I_n&A_{n-1}&\tiny{\begin{array}{cccccc}\Theta_0&-&\Theta_1&-&\ldots&\Theta_i\\
|&&&&&|\\
\Theta_{n-1}&-&\Theta_{n-2}&-&\ldots&\Theta_{i+1}\end{array}}&\begin{array}{c}\Theta_i, i=\\
0,\ldots, n-1\end{array}\\
\hline
I_n^*&D_{n+4}&\tiny{\begin{array}{cccccccccccc}
\Theta_0&&&&&&&&\Theta_{n+3}\\
&\diagdown&&&&&&\diagup&\\
&&\Theta_2&\ldots&\Theta_i -\Theta_{i+1} &\ldots&\Theta_{n+2}&\\
&\diagup&&&&&&\diagdown&\\
\Theta_1&&&&&&&&\Theta_{n+4}
\end{array}}&\tiny{\begin{array}{c}\Theta_i,\\ i=0,1,\\n+3,n+4\end{array}}\\
\hline
IV^*&E_6&\tiny{\begin{array}{ccccccccccccccccc}
\Theta_{0}&-&\Theta_{1}&-&\Theta_{2}&-&\Theta_{3}&-&\Theta_{4}\\
    & &     & &\!\mid & &  & &  & &   & \\
     & &      & & \Theta_5   & &   & &   & &   \\
     & &      & &\mid & &   & &   & &   &\\
     & &      & & \Theta_6   & &   & &   & &   \\
\end{array}}&\begin{array}{c}\Theta_i,\\ i=0,4,6\end{array}\\
\hline
III^*&E_7&\tiny{\begin{array}{cccccccccccccc}
\Theta_{0}&-&\Theta_{2}&-&\Theta_{3}&-&\Theta_{4}&-&\Theta_{5}&-&\Theta_{6}&-&\Theta_{7}\\
      & &     & &&&\mid & &   & &   & &   & \\
      & &      & &&& \Theta_1   & &   & &   & &   & \\
\end{array}}&\begin{array}{c}\Theta_i,\\ i=0,7\end{array}\\
\hline
II^*&E_8&\tiny{\begin{array}{ccccccccccccccc}
\Theta_0&-&\Theta_1&-&\Theta_2&-&\Theta_3&-&\Theta_4&-&\Theta_5&-&\Theta_6&-&\Theta_7\\
      & &      & && &   & &   & &   \mid & &&&\\
      & &      & &    & &   & &   & &\Theta_8   & &&&
\end{array}}&\Theta_0\\
\hline
\end{array}}
\end{eqnarray*}

The set of the sections of an elliptic fibration $\mathcal{E}$ form a group, $\MW(\mathcal{E})$ (the Mordell--Weil group), with the group law which is induced by the one on the smooth fibers. Let $\mathcal{E}$ be an elliptic fibration on the surface $S$ and let $\Red=\{v\in \mathbb{P}^1\,|\, F_v\mbox{ is reducible}\}$. We recall the Shioda--Tate formula: $\rho(S)=\rk(\MW(\mathcal{E}))+2+\sum_{v\in \Red} (r(F_v)-1)$.\\
The trivial lattice $\Tr_S$ of an elliptic fibration on a surface $S$ is the lattice generated by the class of the fiber, the class of the zero section and the classes of the irreducible components of the reducible fibers which do not intersect the zero section. Its rank is $\rho(S)-\rk(\MW(\mathcal{E}))$.

One can define a quadratic form on $\MW(\mathcal{E})$, which takes value in $\Q$, by the height formula: let $P$ be a section of an elliptic fibration $\mathcal{E}$, then we set: $ \langle P, P \rangle=2\chi+2(P\cdot s)-\sum_{v\in \Red}{\rm contr}_v(P)$, where ${\rm contr}_v(P)$ depends on the component of the reducible fibers which meets $P$. The main property of this pairing is that $\langle P, P\rangle=0$ if and only if $P$ is a torsion section. If $P$ meets the component $\Theta_0$, then  ${\rm contr}_v(P) =0$, otherwise, if $P$ meets the component $\Theta_i$, $i\neq 0$ the value of ${\rm contr}_{v(P)}$ is listed in the following table:
$$\begin{array}{c||c|c|c|c}
\mbox{fiber}&I_n, n\geq 2\ (III,IV)&I_n^*&IV^*&III^*\\
\hline
\mbox{Dynkin diagram}&A_{n-1}\ (A_1, A_2)&D_{n+4}&E_6&E_7\\
\hline
{\rm contr}_v(P)&i(n-i)/n&\left\{\begin{array}{ll}1&\mbox{ if }i=1\\1+n/4&\mbox{ if }i=n+3,n+4\end{array}\right.&4/3&3/2\\
\end{array}
$$
The Mordell--Weil group equipped with this bilinear form is called Mordell--Weil lattice and is denoted by $\MWL$.

\subsection{Symplectic involutions on K3 surfaces}\label{subsec: symplectic involutions}
Let $X$ be a K3 surfaces admitting a symplectic involution $\sigma$. Then the anti-invariant lattice $(H^2(X,\Z)^\sigma)^\perp$ is isometric to $E_8(-2)$ and it is primitively embedded in $\NS(X)$. Since $E_8(-2)$ is a negative definite lattice and the N\'eron--Severi group of a projective K3 surface is an hyperbolic lattice, the lattice orthogonal to $E_8(-2)$ in $\NS(X)$ has rank at least 1. Thus the Picard number of a K3 surface with a symplectic involution is at least 9 and the moduli of K3 surfaces with a symplectic involution are 11.\\
The fixed locus $Fix_\sigma(X):=\{x\in X\ |\ \sigma(x)=x\}$ consists of eight isolated points, thus the quotient $X/\sigma$ is a singular surface with 8 singularities of type $A_1$. The desingularization $Y:=\widetilde{X/\sigma}$ is a K3 surface, in general not isomorphic to $X$. The 8 rational curves on $Y$ introduced by the blow up which resolve the singularities of $X/\sigma$ define classes in $\NS(Y)$. The minimal primitive sub--lattice of $\NS(Y)$ containing these classes is the unique even 2-elementary lattice with signature $(0,8)$ and discriminant group $(\Z/2\Z)^6$. It is called Nikulin lattice and is denoted by $N$.\\
A symplectic involution can be realized in several ways, according to the particular projective model considered for a K3 surface. In \cite{vGS} symplectic involutions acting on K3 surfaces realized as quartics in $\mathbb{P}^3$, complete intersections in $\mathbb{P}^4$ and $\mathbb{P}^5$, double covers of $\mathbb{P}^2$ are described. Another way to construct symplectic involutions is to consider K3 surfaces with an elliptic fibration: the translation $\sigma_t$ by a section $t$ on an elliptic K3 surface is a symplectic automorphism, and if the section has order $n$ in the Mordell--Weil group then $\sigma_t$ is a symplectic automorphisms of order $n$. The case $n=2$ is presented in \cite{vGS}, the cases $n>2$ are analyzed in \cite{GS1}, \cite{GS2}.
As in \cite{CD:vGS involutions} we give the following definition:
\begin{definition}The symplectic involutions on an elliptic K3 surface $X$ which are translations by a section of order 2 of an elliptic fibration on $X$ are called van Geemen--Sarti involutions.\end{definition}
The family of K3 surfaces admitting a van Geemen--Sarti involution is presented in \cite[Proposition 4.2]{vGS}: it is the family of the $(U\oplus N)$-polarized K3 surfaces and it has dimension 10.\\
Let $\mathcal{E}$ be an elliptic fibration on a K3 surface $X$ with a 2--torsion section $t$. The van Geemen--Sarti involution $\sigma_t$ is trivial on the base of the fibration and acts only on the fiber. Thus it sends a fiber to itself. This implies that the elliptic fibration is preserved by the van Geemen--Sarti involution and that the elliptic fibration $\mathcal{E}$ induces an elliptic fibration, denoted by $\mathcal{E}/\sigma_t$, over $\widetilde{X/\sigma_t}$.\\
It is well known (cf. \cite[Pag. 79]{ST}) that an elliptic curve with a 2-torsion rational point admits an equation of the form $y^2=x(x^2+ax+b)$ and that the elliptic curve which is its quotient by the translation by the point of order 2 has equation $y^2=x(x^2-2ax+a^2-4b)$. In particular let $\mathcal{E}$ be an elliptic fibration on a K3 surface with a 2-torsion section $t$, then $\mathcal{E}$ and $\mathcal{E}/\sigma_t$ admit the following equations:\begin{equation}\mathcal{E}:\ y^2=x (x^2+a(\tau)x+b(\tau)),\ \ \mathcal{E}/\sigma_t:\ y^2=x(x^2-2a(\tau)x+a(\tau)^2-4b(\tau)) \end{equation} where $a(\tau:\sigma)$ and $b(\tau:\sigma)$ are polynomials of degree 4 and 8 respectively.   
\subsection{Non--symplectic involutions on K3 surfaces}\label{subsec: non symplectic involutions}
In \cite{Nikulin non symplectic} Nikulin classified the non--symplectic involutions on K3 surfaces. 
If $\iota$ is a non--symplectic involution on a K3 surface $X$, the lattice $H^2(X,\Z)^\iota$ is hyperbolic and 2-elementary, thus it is uniquely determined by the invariants $(r,a,\delta)$. We will denote by $N_{(r,a,\delta)}$ the unique (up to isometries) even hyperbolic 2-elementary lattice with invariant $(r,a,\delta)$.\\
If $\iota$ is a non--symplectic involution of a K3 surface, the fixed locus $Fix_\iota(X)$ is smooth and is one of the following: 1) empty; 2) made up of two elliptic curves; 3) made up of $k$ rational curves and one curve of genus $g\geq 0$. 
The fixed locus is associated to the invariant lattice $H^2(X,\Z)^\iota$: the lattice $H^2(X,\Z)^\iota$ determines the fixed locus uniquely; given a fixed locus it determines uniquely the values $r,a$ and thus a particular fixed locus is associated at most to two lattices (the ones with the same values of $r$ and $a$ but with a different value of $\delta$). More precisely:\\
1) if $Fix_\iota(X)$ is empty then the lattice $H^2(X,\Z)^\iota$ is $N_{(10,10,0)}$;\\
2) if $Fix_\iota(X)$ is the disjoint union of two elliptic curves then $H^2(X,\Z)^\iota$ is $N_{(10,8,0)}$;\\
3) if $Fix_\iota(X)$ is the disjoint union of $k$ rational curves and a curve of genus $g$ then $H^2(X,\Z)^\iota$ is $N_{(r,a,\delta)}$ where \begin{equation}\label{eq: g,k,r,a}g=\frac{22-r-a}{2},\ \ \ k=\frac{r-a}{2}.\end{equation}
If $X$ admits a non--symplectic involution $\iota$, then $H^2(X,\Z)^\iota$ is primitively embedded in $\NS(X)$. The K3 surfaces admitting a non--symplectic involution with a given fixed locus appear in families: they are the $N_{(r,a,\delta)}$-polarized K3 surfaces (where $N_{(r,a,\delta)}$ is the lattice associated to the chosen fixed locus). The general member in this family has N\'eron--Severi group isometric to $N_{(r,a,\delta)}$ and will be denoted by $X_{(r,a,\delta)}$.

In Section \ref{subsec: symplectic involutions} we saw an example of symplectic involution (the van Geemen--Sarti involution) preserving an elliptic fibration on a K3 surface. Here we consider three examples of non--symplectic involutions preserving an elliptic fibration on a K3 surface.
\begin{example}\label{ex: involution on the basis}{\rm Let $\varphi:S\ra\mathbb{P}^1$ be an elliptic fibration on a rational surface $S$ with Weierstrass equation $y^2=x^3+A(\tau)x+B(\tau)$. Let $\alpha:\mathbb{P}^1\ra \mathbb{P}^1$ be a double cover. Up to a projectivity, we can assume that $\alpha:\tau\mapsto\tau^2$. We assume that the fibers of $\varphi$ over $\tau=0$ and $\tau=\infty$ are reduced.  After pulling back $\varphi$ via $\alpha$, we obtain an elliptic fibration on a K3 surface $X$ with Weierstrass equation $y^2=x^3+A(\tau^2)x+B(\tau^2)$. It is clear that there exists an involution $i:(x,y,\tau)\mapsto (x,y,-\tau)$ on $X$ (we will refer to this involution as "involution of type a)"). It acts only on the basis of the fibration and is non--symplectic. It sends a fiber to another one and thus the class of the fiber is preserved by $i$. There are two fibers which are fixed, the ones over $\tau=0, \tau=\infty$. Since $i$ does not act on the fibers, all the sections of the fibration are preserved by $i$. We observe that $S=X/i$, by construction.\erem}\end{example}

\begin{example}\label{ex: hyperelliptic involution}{\rm An elliptic curve is a $2:1$ cover of a rational curve branched over 4 points, thus on each elliptic curve there exists an involution which is the cover involution and which fixes the four points $p$ such that $2p=0$. This involution extends to an involution, $i$, of any elliptic fibration  (we will refer to this involution as "the hyperelliptic involution" or as "involution of type b)"). If one consider the equation $y^2=x^3+A(\tau)x+B(\tau)$ of the elliptic fibration, then $i$ is $(x,y,\tau)\mapsto (x,-y,\tau)$.  The involution $i$ fixes the basis of the fibration, thus sends a fiber to itself. The zero section is fixed by $i$. Moreover, the 2-torsion points are fixed on every fiber. Generically the points of order 2 lie on an irreducible curve, which is a trisection, fixed by $i$. However, for certain elliptic fibrations the trisection splits in a section and a bisection or into three sections.\erem}\end{example}

\begin{example}\label{ex: non symplectic involution Q-P}{\rm Let $\mathcal{E}$ be an elliptic fibration on a K3 surface. For each section $P\in \MW(\mathcal{E})$ there exists a non--symplectic involution $\beta_P$.  Indeed a rational point $P$ on an elliptic curve $E$ defines the involution $\beta_P:Q\mapsto P-Q$, $Q\in E$. In particular, such an involution extends to an involution of the fibration if $P$ is a section. This involution is non--symplectic because it is the composition of the hyperelliptic involution (which is non--symplectic) and a translation (which is symplectic). The involution $\beta_P$ preserves the fibers. If $P$ is the zero of the Mordell--Weil group, then $\beta_P$ is the hyperelliptic involution.  
}\end{example}

\section{Double cover of $\mathbb{P}^2$ }\label{section:double cover}

Let us consider double covers of $\mathbb{P}^2$ branched over a sextic $\mathcal{C}_6$. This gives a (singular) model of a K3 surface. The cover involution is a non--symplectic involution with fixed locus the ramification locus. In this way we obtain a pair $(X,i)$ where $X$ is a K3 surface and $i$ a non--symplectic involution on it and thus generically a pair $(X_{(r,a,\delta)},\iota)$, with the notation of the previous section. The aim of this Section is to describe more precisely the relation among the geometric properties of $\mathcal{C}_6$ and the numbers $(r,a,\delta)$ which describe the K3 surface $X$. The main results are stated in Corollary \ref{cor: double cover 1} and Corollary \ref{cor: double cover 18}.\\

The following Proposition resumes some well known results on K3 surfaces which are double covers of $\mathbb{P}^2$ branched over a (possibly singular and reducible) sextic.
\begin{proposition}\label{prop: double cover of P2}
 Let $\mathcal{C}_6$ be a sextic with singular points which are either ordinary double points or ordinary triple points. Assume that $\mathcal{C}_6$ has $\alpha$ components $b_1,\ldots, b_\alpha$, $\gamma$ double points $P_1,\ldots, P_\gamma$ and $\nu$ triple points $Q_1,\ldots, Q_\nu$. Let $X$ be the K3 surface obtained as desingularization of the double cover of $\mathbb{P}^2$ branched over $\mathcal{C}_6$. Let $\iota$ be the involution induced on $X$ by the cover involution. Then the Picard number of $X$ is $\rho(X)\geq 1+\gamma+4\nu$ and $\iota$ fixes $s:=\alpha+\nu$ curves. The involution $\iota$ acts as the identity on $\NS(X)$ if and only if $\rho(X)=1+\gamma+4\nu$.\end{proposition}
\proof We first consider the case $\nu=0$. The double cover $S$ of $\mathbb{P}^2$ branched along $\mathcal{C}_6$ is singular over the singular points of $\mathcal{C}_6$. To resolve these singularities one blows up the double cover. Equivalently one can blow up $\mathbb{P}^2$ in the singular points of $\mathcal{C}_6$ and then consider the double cover of this blow up. More precisely there is the following commutative diagram: 
\begin{equation*}
    \xymatrix{
      S \ar@{->}[d]_{2:1} & X \ar[l]\ar@{->}[d]^{2:1}_\pi  \\
      \mathbb P^2 & \widetilde{\mathbb{P}^2}   \ar[l]_\beta }
\end{equation*} 
where $\beta:\widetilde{\mathbb{P}^2}\ra \mathbb{P}^2$ is the blow up of $\mathbb{P}^2$ in $P_i$, $i=1,\ldots, \gamma$. Let $\widetilde{b}_i\subset \widetilde{\PP^2}$ be the strict transform of $b_i$. Let $h$ be the class of a line in $\mathbb{P}^2$ and $e_i$ be the exceptional divisor over $P_i$, $i=1,\ldots ,\gamma$. One finds $\beta^{-1}(\mathcal{C}_6)=\bigcup_{i=1}^\alpha\widetilde{b}_i\bigcup_{j=1}^\gamma e_j$ and the divisors $e_j$ appear with multiplicity 2 because $P_j$ are double points. Thus the branch locus of the double cover $\pi: X\ra \widetilde{\mathbb{P}}^2$ consists of the disjoint union of $\widetilde{b}_i$, $i=1,\ldots, \alpha$. The classes $h$, $e_i$, $i=1,\ldots, \gamma$ generate $\NS(\widetilde{\mathbb{P}^2})\otimes \Q$.
\newline Since $\left(\NS(X)\otimes \Q\right)^\iota=\pi^*(\NS(\widetilde{\mathbb{P}^2})\otimes \Q)$\mbox{ one obtains }$\rk\left((\NS(X)\otimes \Q)^\iota\right)=\rho(\widetilde{\mathbb{P}^2})=1+\gamma.$  It is clear that $\left(\NS(X)\otimes \Q\right)^\iota$ is embedded in $\NS(X)\otimes \Q$, and thus $\rho(X)\geq 1+\gamma$. Moreover, $\iota$ acts as the identity on $\NS(X)$ if and only if $\left(\NS(X)\otimes \Q\right)^\iota=\NS(X)\otimes \Q$, thus if and only if $\rho(X)=1+\gamma$.\\
Let us now consider the case $\nu>0$. Let $Q_j$ be a triple point on $\mathcal{C}_6$. The blow up of $Q_j$ introduces an exceptional divisor $f_j$ with multiplicity 3. Since it has an odd multiplicity, $f_j$ is a class in the branch locus of the double cover $\pi:X\ra \widetilde{\mathbb{P}^2}$. The curve $f_j$ meets the strict transform of $\mathcal{C}_6$ in 3 points (because $Q$ is a triple point). To obtain a smooth branch locus one has to blow up these three points, introducing the classes $f_{j,1}$, $f_{j,2}$, $f_{j,3}$. Each of them has multiplicity 2 and thus is not contained in the branch locus. Hence over the triple points we found four curves with the following intersection properties: $f_j^2=-4$, $f_{j,h}^2=-1$, $f_jf_{j,h}=1$, $f_{j,i}f_{j,h}=0$, $i,h=1,2,3$, $i\neq h$. Thus  a triple point introduces 4 exceptional divisors in $\NS(\widetilde{\mathbb{P}^2})$ (the curves $f_j$,$f_{j,1}$,$f_{j,2}$,$f_{j,3}$) and one curve ($f_j$) in the fixed locus. One concludes the proof as in case $\nu=0$.\eprf
\begin{definition} We will say that the components $b_i$, $i=1,\ldots, \alpha$ are in general position if $\rho(X)=1+\gamma+4\nu$\end{definition}
If the components $b_i$ are in general position then $\NS(X)\otimes \Q=\pi^*(\NS(\widetilde{\mathbb{P}^2}))\otimes \Q$, thus $\pi^*(h)$, $\pi^*(e_i)$, $i=1,\ldots \gamma$, $	\pi^*(f_j)$, $j=1,\ldots \nu$, $\pi^*(f_{j,h})$ $j=1,\ldots,\nu$, $h=1,2,3$ form a $\Q$-basis of $\NS(X)$.\\ 
The following Corollary follows immediately by Proposition \ref{prop: double cover of P2} and by the classification of the non--symplectic involutions on K3 surfaces (see Section \ref{subsec: non symplectic involutions}).
\begin{corollary}\label{cor: double cover 1} Let $\mathcal{C}_6$ and $b_i$ be as in Proposition \ref{prop: double cover of P2} and let $b_i$ be in general position. Let $g_i=g(b_i)$ be the genus of $b_i$ and $g_i\geq g_{i+1}$. Then: \begin{itemize}\item[i)] if $g_2\neq 0$ then $g_1=g_2=1$,  $\alpha=2$, $\gamma=9$, $\nu=0$;\item[ii)] if $g_2=0$ then $g=g_1$, $k=\alpha+\nu-1$, $r=1+\gamma+4\nu$, $a=21-\gamma-4\nu-2g_1$, and thus  $\gamma=\alpha-3\nu +9-g_1$ {\rm (}with the notation of \eqref{eq: g,k,r,a}{\rm )}.\end{itemize} \end{corollary}
\begin{ass}\label{ass}From now on we always assume that the curves $b_i$ are in general position and rational. \end{ass}

\begin{notation}{\rm We number the components of the branching sextic; if a double point is the intersection between the $i$-th and the $j$-th components, it is denoted by $P_{i,j}$ (resp. $P_{i,j}^k$, $k=1,\ldots, h$ if the intersection consists of $h$ distinct points), the exceptional curve in $\widetilde{\mathbb{P}^2}$ of the point $P_{i,j}$ (resp. $P_{i,j}^k$) will be denoted by $e_{P_{i,j}}$  (resp. $e_{P_{i,j}^k}$) and its pull-back on the K3 surface $X_{(r,a,\delta)}$ is denoted by $E_{P_{i,j}}$ (resp. $E_{P_{i,j}^k}$). If $A$ is a triple point, intersection of the $i$-th, the $j$-th and the $k$-th components, then we will call $e_{A}$, $e_{A,i}$, $e_{A,j}$ and $e_{A,k}$ the exceptional curves over $A$. The intersection properties of these curves are $e_{A}e_{A,i}=e_Ae_{A,j}=e_Ae_{A,k}=1$, $e_{A,i}e_{A,j}=e_{A,i}e_{A,k}=e_{A,j}e_{A,k}=0$, $\{i,j,k\}=\{1,2,3\}$. We denote by $E_{A,i}:=\pi^*(e_{A,i})$, $i=1,2,3$. The pull back on $X_{(r,a,\delta)}$ of the strict transform of components $b_i$ (resp. $l_i$, $c_i$ when they are lines or conics) and the pull back on $X_{(r,a,\delta)}$ of $e_A$ are curves with multiplicity 2 (because of $\beta^*(b_i)$, $\beta^*(l_i)$, $\beta^*(c_i)$ and $e_A$ are in the branch locus). We denote by $\mathcal{B}_i$ (resp. $\mathcal{L}_i$, $\mathcal{C}_i$, $E_A$) the curves on $X$ such that $\pi^*(\beta^*(b_i))=2\mathcal{B}_i$ (resp. $\pi^*(\beta^*(l_i))=2\mathcal{L}_i$, $\pi^*(\beta^*(c_i))=2\mathcal{C}_i$, $\pi^*(e_A)=2E_A$).\erem }\end{notation}

\begin{corollary}\label{cor: double cover 18} Under the Assumption \ref{ass} the double cover of $\mathbb{P}^2$ branched over $\mathcal{C}_6$ is $X_{(r,a,\delta)}$ with $$a=22-r,\ \ \ r=1+\gamma+4\nu,$$ and the cover involution is $\iota$. The fixed locus of $\iota$ on $X_{(r,a,\delta)}$ consists of  $k+1=\alpha+\nu$ rational curves (this implies $\gamma=\alpha-3\nu+9$).\\
Moreover, if $\delta=0$ then $(r,a)=(18,4)$ and $X_{(18,4,0)}$ is the Kummer surface $Km(E\times E')$ where $E$ and $E'$ are two non isogenous elliptic curves.\\
The double cover of $\mathbb{P}^2$ branched along six lines with two triple points $A$, $B$ is either $X_{(18,4,1)}$ or $X_{(18,4,0)}$. In the first case the line through $A$ and $B$ is in the branch locus, in the second case the line through $A$ and $B$ is not in the branch locus.\end{corollary}
\bprf The conditions on $a$, $k$, $r$ and $\gamma$ depend on \eqref{eq: g,k,r,a} and on Corollary \ref{cor: double cover 1}. \\
From the classification of non--symplectic involutions, it follows that the surface $X_{(r,22-r,0)}$ has $r=18$. In particular the K3 surface $X_{(18,4,0)}$ is the Kummer of the product of two non isogenous elliptic curves. Indeed the transcendental lattice of the Kummer surface $Km(E\times E')$ is isometric to $U(2)\oplus U(2)$ (cf.\ \cite{O}). This lattice is a 2-elementary lattice with invariant $a=4$, $\delta=0$ and thus the N\'eron--Severi group of $Km(E\times E')$ is a 2-elementary lattice with invariants $(18,4,0)$.\\
We first consider the sextic which consists of the 6 lines $l_i$, $i=1,\ldots, 6$, with the 2 triple points $A:=l_1\cap l_2\cap l_3$, $B:=l_4\cap l_5\cap l_6$ and 9 double points $P_{1,4}$, $P_{1,5}$, $P_{1,6}$, $P_{2,4}$, $P_{2,5}$, $P_{2,6}$, $P_{3,4}$, $P_{3,5}$, $P_{3,6}$. In this case the line through $A$ and $B$ is not in the branch locus. The invariants associated to this configuration of lines are $\alpha=6$, $\nu=2$, $\gamma=9$. By Corollary \ref{cor: double cover 1} one obtains $r=18$, $a=4$. Hence the double cover of $\mathbb{P}^2$ branched over this sextic is either $X_{(18,4,0)}$ or $X_{(18,4,1)}$. The classes $\mathcal{L}_i$, $i=1,\ldots 6$, $E_A$, $E_B$, $E_{A,i}$, $i=1,2,3$, $E_{B,j}$, $j=4,5,6$, $E_{P_{h,k}}$, $(h,k)=(1,4),(1,5), (1,6), (2,4), (2,5), (2,6), (3,4), (3,5), (3,6)$ form a particular lattice, called double Kummer (cf. \cite{O}), which is typical of the Kummer surface of the product of two elliptic curves (in particular the last 16 classes form a Kummer lattice). Thus, the double cover of $\mathbb{P}^2$ branched over this configuration of 6 lines is $Km(E\times E')$ and hence $X_{(18,4,0)}$.\\ Now we consider the sextic which consists of 6 lines $l_i$, $i=1,\ldots, 6$, with 2 triple points $A:=l_1\cap l_2\cap l_3$, $B:=l_1\cap l_5\cap l_6$ and 9 double points $P_{1,4}$, $P_{2,4}$, $P_{2,5}$, $P_{2,6}$, $P_{3,4}$, $P_{3,5}$, $P_{3,6}$, $P_{4,5}$, $P_{4,6}$. The line through $A$ and $B$ is $l_1$ and is contained in the branch locus. The invariants associated to this configuration of lines are $\alpha=6$, $\nu=2$, $\gamma=9$. By Corollary \ref{cor: double cover 1} one obtains $r=18$, $a=4$. Thus, the double cover of $\mathbb{P}^2$ branched over this sextic is either $X_{(18,4,0)}$ or $X_{(18,4,1)}$.  The class $F:=\mathcal{L}_1+E_{P_{1,4}}+\mathcal{L}_4+E_{P_{4,6}}+\mathcal{L}_6+E_{B,6}+E_{B}+E_{B,5}+\mathcal{L}_5+E_{P_{2,5}}+\mathcal{L}_2+E_{A,2}+E_A+E_{A,1}$ gives an elliptic fibration with section $s:=E_{A,3}$ and with a fiber of type $I_{14}$. By \cite{O}, there is not an elliptic fibration with such a fiber on $X_{(18,4,0)}$. Hence this configuration of lines on $\mathbb{P}^2$ corresponds to the K3 surface $X_{(18,4,1)}$.\eprf

In Corollary \ref{cor: double cover 1} we proved that $g_i$, $\alpha$, $\gamma$, $\nu$ determine $r$ and $a$. However there are several different choices for $\alpha$, $\gamma$, $\nu$, and thus for the sextic $\mathcal{C}_6$, which are related to the same choice of the invariant $(r,a)$, as shown in the following Example. 

\begin{example}\label{ex: branch locus X16,6,1}{\rm Let $X$ be the K3 surface which is the desingularization of the double cover of $\mathbb{P}^2$ branched along six lines in general position. Then $\alpha=6$, $\gamma=15$, $\nu=0$. These numbers determine the invariant $(r,a)=(16,6)$, thus $X$ is $X_{(16,6,1)}$. The surface $X$ is associated also to other sextics. The possible choices for the sextic $\mathcal{C}_6$ are: 
\begin{enumerate}[\it i)]
 \item 1 rational quartic, 2 lines and 3 triple points: $\alpha=3$, $\gamma=3$ $\nu=3$;
 \item 1 rational cubic, 1 conic, 1 line and 3 triple points: $\alpha=3$, $\gamma=3$, $\nu=3$;
 \item 1 rational cubic, 3 lines and 2 triple points: $\alpha=4$, $\gamma=7$, $\nu=2$;
 \item 3 conics and 3 triple points: $\alpha=3$, $\gamma=3$, $\nu=3$;
 \item 2 conics, 2 lines and 2 triple points: $\alpha=4$, $\gamma=7$, $\nu=2$;
 \item 1 conic, 4 lines and 1 triple point: $\alpha=5$, $\gamma=11$, $\nu=1$;
 \item 6 lines: $\alpha=6$, $\gamma=15$, $\nu=0$.
\end{enumerate}

The surface $X_{(r,a,\delta)}/\iota$ does not depend  on the choice of $\mathcal{C}_6$ up to birational transformation. Thus there is a birational map of $\mathbb{P}^2$ to itself which transforms one of the sextics of the previous list to another one. To transform the sextic  $vii)$ to the sextic $vi)$ (resp. $v)$, $iv)$) one considers the Cremona transformation centered in three double points of $\mathcal{C}_6$ contracting 1 (resp. 2, 3) lines of the branch locus. To transform the sextic $vi)$ to the sextic $iii)$ (resp. $ii)$, $i)$) one considers the Cremona transformation centered in three double points of $\mathcal{C}_6$, one (resp. one, zero) of them lying on the conic, and contracting 1 line (resp. 2 lines, 2 lines). \erem}\end{example}

In case $r=18$, different choices of $\mathcal{C}_6$ associated to the same $\alpha$, $\gamma$, $\nu$ are not necessary birational: in Corollary \ref{cor: double cover 18} we show that different configuration of 6 lines with two triple points are associated to different K3 surfaces ($X_{(18,4,0)}$ and $X_{(18,4,1)}$). We analyze the similar situation for other choices of $\mathcal{C}_6$:
\begin{example}\label{ex: branch locus r=18} {\rm Let us fix $(r,a)=(18,4)$. The possible choices for the components of $\mathcal{C}_6$ are:\\ 1) 2 conics, 2 lines and 4 triple points: $\alpha=4$, $\gamma=1$, $\nu=4$; \\ 2) 1 conic, 4 lines and 3 triple points: $\alpha=5$, $\gamma=5$, $\nu=3$;\\ 3) 6 lines and 2 triple points: $\alpha=6$, $\gamma=9$, $\nu=2$.\\ We denote by 3a) the case associated to $X_{(18,4,1)}$ (the line through the triple points is in the branch locus), and by 3b) the other one. Let us denote by $l_i$, $i=1,\ldots, 6$ the 6 lines of the branch locus. In case 3a) we put $A=l_1\cap l_2\cap l_3$ and $B=l_3\cap l_4\cap l_5$. The Cremona  transformation centered in $A$, $B$ and $P_{2,4}$ gives a sextic of type 2) such that every triple point is the intersection of two lines and the conic and every line passes through a triple point.
In case 3b) we put $A=l_1\cap l_2\cap l_3$, $B=l_4\cap l_5\cap l_6$. The Cremona  transformation centered in $A$, $P_{1,4}$ and $P_{2,5}$ gives a sextic of type 2) such that every triple point is the intersection of two line and the conic and there is a line which does not pass through a triple point. The Cremona transformation centered in $B$, $P_{2,4}$, $P_{2,6}$ gives a sextic of type 1). \erem}\end{example}

\section{Elliptic fibration on $X_{(r,22-r,\delta)}$: possibilities}\label{section: possibilities}
In this section we use the results given in \cite{O} and \cite{K} and conditions coming from the lattice properties of the N\'eron--Severi group in order to find a list of the admissible elliptic fibrations on the surfaces $X_{(r,22-r,\delta)}$. This is done in Propositions \ref{prop: classification non finite MW} and \ref{prop: classification finite MW}. 

\begin{definition}\label{def: special curve} A rational curve in $X$ which is in the branch locus (i.e. which is among $\widetilde{b_i}$ or $f_j$ with the notations of the previous section) is called special curve; a rational curve which is not in the branch locus is called ordinary rational curve.\end{definition}
\begin{lemma}\label{lemma:DB=2}{\rm (\cite[Lemma 5.11]{K})} Let $V$ be the branch locus, i.e.\ $V=\cup_i \widetilde{b_i}\cup_j f_j$. Let $D$ be an ordinary curve. Then $D\cdot V=2$.\end{lemma}
\begin{lemma}\label{lemma:D1D2 even}{\rm (\cite[Lemma 1.6]{O}\cite[Lemma 5.12]{K})} Let $D_1$ and $D_2$ be two ordinary curves. Then $D_1D_2\equiv 0\mod 2$.\end{lemma}

\begin{proposition}\label{prop: possible singular fibers}{\rm (cf.\ \cite[Proposition 7.4]{K})} Let $\mathcal{E}:X_{(r,22-r,\delta)}\ra \mathbb{P}^1$ be an elliptic fibration. We recall that $s=r-10$ is the number of curves fixed by the cover involution, i.e.\ the number of special curves. Then the Kodaira type of the reducible fibers of $\mathcal{E}$ is contained in the following list:
\scriptsize
\begin{align}\label{table: possible fibers}
\begin{array}{|c|c|c|c|c|c|}
\hline\mbox{fiber } F_i&r_i&s_i&\mbox{special simple components}&\delta_i&|d_i|\\
\hline
III&2&1,0&1,0&3&2\\
\hline
I_2&2&0&0&2&2\\
\hline
I_{2k}, k\geq1&2k&k&k&2k&2k\\
k\leq r-10, k\leq (r-1)/2&&&&&\\
\hline
I_{2h}^*, h\geq 0&2h+5&h+1&0&2h+6&4\\
h\leq r-11,h\leq r/2-3&&&&&\\
\hline
IV^*&7&4&3&8&3\\
\hline
III^*&8&3&0&9&2\\
\hline
II^*&9&4&0&10&1\\
\hline
\end{array}\end{align}
\normalsize
where $r_i$ (resp. $d_i$, $\delta_i$) is $r(F_i)$ (resp.  $d(F_i)$, $\delta(F_i)$) and $s_i=s(F_i)$ is the number of special curves contained in the fiber $F_i$.  
\end{proposition}\proof Here we sketch the proof, which is similar to the one of \cite[Proposition 7.2]{K}. There are no fibers of type $IV$, $I_{2k+1}$ and $I_{2h+1}^*$ (see \cite[Proposition 7.4]{K}), essentially because Lemmas \ref{lemma:DB=2} and \ref{lemma:D1D2 even} imply an alternation among ordinary and special curves in the components of a reducible fiber. Analogously one obtains the number of special curves and of simple special curves given in \eqref{table: possible fibers} for the fibers $III$, $I_{2k}$, $I_{2h}^*$, $IV^*$, $III^*$, $II^*$. Clearly the number of special curves contained in a fiber can not be greater then the number of special curves $s=r-10$ on the surface. Analogously the number of components of a reducible fiber can not be greater then $r-2=\rho(X_{(r,22-r,\delta)})-2$. These two conditions give the restrictions on $k$ and $h$ in \eqref{table: possible fibers}. The others columns of Table \eqref{table: possible fibers} are standard in theory of the elliptic fibrations (cf. Section \ref{subsec: elliptic fibration} and \cite{SS}).\eprf

\begin{proposition}\label{prop: classification non finite MW}{\rm (cf.\ \cite[Proposition 8.1]{K})} Let $\mathcal{E}_r:X_{(r,22-r,\delta)}\ra \mathbb{P}^1$ be an elliptic fibration. Let us assume that all the special curves are contained in the reducible fibers of the fibration. Then $\mathcal{E}_r$ is one of the following fibrations:
\scriptsize
\begin{align}\label{table: possible fiber non trivial MW}\begin{array}{|c|c|c|c|c}
\hline
r&\mbox{reducible fibers}&\mbox{singular irreducible fibers}&\rk(\MW(\mathcal{E}_r))\\
\hline
12\leq r \leq 20&I_{2k_1}+I_{2k_2}&mI_1+nII&20-r\\
&k_1+k_2=r-10, 1\leq k_1\leq k_2&m+2n=44-2r&\\
\hline
15\leq r\leq 20&IV^*+I_{2k}&mI_1+nII&21-r\\
&k=r-14&m+2n=16-2k&\\
\hline
r=18&2IV^*&mI_1+(4-m/2)II&4\\
\hline
\end{array}\end{align}
\normalsize
\end{proposition}
\proof If all the special curves are contained in the reducible fibers of the fibration, the zero section is an ordinary curve. This implies (cf.\ \cite[Proposition 8.1]{K}) that the elliptic fibration $\mathcal{E}_r$ has only 2 reducible fibers and the zero section meets these fibers in a special curve. Thus there are no fibers of type $I_{2h}^*$, $III^*$, $II^*$, because the simple components of these fibers are all ordinary curves.\\
We number the 2 reducible fibers and we denote by $s_i$, $r_i$, $\delta_i$ the invariants defined in Proposition \ref{prop: possible singular fibers}. Since all the special curves are contained in the reducible fibers, the condition $s_1+s_2=s=r-10$ holds. Moreover, $s_1,s_2>1$ because the zero section meets a special curve in each reducible fiber. If the first fiber is of type $IV^*$, $s_1=4$ and thus if this fiber occurs then $s\geq 5$ hence $r\geq 15$. If both the reducible fibers are of type $I_{2k}$, then $s_1+s_2=k_1+k_2$; if one of them is $IV^*$ then $s_1+s_2=4+s_2$. These equalities together with the condition $s_1+s_2=s=r-10$ give the condition on $k_1$, $k_2$.\\
Since $X_{(r,22-r,\delta)}$ has Euler characteristic 24, $\sum_i\delta_i=24$. Since $\delta(I_1)=1$, $\delta(II)=2$, we have  $\sum_i\delta_i=\delta_1+\delta_2+m+2n$ which implies the condition on $m+2n$.\\
The Picard number of an elliptic surface is the sum of the rank of the trivial lattice and the one of the Mordell--Weil group, thus $\rk(\MW(\mathcal{E}_r))=r-(r_1+r_2)$. If both the fibers are of type $I_{2k}$, $r_1+r_2=2r-20$, if one of them is of type $IV^*$ then $r_1+r_2=7+2k_2=7+2(r-14)$.\eprf

\begin{lemma}\label{lemma: MW
 finite}{\rm (\cite[Lemma 2.4]{O}, \cite[Lemma 9.1]{K}, \cite[Lemma 9.2]{K})} Let $\mathcal{E}_r:X_{(r,22-r,\delta)}\ra\mathbb{P}^1$ be an elliptic fibration on $X_{(r,22-r,\delta)}$. Let us assume that there is at least a special curve not contained in the reducible fibers. Then $\MW(\mathcal{E}_r)$ is finite, the only fibers of type $I_{2k}$ are of type $I_2$ and both the fibers of type $I_2$ and $III$ do not contain special curves.\end{lemma}
\begin{lemma}\label{lemma: discriminant NS} Let $\mathcal{E}_r:X_{(r,22-r,\delta)}\ra\mathbb{P}^1$ be an elliptic fibration on $X_{(r,22-r,\delta)}$. Let us assume that there is a special curve which is not contained in the reducible fibers. Then there are no fibers of type $IV^*$. The discriminant of the trivial lattice of $\mathcal{E}_r$ is $2^c$ and the possible choices for $c$ are: $i)$ $c=a=22-r$, in this case $\MW(\mathcal{E}_r)=\{1\}$; $ii)$ $c=a+2=24-r$, in this case $\MW(\mathcal{E}_r)=\Z/2\Z$; $iii)$ $c=a+4=26-r$, in this case $\MW(\mathcal{E}_r)=(\Z/2\Z)^2$.\end{lemma}
\proof If the Mordell--Weil group is finite, the discriminant of the N\'eron--Severi group is $d(\Tr_{\mathcal{E}_r})/|\MW(\mathcal{E}_r)|^2$. We recall that $d(\NS(X_{(r,22-r,\delta)}))=2^{22-r}$. Let $d(\Tr_{\mathcal{E}_r})=2^c\prod_i p_i^{h_i}$, with $p_i$ prime numbers. We obtain \begin{equation}\label{eq: a and MW }|\MW(\mathcal{E}_r)|^2=d(\Tr_{\mathcal{E}_r})/d(\NS(X_{(r,22-r,\delta)}))=2^{c-(22-r)}\prod_i p_i^{h_i},\end{equation} thus $h_i$ and $c-(22-r)$ are non negative and even for each $i$. We recall that $d(\Tr_{\mathcal{E}_r})=\prod_i d(F_i)$ where $F_i$ are the reducible fibers of $\mathcal{E}_r$. By Proposition \ref{prop: possible singular fibers} and \ref{lemma: MW finite} the only reducible fiber on $\mathcal{E}_r$ with discriminant which is not a power of 2 is $IV^*$ and $d(IV^*)=3$, thus the only odd prime in \eqref{eq: a and MW } is $3$. Hence, if there is a fiber of type $IV^*$ then there are an even number $2x$ of fibers of type $IV^*$, $|\MW(\mathcal{E}_r)|^2=3^{2x}2^{c-(22-r)}$ and there are $x$ independent 3-torsion sections. By the height formula, if there is a 3-torsion section then there are at least 3 fibers of type $IV^*$. Since $r(IV^*)=6$  and $\sum_i r(F_i)+2=\rk(\Tr_{\mathcal{E}_r})\leq 20$, there are at most 3 fibers of type $IV^*$. Thus, if there is a fiber of type $IV^*$ then there are exactly 3 fibers of type $IV^*$, but this is impossible, since 3 is not even. So there are no fibers of type $IV^*$. Hence the reducible fibers are of type $III$, $I_2$, $I_{2h}^*$, $III^*$, $II^*$ and this implies that the Mordell--Weil group is trivial, or $\Z/2\Z$ or $(\Z/2\Z)^2$. So $|\MW(\mathcal{E}_r)|^2$ is $1$ or $2^2$, $2^4$. By \eqref{eq: a and MW }, according to these possibilities one finds  $c-22+r=0$, $c-22+r=2$, $c-22+r=4$ respectively.\eprf

\begin{proposition}\label{prop: classification finite MW}{\rm (cf.\ \cite[Theorem 1.2]{K})} Let $\mathcal{E}_r:X_{(r,22-r,\delta)}\ra \mathbb{P}^1$ be an elliptic fibration on $X_{(r,22-r,\delta)}$. Let us assume that there is at least one special curve which is not contained in the reducible fibers of the fibration. Then $\mathcal{E}_r$ is one of the elliptic fibrations listed in the tables in Section \ref{subsec: existence case b}: the singular fibers are $m_1I_2+m_2III+n_1I_1+n_2II$ and the ones listed in the second column, the values of $m_1,m_2,n_1,n_2$  are non--negative and have the properties listed in the third and fourth column.\\
In case $r=11$ there are no elliptic fibrations such that the special curve is not contained in a fiber.
\end{proposition}
\proof By Lemma \ref{lemma: MW
 finite} there are no fibers of type $I_{2\nu}$, $\nu\geq2$ and by Lemma \ref{lemma: discriminant NS} there are no fibers of type $IV^*$. By Lemma \ref{lemma: MW
 finite} the rank of the Mordell--Weil is zero, thus $\rk(\NS(\mathcal{E}_r))=\rk(\Tr_{\mathcal{E}_r})$. Hence the reducible fibers $F_i$ satisfy the condition $\sum_i(r(F_i)-1)=r-2$. There is at least one special curve which is not contained in the reducible fibers, thus $\sum_i s(F_i)\leq s-1=r-11$, where $s(F_i)=0$ if $F_i$ is of type $I_2$ or $III$, by Lemma \ref{lemma: MW
 finite}, and the values of $s(F_i)$ are listed in Table \eqref{table: possible fibers} for the other cases. Let $\delta$ be the sum of the Euler characteristics of the reducible fibers. Since $X$ is a K3 surface and the Euler numbers of the $I_1$ and $II$ fiber are 1, respectively 2, we have that $24-\delta=n_1+2n_2$.
By Lemma \ref{lemma: MW finite}, $\sum_i d(F_i)=2^{c}$ with $c=22-r,\ 24-r,\ 26-r$. To recap the fibrations listed in the statement satisfy the following conditions:\\
$\bullet$ the reducible fibers are of type $II^*$, $III^*$, $I_{2k}^*$, $I_2$, $III$ (by Proposition \ref{prop: possible singular fibers} and Lemma \ref{lemma: MW
 finite});\\
 $\bullet$ $\sum_i(r(F_i)-1)=r-2$;\\ 
 $\bullet$ $\sum_i s(F_i)\leq s-1=r-11$ ($s(I_2)=s(III)=0$);\\
 $\bullet$ $\sum_i \delta(F_i)\leq 24$; \\
 $\bullet$ $\sum_i d(F_i)=2^{c}$ with $c=22-r,\ 24-r,\ 26-r$.\\
Moreover, in cases $c=24-r$ or $c=26-r$ there is another condition. Indeed, if $c=24-r$ (resp. $26-r$) then $\MW=\Z/2\Z$ (resp. $\MW=(\Z/2\Z)^2$). The height formula (cf. Section \ref{subsec: elliptic fibration}) imposes conditions to the reducible fibers to have torsion sections. For example in case $r=19$ the configuration $II^*+I_2^*+3I_2$ satisfies the conditions listed before with $c=5=24-r$. By the height formula such an ellliptic fibration has no torsion sections hence we can exclude this case. The height formula implies the properties of the intersection between the torsion sections and the reducible fibers given in the last column of the tables. \eprf
\begin{rem}{\rm Propositions \ref{prop: classification finite MW}, \ref{prop: classification non finite MW} imply that there are no elliptic fibrations on $X_{(11,11,1)}$.\erem}\end{rem}
\begin{proposition}\label{prop: iota on elliptic fibration} Let $\mathcal{E}_r:X_{(r,22-r,\delta)}\ra \mathbb{P}^1$ be an elliptic fibration on $X_{(r,22-r,\delta)}$. The involution $\iota$ acts on $\mathcal{E}_r$ either as an involution of type a) (cf.\ Example \ref{ex: involution on the basis}) or as an involution of type b) (cf.\ Example \ref{ex: hyperelliptic involution}). In particular, $\iota$ is of type a) if and only if all the special curves are contained in the reducible fibers and $\iota$ is of type  b) if and only if there is a special curve not contained in the reducible fibers. \end{proposition}
\proof Let $\mathcal{E}_r:X_{(r,22-r,\delta)}\ra \mathbb{P}^1$ be the elliptic fibration induced by the divisor $F$ (i.e.\ $F$ is the class of the fiber of this fibration or equivalently the map associated to the nef divisor $F$ is $\mathcal{E}_r:X_{(r,22-r,\delta)}\ra \mathbb{P}^1$). Since the involution $\iota$ acts as the identity on the N\'eron--Severi group, $\iota^*(F)=F$. Thus a fiber is sent to a fiber by $\iota$. If each fiber is sent to itself, the involution $\iota$ acts trivially on the base of the elliptic fibration and restricts to an involution of each fiber.  The generic fiber of $\mathcal{E}_r$ is an elliptic curve, so its automorphisms are the composition of the translations by rational points and the hyperelliptic involution. The hyperelliptic involution fixes the zero section, the translations send it to another section. Since $\iota$ acts trivially on the N\'eron--Severi group, we conclude that $\iota$ is the hyperelliptic involution, i.e.\ of type b). On the other hand, if the generic fiber is sent to another fiber then the involution acts on the basis as an involution of $\mathbb{P}^1$ (and up to a projectivity we can assume that it acts on $\mathbb{P}^1$ fixing the points $0$ and $\infty$). In particular it restricts to an automorphism of the two fibers over $0$ and $\infty$. On these fibers it acts as the identity or as a translation by a point of order 2 or as a non--symplectic automorphism $\beta$ (like the ones describe in Examples \ref{ex: hyperelliptic involution}, \ref{ex: non symplectic involution Q-P}). In the first case, the involution acts only on the base of the fibration, and is an involution of type a). If it acts as a translation by a point of order 2 on the fixed fibers, then it sends the zero section to a section of order 2, but these two sections are different classes in the N\'eron--Severi group, which is a contradiction, because $\iota$ acts as the identity on the N\'eron--Severi group. If it acts as a non--symplectic involution $\beta$, then $\iota$ is the composition of two commuting non--symplectic involutions. But the composition of two commuting non--symplectic involutions is a symplectic involution and $\iota$ is non--symplectic.\\
If $\iota$ is of type a), it fixes the fibers over $\tau=0$, $\tau=\infty$ and switches the fibers over $\tau$ and $-\tau$. The curves which are not contained in an invariant fiber are clearly not fixed. Thus, if $\iota$ is of type a) then the special curves are all contained in two fibers, which are the only reducible fibers of the fibration.\\
If $\iota$ is the hyperelliptic involution (i.e. of type b)) on $\mathcal{E}_r$, then the zero section is a fixed curve. By Definition \ref{def: special curve} the zero section is a special curve and, since it is a section, it is not contained in the reducible fibers. Thus, if $\iota$ is the hyperelliptic involution on $\mathcal{E}_r$ then there is at least one special curve not contained in the reducible fibers.\eprf
\begin{rem}{\rm If an elliptic fibration admits an involution of type a) then the type of the fiber over $\overline{\tau}\neq 0,\infty$ is equal to the type of the fiber over $-\overline{\tau}$, because they are exchanged by the involution. In particular this means that the numbers $m$ and $n$ in Table \eqref{table: possible fiber non trivial MW} are even.\erem}\end{rem}

\section{Elliptic fibrations on $X_{(r,22-r,\delta)}$: existences in case a)}\label{section:existence case 1}
In Table \eqref{table: possible fiber non trivial MW} we listed the possible elliptic fibrations on $X_{(r,22-r,\delta)}$ such that all the special curves are contained in reducible fibers. To prove the existence of these elliptic fibrations it suffices to find a nef class $F$ in $\NS(X_{(r,22-r,\delta)})$ such that the induced map, $X_{(r,22-r,\delta)}\rightarrow \mathbb{P}(H^0(X_{(r,22-r,\delta)}, F))$, is an elliptic fibration with the required properties. In particular $F$ turns out to be the class of a fiber of the fibration and has a trivial self intersection.
%
In order to find $F$, it is necessary to give a description of $\NS(X_{(r,22-r,\delta)})$. For each $r$ we describe $X_{(r,22-r,\delta)}$ as double cover of $\mathbb{P}^2$ branched over a chosen particular sextic and this give a $\Q$-basis of $\NS(X_{(r,22-r,\delta)})$. We will identify the nef divisor $F$ as linear combinations of this $\Q$-basis.\\ 
This section is devoted to the proof of the following proposition.
\begin{proposition}\label{prop: existence case 1} All the elliptic fibrations listed in Table \eqref{table: possible fiber non trivial MW} appear as elliptic fibrations on $X_{(r,22-r,\delta)}$ except: $r=20:$ $I_{14}+I_6$, $I_{12}+I_8$, $IV^*+I_{12}$.\\
The quotient $X_{(r,22-r,\delta)}/\iota$ is a rational elliptic surface admitting the following elliptic fibrations:
\begin{align}\label{table: existence fibration non trivial MW}\begin{array}{|c|c|c|c|c}
\hline
r&\mbox{reducible fibers}&\mbox{singular irreducible fibers}&\rk(\MW(\mathcal{E}_r))\\
\hline
12\leq r \leq 20&I_{k_1}+I_{k_2},&(m/2)I_1+(n/2)II&20-r\\
&k_1+k_2=r-10, 1\leq k_1\leq k_2&&\\
\hline
15\leq r\leq 20&IV+I_{k},&(m/2)I_1+(n/2)II&21-r\\
&k=r-14&&\\
\hline
r=18&2IV&(m/2)I_1+(n/2)II&4\\
\hline
\end{array}\end{align}
where $m$ and $n$ are the numbers given in Table \eqref{table: possible fiber non trivial MW}. 
\end{proposition}

\subsection{$r=20$}\label{section: non trivial MW r=20} In this case, the branch sextic consists of six lines with 4 triple points. Thus the components of the sextic are $l_i$, $i=1,\ldots,6$ and we assume that the triple points are $A:=l_1\cap l_2\cap l_3$, $B:=l_1\cap l_5\cap l_6$, $C:=l_3\cap l_4\cap l_5$, $D:=l_2\cap l_4\cap l_6$. The double points are $P_{1,4}$, $P_{2,5}$, $P_{3,6}$. The N\'eron--Severi group is generated by $\mathcal{L}_i$, $i=1,\ldots, 6$, $E_{A}$, $E_{A,1}$, $E_{A,2}$, $E_{A,3}$, $E_B$, $E_{B,1}$, $E_{B,5}$, $E_{B,6}$, $E_C$, $E_{C,3}$, $E_{C,4}$, $E_{C,5}$, $E_D$, $E_{D,2}$, $E_{D,4}$, $E_{D,6}$, $E_{P_{h,k}}$, $(h,k)=\{(1,4),\ (2,5),\ (3,6)\}$.\\
An elliptic fibration with a fiber of type $I_{18}$ is associated to the divisor $F:=\mathcal{L}_1+E_{A_1}+E_A+E_{A,2}+\mathcal{L}_2+E_{D,2}+E_{D}+E_{D,6}+\mathcal{L}_6+E_{B,6}+E_B+E_{B,5}+\mathcal{L}_5+E_{C,5}+E_C+E_{C,4}+\mathcal{L}_4+E_{P_{1,4}}$. The classes $E_{A,3}$, $E_{C,3}$, $E_{P_{3,6}}$ are the classes of sections. The fibration $\phi_F:X_{(20,2,1)}\ra \mathbb{P}^1$ has a fiber of type $I_{18}$ (whose components are the summands $\mathcal{L}_1$, $E_{A_1}$,\ldots, $E_{P_{1,4}}$) and a fiber of type $I_2$ (because we already proved that the unique admissible elliptic fibration on $X_{(20,2,1)}$ with a fiber of type $I_{18}$ has $I_{18}+I_2$ as reducible fibers). The Mordell--Weil group of this elliptic fibration is $\Z/3\Z$, indeed we proved that the rank of the Mordell--Weil is 0, thus $4=d(\NS(X_{(22,2,1)}))=(18\cdot 2)/|\MW|^2$ hence $|\MW|=3$. \\
Similarly, we prove that there exist the elliptic fibrations:\\
$I_{16}+I_4$, given by the divisor $F:=\mathcal{L}_1+E_{P_{1,4}}+\mathcal{L}_4+E_{C,4}+E_C+E_{C,5}+\mathcal{L}_5+E_{B,5}+E_{B}+E_{B,6}+\mathcal{L}_6+E_{P_{3,6}}+\mathcal{L}_3+E_{A,3}+E_A+E_{A,1}$, zero section $s:=E_{P_{2,5}}$ and Mordell--Weil group $\Z/4\Z$, generated by $E_{A,2}$;\\
$2I_{10}$, given by the divisor $F:=\mathcal{L}_1+E_{A,1}+E_A+E_{A,2}+\mathcal{L}_2+E_{P_{2,5}}+\mathcal{L}_5+E_{B,5}+E_B+E_{B,1}$, zero section $s:=E_{A,3}$ and Mordell--Weil group $\Z/5\Z$, generated by $E_{D,2}$.\\
There is no elliptic fibration with fibers configuration $I_{14}+I_6$, since if it exists then $|\MW|^2=(14\cdot 6)/d(\NS(X_{20,2,1}))=21$, which is a contradiction, because 21 is not a square. Similarly there are no elliptic fibrations with fibers $I_{12}+I_8$, $IV^*+I_{12}$.

\subsection{$r=19$}\label{section: existence trivial MW r=19} The sextic we consider consists of six lines $l_i$, $i=1,\ldots, 6$ with 3 triple points $A:=l_1\cap l_2\cap l_3$, $B:=l_1\cap l_5\cap l_6$, $C:=l_3\cap l_4\cap l_5$. The double points are $P_{1,4}$, $P_{2,4}$, $P_{2,5}$, $P_{2,6}$, $P_{3,6}$, $P_{4,6}$. On $X_{(19,3,1)}$ there exist the following elliptic fibrations: 
\scriptsize{
\begin{align*}
\begin{array}{|c|c|c|c|}
\hline
\mbox{reducible fibers}&F&s\\
\hline
I_{16}+I_2&\mathcal{L}_1+E_{P_{1,4}}+\mathcal{L}_4+E_{C,4}+E_C+E_{C,5}+\mathcal{L}_5+E_{B,5}&E_{P_{2,5}}\\
&+E_{B}+E_{B,6}+\mathcal{L}_6+E_{P_{3,6}}+\mathcal{L}_3+E_{A,3}+E_A+E_{A,1}&\\
\hline
I_{14}+I_4&\mathcal{L}_1+E_{P_{1,4}}+\mathcal{L}_4+E_{P_{4,6}}+\mathcal{L}_6+E_{B,6}+E_{B}&E_{A_3}\\
&+E_{B,5}+\mathcal{L}_5+E_{P_{2,5}}+\mathcal{L}_2+E_{A,2}+E_A+E_{A,1}&\\
\hline
I_{12}+I_6&\mathcal{L}_1+E_{P_{1,4}}+\mathcal{L}_4+E_{P_{2,4}}+\mathcal{L}_2+E_{P_{2,6}}&E_{C,4}\\
&+\mathcal{L}_6+E_{P_{3,6}}+\mathcal{L}_3+E_{A,3}+E_A+E_{A,1}&\\
\hline
I_{10}+I_8&\mathcal{L}_1+E_{P_{1,4}}+\mathcal{L}_4+E_{P_{2,4}}+\mathcal{L}_2+E_{A,2}+E_A+E_{A,1}&E_{A,3}\\
\hline
IV^*+I_{10}& \mathcal{L}_1+\mathcal{L}_3+\mathcal{L}_2+2(E_{A,1}+E_{A_2}+ E_{A,3})+3E_A&E_{P_{2,4}}\\
\hline
\end{array}
\end{align*}}
\normalsize
In case $I_{10}+I_8$ the fiber associated to $F$ is $I_8$ and in case $IV^*+I_{10}$ it is $IV^*$.

\subsection{$r=18$}\label{section: non trivial MW r=18} In case $r=18$ the invariants $(r,a)$ do not identify uniquely the lattice $N_{(18,4,\delta)}$, indeed there are two possible values for $\delta$ (cf. Corollary \ref{cor: double cover 18}).\\
In \cite{O} it is proved that there exist the elliptic fibrations with reducible fibers $I_{12}+I_4$, $2I_8$, $2IV^*$ on $X_{(18,4,0)}$.\\
By Corollary \ref{cor: double cover 18}, we can associate the K3 surface $X_{(18,4,1)}$ to the sextic $B_{(18,4,1)}$ which consists of the lines $l_i$, $i=1,\ldots, 6$, with triple points $A:=l_1\cap l_2\cap l_3$, $B:=l_1\cap l_5\cap l_6$ and double points $P_{1,4}$, $P_{2,4}$, $P_{2,5}$, $P_{2,6}$, $P_{3,4}$, $P_{3,5}$, $P_{3,6}$, $P_{4,5}$, $P_{4,6}$. 
On $X_{(18,4,1)}$ there exist the following elliptic fibrations:
\scriptsize
\begin{align*}
\begin{array}{|c|c|c|c|}
\hline
\mbox{reducible fibers}&F&s\\
\hline
I_{14}+I_2&\mathcal{L}_1+E_{P_{1,4}}+\mathcal{L}_4+E_{P_{4,6}}+\mathcal{L}_6+E_{B,6}+E_{B}&E_{A_3}\\
&+E_{B,5}+\mathcal{L}_5+E_{P_{2,5}}+\mathcal{L}_2+E_{A,2}+E_A+E_{A,1}&\\
\hline
I_{12}+I_4&\mathcal{L}_1+E_{P_{1,4}}+\mathcal{L}_4+E_{P_{2,4}}+\mathcal{L}_2+E_{P_{2,6}}&E_{P_{4,5}}\\
&+\mathcal{L}_6+E_{P_{3,6}}+\mathcal{L}_3+E_{A,3}+E_A+E_{A,1}&\\
\hline
I_{10}+I_6&\mathcal{L}_1+E_{A,1}+E_A+E_{A,2}+\mathcal{L}_2+E_{P_{2,5}}+\mathcal{L}_5+E_{B,5}+E_B+E_{B,1}&E_{A,3}\\
\hline
2I_8&\mathcal{L}_1+E_{P_{1,4}}+\mathcal{L}_4+E_{P_{2,4}}+\mathcal{L}_2+E_{A,2}+E_A+E_{A,1}&E_{A,3}\\
\hline
IV^*+I_{8}& \mathcal{L}_1+\mathcal{L}_2+\mathcal{L}_3+2(E_{A,1}+E_{A_2}+E_{A,3})+3E_A&E_{P_{2,4}}\\
\hline 
\end{array}
\end{align*}
\normalsize
In case $2I_8$ the fiber associated to $F$ is one of the fiber of type $I_8$. To prove that it gives the fibrations with fibers $2I_8$ and not $IV^*+I_8$ it is necessary to identify also the classes which form the second fiber of type $I_8$ (they are $\mathcal{L}_i$, $i=3,5,6$, $E_{P_{3,5}}$, $E_{P_{3,6}}$, $E_{B,5}$, $E_B$, $E_{B,6}$). 
\subsection{$r=17$} The sextic we consider consists of six lines $l_i$, $i=1,\ldots 6$ with a triple point $A:=l_1\cap l_2\cap l_3$. There are 12 double points: $P_{1,4}$, $P_{1,5}$, $P_{1,6}$, $P_{2,4}$, $P_{2,5}$, $P_{2,6}$, $P_{3,4}$, $P_{3,5}$, $P_{3,6}$, $P_{4,5}$,$P_{4,6}$, $P_{5,6}$. 
On $X_{(17,5,1)}$ there exist the following elliptic fibrations:
\scriptsize
\begin{align*}
\begin{array}{|c|c|c|c|}
\hline
\mbox{reducible fibers}&F&s\\ \hline
I_{12}+I_2&\mathcal{L}_1+E_{P_{1,4}}+\mathcal{L}_4+E_{P_{2,4}}+\mathcal{L}_2+E_{P_{2,6}}&E_{P_{4,5}}\\
&+\mathcal{L}_6+E_{P_{3,6}}+\mathcal{L}_3+E_{A,3}+E_A+E_{A,1}&\\
\hline
I_{10}+I_4&\mathcal{L}_1+E_{A,1}+E_A+E_{A,2}+\mathcal{L}_2+E_{P_{2,5}}+\mathcal{L}_5+E_{P_{5,6}}+\mathcal{L}_6+E_{P_{1,6}}&E_{A,3}\\
\hline
I_8+I_6&\mathcal{L}_1+E_{P_{1,4}}+\mathcal{L}_4+E_{P_{2,4}}+\mathcal{L}_2+E_{A,2}+E_A+E_{A,1}&E_{A,3}\\
\hline
IV^*+I_{6}& \mathcal{L}_1+\mathcal{L}_2+\mathcal{L}_3+2(E_{P_{1,4}}+E_{P_{2,4}}+E_{P_{3,4}})+3\mathcal{L}_4&E_{A,1}\\
\hline 
\end{array}
\end{align*}
\normalsize
\subsection{$r=16$} This case was already analyzed in \cite{K}, where the sextic considered consists of 6 lines. Here we report the fibers given in \cite{K}:
\scriptsize
\begin{align*}
\begin{array}{|c|c|c|c|}
\hline
\mbox{reducible fibers}&F&s\\
\hline
I_{10}+I_2&\mathcal{L}_1+E_{P_{1,2}}+\mathcal{L}_2+E_{P_{2,3}}+\mathcal{L}_3+E_{P_{3,4}}+\mathcal{L}_4+E_{P_{4,5}}+\mathcal{L}_5+E_{P_{1,5}}&E_{P_{1,6}}\\
\hline
I_8+I_4&\mathcal{L}_1+E_{P_{1,2}}+\mathcal{L}_2+E_{P_{2,3}}+\mathcal{L}_3+E_{P_{3,4}}+\mathcal{L}_4+E_{P_{1,4}}&E_{P_{1,6}}\\
\hline
2I_6&\mathcal{L}_1+E_{P_{1,2}}+\mathcal{L}_2+E_{P_{2,3}}+\mathcal{L}_3+E_{P_{1,3}}&E_{P_{1,6}}\\
\hline
IV^*+I_{4}& \mathcal{L}_1+\mathcal{L}_3+\mathcal{L}_4+2(E_{P_{1,2}}+E_{P_{2,3}}+E_{P_{2,4}})+3\mathcal{L}_2&E_{P_{1,6}}\\
\hline 
\end{array}
\end{align*}
\normalsize
\subsection{$r=15$} Here we consider a sextic which consists of 4 lines $l_i$, $i=1,2,3,4$  and 1 conic $c_5$, without triple points. The double points are $\{R^1_{j,5}, R^2_{j,5}\}=l_j\cap c_5$, $j=1,2,3,4$, $P_{h,k}=l_h\cap l_k$, $h,k=1,2,3,4$, $h\neq k$. 
On $X_{(15,7,1)}$ there exist the following elliptic fibrations:
\scriptsize
\begin{align*}
\begin{array}{|c|c|c|c|}
\hline
\mbox{reducible fibers}&F&s\\
\hline
I_8+I_2&\mathcal{L}_1+E_{P_{1,2}}+\mathcal{L}_2+E_{P_{2,3}}+\mathcal{L}_3+E_{P_{3,4}}+\mathcal{L}_4+E_{P_{1,4}}&E_{R^1_{1,5}}\\
\hline
I_6+I_4&\mathcal{L}_1+E_{P_{1,2}}+\mathcal{L}_2+E_{P_{2,3}}+\mathcal{L}_3+E_{P_{1,3}}&E_{R^1_{1,5}}\\
\hline
IV^*+I_{2}& \mathcal{L}_1+\mathcal{L}_3+\mathcal{L}_4+2(E_{P_{1,2}}+E_{P_{2,3}}+E_{P_{2,4}})+3\mathcal{L}_2&E_{R^1_{1,5}}\\
\hline 
\end{array}
\end{align*}
\normalsize

\subsection{$r=14$} Here we consider a sextic made up of 2 lines $l_i$, $i=1,2$ and 2 conics $c_j$, $j=3,4$, without triple points. The double points are $\{R^1_{j,i}, R^2_{j,i}\}=l_j\cap c_i$, $j=1,2$, $i=3,4$, $P_{1,2}$, $\{Q^1_{3,4}, Q^2_{3,4},Q^3_{3,4}, Q^4_{3,4}\}=c_3\cap c_4$. On $X_{(14,8,1)}$ there exist the following elliptic fibrations:
\scriptsize
\begin{align*}
\begin{array}{|c|c|c|c|}
\hline
\mbox{reducible fibers}&F&s\\
\hline
I_6+I_2&\mathcal{L}_1+E_{P_{1,2}}+\mathcal{L}_2+E_{R^1_{(2,3)}}+\mathcal{C}_3+E_{R^1_{1,3}}&E_{R^1_{1,4}}\\
\hline
2I_4& \mathcal{L}_1+E_{R^1_{1,3}}+\mathcal{C}_3+E_{R^2_{1,3}}&E_{R^1_{1,4}}\\
\hline 
\end{array}
\end{align*}
\normalsize
\subsection{$r=13$} We consider a sextic which consists of 3 conics $c_i$, $i=1,2,3$, without triple points. We put $\{Q^1_{i,j},Q^2_{i,j}, Q^3_{i,j},Q^4_{i,j}\}=c_i\cap c_j$, $i\neq j$, $i,j=1,2,3$. There is an elliptic fibration with reducible fibers $I_4+I_2$ associated to the divisor $F:=\mathcal{C}_1+E_{Q^1_{1,2}}+\mathcal{C}_2+E_{Q^2_{1,2}}$. 
\subsection{$r=12$}  We consider a sextic made up of a quartic $q$ with 3 singular points $A$, $B$, $C$ and a conic $c$. There is an elliptic fibration with reducible fibers $2I_2$ associated to the divisor $F:=\mathcal{Q}+E_{A}$. 

\subsection{Quotients by $\iota$} We proved in Proposition \ref{prop: iota on elliptic fibration} that the involution $\iota$ acts on the elliptic fibrations of Table \ref{table: possible fiber non trivial MW} as an involution on the base of the fibration. As in Example \ref{ex: involution on the basis}, in this case the elliptic fibration, $\mathcal{E}_{X_{r}}$, on $X_{(r,22-r,\delta)}$ can be obtained by an elliptic fibration $\mathcal{E}_{R_r}$ on a rational elliptic surface $R_r$, by a base change of order 2: $$\begin{array}{ccc}R_r&\la &X_{(r,22-r,\delta)}\\
\downarrow&&\downarrow\\
\mathbb{P}^1_t&\la&\mathbb{P}^1_\tau
\end{array}$$
where $t=\tau^2$. The fiber over $t=0$ (resp. $t=\infty$) corresponds to the fiber over $\tau=0$ (resp. $\tau=\infty$). If the fiber over $t=0,\infty$ is of type $I_{h}$ or $IV$ the fiber over $\tau=0,\infty$ is of type $I_{2h}$ or $IV^*$ respectively. Each fiber over $t\neq 0,\infty$ corresponds to two fibers of the same type over $\tau$. The quotient $X_{(r,22-r,\delta)}/\iota$ is $R_r$ and this implies the description of the reducible fibers on $R_r$ given in the statement of Proposition \ref{prop: existence case 1}. Since the sections of $\mathcal{E}_{R_{r}}$ lift to sections of $\mathcal{E}_{X_r}$ generically the Mordell--Weil group of the two elliptic fibrations is the same. In our case one can explicitly verify that the rank of the Mordell--Weil group of $\mathcal{E}_{R_r}$ and $\mathcal{E}_{X_r}$ are the same recalling that the Picard number of $\mathcal{E}_{R_r}$ is 10 and thus $\rk(\MW(\mathcal{E}_{R_r}))=10-k_1-k_2$, (resp. $10-3-k$, $10-6$) if the reducible fibers of $\mathcal{E}_{X_r}$ are $I_{2k_1}+I_{2k_2}$ (resp. $IV^*+I_{2k}$, $2IV^*$).

For certain values of $r$ ($12\leq r\leq 15$) the Mordell--Weil group of $\mathcal{E}_{R_r}$ is known as well as a geometrical construction for the pencil of cubics in $\mathbb{P}^2$ associated to this elliptic fibration (\cite[Theorem 10.4]{S}, \cite{F}, \cite{FT}, \cite{Sa}). Thus the K3 surface $X_{(r,22-r,\delta)}$ is obtained from a known elliptic fibration by a base change of order 2. In particular, if $r=12$ then $\mathcal{E}_{R_{12}}$ is an elliptic fibration without reducible fibers. In this case $\MWL(\mathcal{E}_{R_{12}})\simeq E_8$ and the pencil of cubics associated to this elliptic fibration is the generic one. If $r=13$, $\MWL(\mathcal{E}_{R_{13}})\simeq E_7^*$;\\ if $r=14$ and the reducible fibers of $\mathcal{E}_{X_{14}}$ are $I_6+I_2$  then $\MWL(\mathcal{E}_{R_{14}})=E_6^*$;\\
if $r=14$ and the reducible fibers of $\mathcal{E}_{X_{14}}$ are $I_4+I_4$  then $\MWL(\mathcal{E}_{R_{14}})=D_6^*$;\\
if $r=15$ and the reducible fibers of $\mathcal{E}_{X_{15}}$ are $I_8+I_2$  then $\MWL(\mathcal{E}_{R_{14}})=D_5^*$;\\
if $r=15$ and the reducible fibers of $\mathcal{E}_{X_{15}}$ are $I_6+I_4$  then $\MWL(\mathcal{E}_{R_{14}})=A_5^*$.

\section{Elliptic fibrations on $X_{(r,22-r,\delta)}$: existence and equations in case b)}\label{sec: existence case 2)}
In this section we prove the existence of the elliptic fibrations listed in Proposition \ref{prop: classification finite MW}  (these are the elliptic fibrations with at least one special curve not contained in the reducible fibers).  As in case a) (Section \ref{section:existence case 1}), to prove the existence of these fibrations, it suffices to find a nef divisor $F$ in $\NS(X_{(r,22-r,\delta)})$ whose associated map is an elliptic fibration with the required properties. The class $F$ is obtained in a very geometrical way from the properties of the double cover $X_{(r,22-r,\delta)}\ra \mathbb{P}^2$ (see Construction \ref{construction} and Examples \ref{ex: pencil-17-5-1}, \ref{example: 17-5-1 not line}). The geometrical construction of $F$ allows us to write explicit equations for the classified elliptic fibrations (see Example \ref{ex: B13}).
%
%

\begin{construction}\label{construction}{\rm Let us denote by $\mathcal{C}_6$ the branch sextic of the double cover $X_{(r,22-r,\delta)}\rightarrow \mathbb{P}^2$. Let us consider a rational curve $\gamma\subset\mathbb{P}^2$ such that the intersection of $\gamma$ with $\mathcal{C}_6$ consists of exactly four smooth points of $\mathcal{C}_6$ and a certain numbers of double points of $\mathcal{C}_6$. The pull-back of $\gamma$ to $X_{(r,22-r,\delta)}$ is a double cover of $\PP^1\simeq \gamma$ branched over 4 points, i.e.\ an elliptic curve. Let $\mathcal{P}$ be a pencil of rational curves in $\mathbb{P}^2$ whose general member intersects $\mathcal{C}_6$ in the base locus of the pencil and in other 4 smooth points. It gives an elliptic fibration on $X_{(r,22-r,\delta)}$. Indeed $\mathcal{P}$ is a 1-dimensional family of curves, parametrized by $\mathbb{P}^1$, and it is associated to a map $\mathcal{E}_{\mathcal{P}}: X\rightarrow \PP^1$ such that the general fiber is an elliptic curve. Reducible fibers of the fibration correspond to curves in $\mathcal{P}$ through singular points of the sextic or to reducible members of $\mathcal{P}$.

Particular choices for $\mathcal{P}$ are the pencils of lines:  let $P$ be a double point of $\mathcal{C}_6$ and let $\mathcal{P}$ be the pencil of lines of $\mathbb{P}^2$ through $P$. A generic line $l$ of $\mathcal{P}$ intersects $\mathcal{C}_6$ in $P$ and in 4 smooth points. In order to construct a smooth model of $X_{(r,22-r,\delta)}$ we blow up $\mathbb{P}^2$ in the singular points of $\mathcal{C}_6$ as described in Proposition \ref{prop: double cover of P2}. Let $\widetilde{l}$ be the strict transform of $l$ and $\widetilde{\mathcal{C}_6}$ be the strict transform of $\mathcal{C}_6$. Let $\mathcal L$ be the pull-back of $\widetilde{l}$ on $X_{(r,22-r,\delta)}$. The curve $\mathcal{L}$ is a smooth fiber of the fibration $\mathcal{E}_{\mathcal{P}}$. The curve $E_P$ (with the notation of Section \ref{section:double cover}) is a bisection of the fibration $\mathcal{E}_{\mathcal{P}}$, indeed $\widetilde{l}$ meets $e_P$ in a point $Q\in\widetilde{\mathbb{P}^2}$ that is not in the branch locus, hence it corresponds to 2 points $Q_1,Q_2$ in $X_{(r,22-r,\delta)}$ and $\mathcal L\cap E_P=\{Q_1,Q_2\}$.

Let us now assume that $P$ is a triple point of $\mathcal{C}_6$ and $\mathcal{P}$ is again the pencil of lines of $\mathbb{P}^2$ through $P$. Let $l$ be a generic line of the pencil $\mathcal{P}$: $l$ intersects $\mathcal{C}_6$ in $P$ and in three smooth points. In order to construct a smooth model of $X_{(r,22-r,\delta)}$ we blow up $\mathbb{P}^2$ in the singular points of $\mathcal{C}_6$ and in particular we introduce four rational curves on $P$. One of the exceptional curves on $P$, $e_P$, is in the branch locus of the double cover $X_{(r,22-r,\delta)}\rightarrow \widetilde{\mathbb{P}}^2$. The intersection of $\widetilde{l}$ with the branch locus of the double cover $X_{(r,22-r,\delta)}\rightarrow \widetilde{\mathbb{P}}^2$ consists of four points, three of them are the intersection of $\widetilde{l}$ with $\widetilde{C_6}$ and the other one is the intersection of $\widetilde{l}$ with $e_P$. Thus, for a generic line $l\in\mathcal{P}$, $\widetilde{l}$ intersects the branch locus of the double cover $X_{(r,22-r,\delta)}\rightarrow \widetilde{\mathbb{P}^2}$ in four points. Hence, the pencil $\mathcal{P}$ induces an elliptic fibration $\mathcal{E}_{\mathcal{P}}$ on $X_{(r,22-r,\delta)}$.  The curve $E_P$ is a section of the fibration $\mathcal{E}_{\mathcal{P}}$. Indeed, $\widetilde{l}$ meets $e_P$ in a point $Q$  which is in the branch locus, since $e_P$ is in the branch locus. Thus, on $X_{(r,22-r,\delta)}$, the intersection $\mathcal L\cap E_P$ consists of a single point.\erem}\end{construction}

Our purpose is to find pencils $\mathcal{P}$ of rational curves in $\mathbb{P}^2$ which correspond to the elliptic fibrations listed in Proposition \ref{prop: classification finite MW}. Once the branch sextic associated to the K3 surface $X_{(r,22-r,\delta)}$ is fixed, one has to find a pencil $\mathcal{P}$ such that special members of $\mathcal{P}$ give certain chosen reducible fibers on $\mathcal{E}_{\mathcal{P}}$. In the following examples we clarify this strategy showing how to choose the pencil $\mathcal{P}$ in order to obtain particular reducible fibers.

\begin{example}\label{ex: pencil-17-5-1}{\rm
Let us consider the sextic made up of six lines $l_1,\ldots,l_6$ with a triple point $A=l_1\cap l_2\cap l_3$ and 12 double points.
The K3 surface $X$ obtained as double cover of $\mathbb P^2$ branched over this sextic has invariants $(r,a,\delta)=(17,5,1)$ (see Corollary \ref{cor: double cover 1}), thus $X=X_{(17,5,1)}$.\\
We consider the pencil of lines $\mathcal P$ through the double point $P_{5,6}$: a generic line of $\mathcal P$ meets the sextic in 4 points (outside the base point of the pencil $P_{5,6}$). Hence the pull-back of the generic line to the double cover $X$ is an elliptic curve and $\mathcal{P}$ induces an elliptic fibration $\mathcal{E}_{\mathcal{P}}$ on $X$. 
The line $m$ through $P_{5,6}$ and $P_{1,4}$ is not a component of the sextic and gives a reducible fiber of type $I_2$, whose components are $E_{P_{1,4}}$ and the pull back on $X_{(17,5,1)}$ of the strict transform on $\widetilde{\mathbb{P}^2}$ of the line $m$. This is easely proved blowing up $\mathbb P^2$ in the singular locus of the sextic and constructing explicitly $X_{(17,5,1)}$. Similarly the lines through $P_{5,6}$ and $P_{2,4}$ and through $P_{5,6}$ and $P_{3,4}$ correspond to fibers of type $I_2$.
The line $l_5$ is a line of the pencil and a component of the sextic and it passes through $P_{5,6}$ and other 4 double points of the sextic: $P_{1,5}, P_{2,5},P_{3,5},P_{4,5}$. It gives a reducible fiber of type $I_0^*$, whose components are $E_{P_{1,5}}$, $E_{P_{2,5}}$, $\mathcal{L}_5$, $E_{P_{3,5}}$, $E_{P_{4,5}}$. The same happens for the line $l_6$. 
The line $m'$ through $P_{5,6}$ and the triple point $A$ is not a component of the sextic and gives a reducible fiber of type $I_0^*$ whose components are $E_{A_1}$, $E_{A_2}$, $E_A$, $E_{A,3}$, $\mathcal{M}'$, where $\mathcal{M}'$ is the pull back on $X_{(17,5,1)}$ of the strict transform of $m'$ on $\widetilde{\mathbb{P}^2}$. Thus, the pencil of lines through $P_{5,6}$ gives an elliptic fibration on $X_{(17,5,1)}$ with reducible fibers $3I_0^*+3I_2$.\erem}\end{example}

\begin{example}\label{example: 17-5-1 not line}{\rm We consider again the surface $X_{(17,5,1)}$ with the same branch curve as in the previous example, but different pencil of rational curves, in order to obtain different reducible fibers. In particular for certain fibers we need to consider a pencil of rational curves of degree greater than 1.\\
Let us consider a pencil of conics through the four points $P_{1,4},P_{2,6},P_{3,4},P_{3,6}$. The conic of the pencil through $P_{4,6}$ splits in two lines $l_4\cup l_6$. It gives a fiber of type $I^*_2$ whose components are   $E_{P_{2,4}}$, $E_{P_{4,5}}$, $\mathcal{L}_4$, $E_{P_{4,6}}$, $\mathcal{L}_6$, $E_{P_{1,6}}$, $E_{P_{5,6}}$.\\ 
Let us consider the pencil of conics through $P_{1,4},P_{1,5},P_{3,5},P_{3,6}$. The conic of the pencil through $A$ splits in $l_1\cup l_3$ and corresponds to a fiber of type $III^*$, whose components are $E_{P_{1,6}}$, $\mathcal{L}_1$, $E_{A,1}$, $E_A$, $E_{A,2}$, $E_{A,3}$, $\mathcal{L}_3$, $E_{P_{3,4}}$.\\ 
Similarly, we can obtain fibers of type $II^*$ and $I^*_i$, $i=1,\ldots,14$ with a pencil of rational curves of degree greater than 1. \erem}\end{example}

By the previous examples it is clear that there is not only a choice of the pencil to obtain a certain type of fiber. Indeed the reducible fibers of type $I_0^*$ are obtained in two different ways in the example: as pull-back of a line of the pencil which is contained in the branch locus and as pull-back of a line of the pencil which is not in the branch locus and passes through a triple point. Moreover, it is clear that one can use different choices for the branch sextic associated to the surface $X_{(r,22-r,\delta)}$ (cf. Section \ref{section:double cover}) and this gives more possibilities in the choice of the pencil $\mathcal{P}$. In the following we use pencil of lines, when this is possible. Sometimes this implies a different choice for the branch sextic.

If $\mathcal{P}$ is a pencil of lines, a reducible fiber of the associated fibration corresponds to a line $l$ of the $\mathcal{P}$ through a singular point of the sextic and the type of the fiber is as follows:
\begin{enumerate}
	\item type $I_2$, if $l$ is not a component of the sextic and pass through a double point;
	\item type $I^*_0$, if $l$ is not a component of the sextic and pass through a triple point;
	\item type $I^*_0$, if $l$ is a component of the sextic and pass through 4 double points;
	\item type $I^*_2$, if $l$ is a component of the sextic and pass through a triple point;
	\item type $I^*_4$, if $l$ is a component of the sextic and pass through 2 triple points.
\end{enumerate}

\begin{rem} \rm{
When the reducible fibers are $I_2,I^*_0,I^*_2,I^*_4$, we will always use a pencil of lines, except for case 13) for $r=18$. In this case the reducible fibers are $2I^*_4$ and, for geometric reasons, there is no branch sextic associated to $r=18$ such that a pencil of lines gives this fibration.}
\end{rem}

Thanks to the construction of the elliptic fibration starting from a pencil of lines (and more in general from a pencil of rational curves) we can give an explicit equation of the elliptic fibration. We consider as an example the elliptic fibration on $X_{(13,9,1)}$ with reducible fibers $11I_2+2I_1$.

\begin{example}\label{ex: B13}{\rm The surface $X_{(13,9,1)}$ is obtained as double cover of $\PP^2$ branched along the sextic $B_{(13,9,1)}$ which consists of 3 conics. Up to a choice of coordinates of $\PP^2$ we can assume that
\begin{equation}\label{eq: B13} B_{(13,9,1)}=V((xy+axz+(-1-a)yz)(xy+bxz+(-1-b)yz)(cx^2+dxy+ey^2+fxz+gyz+z^2)),
\end{equation}which is a 7-dimensional family of reducible sextics. 

We consider the pencil of lines of $\mathbb{P}^2$ passing through $P=(0:0:1)$; the pencil $\mathcal{P}$ is $y=\tau x$ and the elliptic fibration is 
\begin{equation}\label{eq: E13,9,1}
w^2=(x +\frac{a}{\tau}z+(-1-a) z)(x +\frac{b}{\tau}z+(-1-b) z)(cx^2+d\tau x^2+e x^2\tau^2+fxz+g\tau x z+z^2)
\end{equation} 
obtained substituting $y=\tau x$ in \eqref{eq: B13} and considering the change of coordinate $w\mapsto w\tau x$.

Applying standard transformations to this equation one finds the following ``more canonical" equation $w^2=x(x^2+A(\tau)x+B(\tau))$ with: 
\begin{eqnarray*}\label{formula Weierstrass E13,9,1}
\begin{array}{ll}
A(\tau):=&2\left(c+d\tau+e\tau^2\right)\left(\tau+b\tau-b\right)\left(\tau+a\tau-a\right)+\tau\left(f+g\tau\right)\left(2\tau+b\tau-b+a\tau-a\right)+2\tau^2\\
B(\tau):=&\left(\left(c+d\tau+e\tau^2\right)\left(\tau+b\tau-b\right)^2+\tau\left(f+g\tau\right)\left(\tau+b\tau-b\right)+\tau^2\right) \\ &\left(\left(c+d\tau+e\tau^2\right)\left(\tau+a\tau-a\right)^2+\tau\left(f+g\tau\right)\left(\tau+a\tau-a\right)+\tau^2\right).\end{array}
\end{eqnarray*}
The values of $\tau$ such that $\Delta(\tau)=0$ correspond to the lines of the pencil through the singular points of the sextic and give the reducible fibers of the fibration. See \cite[Table IV.3.1]{Mi} for the type of fiber in relation to the vanishing of $A,B,\Delta$.\erem}\end{example}

\begin{rem}{\rm In \cite{GS3} it is shown that the surface $X_{(r,a,1)}$ specializes to the surface $X_{(r+1,a-1,1)}$, because the transcendental lattice of $X_{(r,a,1)}$ contains the transcendental lattice of $X_{(r+1,a-1,1)}$ as a primitive sublattice. Thus, a condition on parameters of the equation of $B_{(13,9,1)}$ specializes it to a sextic with associated invariants $(14,8,1)$. Specializing the sextic $B_{(13,9,1)}$, we automatically specialize the elliptic fibration, thus we find an elliptic fibration on $X_{(14,8,1)}$. The equation of this elliptic fibration can be obtained substituting the relations on parameters in the equation of the elliptic fibration on $X_{(13,9,1)}$. For example, putting $c=0$ in the Equation \ref{eq: B13} one obtains a sextic which splits in three conics with a triple point $(1:0:0)$. The double cover of $\mathbb{P}^2$ branched over this sextic is $X_{(14,8,1)}$ (cf.\ Proposition \ref{prop: double cover of P2}). Putting $c=0$ in Equation \eqref{eq: E13,9,1} one obtains an elliptic fibration with singular fibers $I_0^*+8I_2+4I_1$. It is an elliptic fibration on $X_{(14,8,1)}$ which specializes the elliptic fibration \eqref{eq: E13,9,1}. Similarly, specializations of the sextic $B_{(r,22-r,\delta)}$  (defined below) to $B_{(r+1,21-r,\delta)}$ are associated to specializations of elliptic fibrations.}
\end{rem}

The rest of this Section is devoted to the proof of the following Proposition
\begin{proposition}\label{prop: existence case b)}
The fibrations in tables of Section \ref{subsec: existence case b} appear as elliptic fibrations on $X_{(r,22-r,\delta)}$, in the case $m_2=n_2=0$. 
\end{proposition}
\begin{rem}\label{rem: more then one elliptic fibration for each line of the Tables}{\rm In the tables of Section \ref{subsec: existence case b} we allow both the fibers of type $I_2$ and of type $III$ as fibers with two components and both the fibers of type $I_1$ and of type $II$ as singular not reducible fibers. However, generically all the reducible fibers with two components are of type $I_2$ and all the singular not reducible fibers are of type $I_1$. In the following we give examples of elliptic fibrations on $X_{(r,22-r,\delta)}$ without fibers of type $III$ and $II$. Of course it is possible that, on the surface $X_{(r,22-r,\delta)},$ there exist elliptic fibrations with fibers of type $III$ and $II$. More in general, it is clear that in this section we give explicit examples of an elliptic fibration for each possibility listed in tables of Section \ref{subsec: existence case b}, 
but it is possible that there exist more then one (non--isomorphic) elliptic fibration for each line of these Tables. This surely happens in case on $X_{(18,4,0)}$, indeed in \cite{O} it is proved, for example, that there exist 9 non--isomorphic elliptic fibrations with reducible fibers $2I_2^*+4I_2$.\erem}\end{rem}

\subsection{$r=12$}
The sextic we consider is made up of a rational quartic $q$ with 3 double points $Q,R,T$ and a conic $c$. Up to projective transformation of $\mathbb{P}^2$ the equation of the branch sextic is \begin{equation}\label{eq: r=12}\begin{split}
B_{(12,10,1)}&:=(Ax^2y^2+Bx^2yz+Cxy^2z+Dx^2z^2+(1-4A-3B-2C-2D)xyz^2+\\&(A+B+D-1)y^2z^2+(2A+B+C-1)xz^3+yz^3)(Ex^2+xy+Fy^2+Gxz+Hyz)\end{split}\end{equation}
The pencil of lines through $P:=(0:0:1)\in q\cap c$ gives an elliptic fibration on $X_{(12,10,1)}$ and the reducible fibers are associated  to the lines of the pencil through $Q,R,T$ and through the 7 points in $\{q\cap c\} -\{P\}$, thus the reducible fibers are $10I_2$.\\

In the tables of the following section notations are the same as in tables of Section \ref{subsec: existence case b}, where the reader can find the list of reducible fibers of each fibration. 
To show the existence (and the equation) of the fibrations listed in tables of Section \ref{subsec: existence case b}, we consider pencils of rational curves through the singular points of the branch sextic $B_{(r,22-r,\delta)}$ described below.
In the fifth column of the table we list the sections. 

In order to use a pencil of lines, sometimes we need to change the branch sextic: in this case the sextic is obtained as specialization of $B_{(12,10,1)}$ or of $B_{(13,9,1)}$, imposing conditions on the parameters of the equations of these sextics (Equations \eqref{eq: B13}, \eqref{eq: r=12}). 
Thus the sextic we specialize is indicated in the third column, while in the fourth column we list the conditions on the parameters.
When the base point of the pencil is $(0:0:1)$ and the sextic is a specialization of $B_{(13,9,1)}$ then to obtain the equation of the elliptic fibration it suffices to substitute the conditions on the parameters in the equation of $\mathcal E_{(13,9,1)}$. In this case we denote by $t$ the 2-torsion section $t:\tau\mapsto (0,0;\tau)$.

\subsection{$r=13$}
In case 1) we consider the sextic obtained as specialization of the sextic $B_{(12,10,1)}$, while in case 2) we consider the sextic $B_{(13,9,1)}$ (see Example \ref{ex: B13}) made up of 3 conics without triple points, whose equation is \eqref{eq: B13}. 
\scriptsize{
\[\begin{tabular}{|l|l|l|l|}
\hline
 1)&$|H-E_{(0:0:1)}|$&$B_{(12,10,1)}$&$E=0$\\ \hline
 2)&$|H-E_{(0:0:1)}|$&$B_{(13,9,1)}$&\\ \hline
 \end{tabular}\]
}\normalsize

\subsection{$r=14$}
The branch sextics are obtained as specializations of the sextic $B_{(12,10,1)}$ and $B_{(13,9,1)}$. 
\scriptsize{
\[\begin{tabular}{|l|l|l|l|}
\hline
1)&$|H-E_{(0:0:1)}|$&$B_{(12,10,1)}$&$E=G=0$\\ \hline
2)&$|H-E_{(0:0:1)}|$&$B_{(12,10,1)}$&$E=F=0$\\ \hline
3)&$|H-E_{(0:0:1)}|$&$B_{(13,9,1)}$&$c=0$\\ \hline	
4)&$|H-E_{(1:0:0)}|$&$B_{(13,9,1)}$&$c=0$\\
\hline
\end{tabular}\]
}\normalsize

\subsection{$r=15$}
The sextic $B_{(15,7,1)}$ consists of two lines $l_1,l_2$ and two conics $c_3,c_4$ with a triple point $A:=l_1\cap c_3\cap c_4$. The 10 double points are $P_{1,2}$, $R_{1,j}=l_1\cap c_j$, $j=3,4$, $\{R^1_{2,j}, R^2_{2,j}\}=l_2\cap c_j$, $j=3,4$,  $\{Q^1_{3,4}, Q^2_{3,4},Q^3_{3,4}\}\in c_3\cap c_4$.
\scriptsize{
\[\begin{tabular}{|l|l|l|l|l|}
\hline
1)&$|3H-2E_{R_{1,3}}-E_{R_{1,4}}-E_{R^1_{2,3}}-E_{R^2_{2,3}}-E_{Q^1_{3,4}}-E_{Q^2_{3,4}}|$&$B_{(15,7,1)}$&&$\mathcal L_2$\\ \hline
2)&$|H-E_{(C,-B:0)}|$&$B_{(12,10,1)}$&$A=F=E=0$&\\ \hline
3)&$|H-E_{(0:0:1)}|$&$B_{(12,10,1)}$&$E=F=G=0$&\\ \hline
4a)&$|H-E_{(0:0:1)}|$& $B_{(13,9,1)}$&$b=0, f=-1-c$&\\ \hline
4b)&$|H-E_{(1:0:c)}|$&$B_{(13,9,1)}$&$b=0, f=-1-c$&\\ \hline
5)&$|H-E_{(0:0:1)}|$&$B_{(12,10,1)}$&$E=F=0,$&\\&&&$ H=-1-G$&\\ \hline
6)&$|H-E_{(0:0:1)}|$&$B_{(13,9,1)}$&$a=-1, b=0$&\\ \hline
7)&$|H-E_{(0:1:0)}|$&$B_{(13,9,1)}$&$c=e=0$&\\ \hline	
\end{tabular}\]
}\normalsize

In case 1) we consider a pencil of rational cubics that pass through $R_{1,4},R^1_{2,3},R^2_{2,3},$ $Q^1_{3,4},Q^2_{3,4}$ with a double point at $R_{1,3}$. The only section of the fibration is $\mathcal{L}_2$ and the fiber of type $III^*$ corresponds to the cubic that splits in $l_1\cup c_3$. 

\subsection{$r=16$}
We consider the sextic $B_{(16,6,1)}$ made up of four lines $l_1,\ldots,l_4$ and a conic $c_5$ with a triple points $A=l_1\cap l_2 \cap c_5$. The double points are $P_{h,k}$, $h=1,2$, $k=3,4$, $P_{3,4}$, $R_{1,5}$, $R_{2,5}$ and $\{R^1_{j,5},R^2_{j,5}\}$, $j=3,4$. 
\scriptsize{
\[\begin{tabular}{|l|l|l|l|l|}
\hline
1)&$|4H-2E_{P_{1,3}}-2E_{P_{2,4}}-2E_{R_{2,5}}-E_{R_{1,5}}$& $B_{(16,6,1)}$&&$\mathcal L_4$\\ &$-E_{P_{1,4}}-E_{R^1_{3,5}}-E_{R^2_{3,5}}|$&&&\\ \hline 
2)&$|2H-E_{P_{1,4}}-E_{R_{1,5}}-E_{P_{2,4}}-E_{P_{2,3}}|$&$B_{(16,6,1)}$&&$\mathcal L_3$ \\ \hline
3)&$|4H-2E_{P_{2,3}}-2E_{R^2_{3,5}}-2E_{P_{1,4}}-E_{R^2_{4,5}}-E_{R_{1,5}}$&$B_{(16,6,1)}$&&$\mathcal L_2, E_A$\\ &$-2E_A-E_{A,1}-E_{A,2}-E_{A,5}-E_{R^1_{4,5}}|$&&&\\ \hline
4)&$|3H-2E_{P_{2,3}}-E_{P_{1,3}}-E_{R^1_{4,5}}-E_{R^2_{4,5}}-E_{R_{2,5}}$& $B_{(16,6,1)}$&&$\mathcal L_1$\\ &$-E_{P_{1,4}}|$&&&\\ \hline
5)&$|H-E_{(0:1:1-A-B-D)}|$&$B_{(12,10,1)}$&$E=F=G=H=0$&\\ \hline
6a)&$|H-E_{(0:0:1)}|$& $B_{(13,9,1)}$&$c=b=0, f=-1$&\\ \hline
6b)&$|H-E_{(0:0:1)}|$&$B_{(12,10,1)}$&$A=0, E=F=-1/2, $&\\ &&&$H=-G$&\\\hline
7)&$|H-E_{(0:-G,1)}|$&$B_{(12,10,1)}$&$A=E=F=H=0$&\\ \hline
8)&$|H-E_{(0:0:1)}|$&$B_{(12,10,1)}$&$E=F=G=0,H=-1$&\\ \hline
9a)&$|H-E_{(0:0:1)}|$ &$B_{(13,9,1)}$&$a=c=e=0$&\\ \hline
9b)&$|H-E_{(1:0:c)}|$&$B_{(13,9,1)}$&$b=0,f=-1-c,$&\\ &&&$g=-d-e$&\\ \hline
10)&$|H-E_{(0:1:-g)}|$ &$B_{(13,9,1)}$&$a=e=0,b=-1$&\\ \hline
11)&$|H-E_{(0:0:1)}|$ &$B_{(13,9,1)}$&$c=e=0, g=-d-f-1$&\\ \hline
12)& $|H-E_{(1:0:0)}|$&$B_{(13,9,1)}$&$c=e=0,g=-d-f-1$&\\ \hline
\end{tabular}\]
}\normalsize
The quartics of case 1) pass through $R_{1,5},P_{1,4},R^1_{3,5},R^2_{3,5}$ and have double points in $P_{1,3},P_{2,4},R_{2,5}$, so they are rational; the quartic $l_1\cup l_2\cup l_2\cup l_3$ is a quartic of the pencil and corresponds to the fiber of type $II^*$.
The pencil of conics in case 2) consists of conics passing through $P_{1,4},R_{1,5},P_{2,4},P_{2,3}$. The fiber of type $III^*$ comes from the reducible conic $l_1\cup l_2$, while the fiber of type $I^*_0$ corresponds to the conic $l_4\cup m$ where $m$ is the line through $R_{1,5}$ and $P_{2,3}$.
The pencils described in the other cases are similar to the ones already shown.

In cases 1)--4) the branch sextic is the union of 4 lines and 1 conic (cf.\ Example \ref{ex: branch locus X16,6,1}, case $vi)$). In cases 5)--12) we use a different sextics in order to consider pencil of lines: in cases 6a), 9a) and 9b) the sextic is the union of 2 conics and 2 lines with 2 triple points (cf.\ Example \ref{ex: branch locus X16,6,1}, case $v)$); in case 10) the sextic is the union of 1 conic and 4 lines with 1 triple point (cf.\ Example \ref{ex: branch locus X16,6,1}, case $vi)$); in cases 11) and 12) the sextic is the union of 3 conics with 3 triple points (cf.\ Example \ref{ex: branch locus X16,6,1}, case $iv)$); in cases 5) and 8) the sextic is the union of 1 rational quartic and 2 lines with 3 triple points (cf.\ Example \ref{ex: branch locus X16,6,1}, case $i)$); in cases 6b) and 7) the sextic is the union of 1 rational cubic and 3 lines with 2 triple points (cf.\ Example \ref{ex: branch locus X16,6,1}, case $iii)$).
\begin{rem}{\rm In \cite{CD:vGS involutions} the elliptic fibration 6b) is associated to a certain pencil of rational quintics in $\mathbb{P}^2$ and the branch sextic is the union of six lines in general position. Here we consider the union of a rational cubic and three lines as branch sextic and thanks to this choice we obtain the same elliptic fibration from a pencil of lines.\erem}\end{rem}
\subsection{$r=17$}
We consider the sextic $B_{(17,5,1)}$ which consists of six lines $l_1,\ldots, l_6$ with a triple point $A:=l_1\cap l_2\cap l_3$ and 12 double points $P_{1,4}$, $P_{1,5}$, $P_{1,6}$, $P_{2,4}$, $P_{2,5}$, $P_{2,6}$, $P_{3,4}$, $P_{3,5}$, $P_{3,6}$, $P_{4,5}$,$P_{4,6}$, $P_{5,6}$. 
\scriptsize{
\[\begin{tabular}{|l|l|l|l|l|}
\hline
1)&$|6H-3E_{P_{3,5}}-3E_{P_{2,4}}-3E_{P_{4,6}}-2E_{P_{3,6}}-2E_{A}$&$B_{(17,5,1)}$&&$\mathcal L_2$\\ 
&$-2E_{A,1}-E_{A,2}-E_{A,3}-E_{P_{1,5}}-E_{P_{2,5}}-E_{P_{5,6}}|$&&&\\ \hline
2)&$|3H-2E_{A}-2E_{A,2}-E_{A,1}-E_{A,3}-2E_{P_{1,4}}-E_{P_{2,5}}$&$B_{(17,5,1)}$&&$\mathcal L_3$\\ &$-E_{P_{3,5}}-E_{P_{4,5}}|$&&&\\ \hline
3)&$|4H-2E_{A}-E_{A,1}-E_{A,2}-E_{A,3}-2E_{P_{1,5}}-2E_{P_{2,5}}$&$B_{(17,5,1)}$&&$E_A$\\ 
&$-2E_{P_{3,6}}-E_{P_{3,4}}-E_{P_{1,6}}-E_{P_{2,6}}|$&&&\\ \hline
4)&$|2H-E_{P_{1,4}}-E_{P_{1,5}}-E_{P_{3,6}}-E_{P_{3,5}}|$&$B_{(17,5,1)}$&&$\mathcal L_4,\mathcal L_6$\\ \hline
5)&$|5H-3E_{P_{1,4}}-2E_{P_{1,5}}-2E_{P_{3,6}}-2E_{P_{5,6}}-E_{P_{2,6}}$&$B_{(17,5,1)}$&&$\mathcal L_3$\\ &$-E_{P_{2,4}}-E_{P_{3,4}}-E_{P_{3,5}}|$&&&\\ \hline
6)&$|4H-2E_{P_{3,5}}-2E_{P_{2,4}}-2E_{P_{3,6}}-E_{P_{4,5}}-E_{P_{1,4}}$&$B_{(17,5,1)}$&&$\mathcal L_2$\\ &$-E_{P_{5,6}}-E_{P_{2,6}}|$&&&\\ \hline
7)&$|3H-2E_{P_{3,4}}-E_{P_{3,5}}-E_{P_{5,6}}-E_{P_{2,6}}-E_{P_{1,6}}$&$B_{(17,5,1)}$&&$\mathcal L_1, \mathcal L_5$\\ &$-E_{P_{1,4}}|$&&&\\\hline
8)&$|H-E_{(0:1:1)}|$&$B_{(12,10,1)}$&$A=0,E=0, G=-1,H=F=0$&\\ \hline
9)&$|H-E_{(0:0:1)}|$&$B_{(13,9,1)}$&$c=b=0,a=-1,f=-1$&\\ \hline
10)&$|H-E_{(1:1:1)}|$&$B_{(13,9,1)}$&$a=e=0,f=-1-c, d=-cg$&\\ \hline
11)&$|H-E_{(0:0:1)}|$&$B_{(13,9,1)}$&$c=e=a=0, g=-d-f-1$&\\ \hline
12)&$|H-E_{(0:0:1)}|$&$B_{(13,9,1)}$& $c=e=a=d=0$&\\ \hline
13)&$|H-E_{(0:0:1)}|$&$B_{(13,9,1)}$&$c=e=0,f=-1,g=-d$&\\ \hline
14a)&$|H-E_{(0:0:1)}|$&$B_{(13,9,1)}$&$c=b=0,a=-1,g=-1-e$&\\ \hline
14b)&$|H-E_{(1:0:0)}|$&$B_{(12,10,1)}$&$A=0, E=H=F=0, G=-1$&\\ \hline
\end{tabular}\]
}\normalsize
In case 1) we consider a pencil of rational sextics through $A,P_{1,5},P_{2,5},P_{5,6}$, with a double point in $P_{3,6}$ and triple points in $P_{3,5},P_{2,4},P_{4,6}$ and such that the tangent line to the sextic at $A$ is $l_1$. The reducible sextic $l_5\cup l_3\cup l_3\cup l_4\cup l_4\cup l_4$ corresponds to the singular fiber of type $II^*$.
The quintics of the pencil of case 5) pass through $P_{2,6},P_{2,4},P_{3,4},P_{3,5}$, have double points in $P_{1,5},P_{5,6},P_{3,6}$ and a triple point in $P_{1,4}$, so they are rational; the fiber of type $I^*_8$ is given by the quintic that splits in $l_1\cup l_4\cup l_5\cup l_6\cup m$, where $m$ is the line through $P_{1,4}$ and $P_{3,6}$.
The other cases are similar to the ones analyzed for lower values of $r$.

\subsection{$r=18$}\label{subsec: existence a 18}
We consider the two different sextics $B_{(18,4,1)}$, described in Section \ref{section: non trivial MW r=18}, and $B_{(18,4,0)}$ which consists of the lines $l_i$, $i=1,\ldots, 6$, with triple points $A:=l_1\cap l_2\cap l_3$, $B:=l_4\cap l_5\cap l_6$ and double points $P_{1,4}$, $P_{1,5}$, $P_{1,6}$, $P_{2,4}$, $P_{2,5}$, $P_{2,6}$, $P_{3,4}$, $P_{3,5}$, $P_{3,6}$. The values of $\delta$ for each fibration is indicated in the sixth column. The pencils are similar to the ones described for lower values of $r$.
\scriptsize{
\[\begin{tabular}{|l|l|l|l|l|l|}
\hline
1)&$|3H-2E_{A}-2E_{A,3}-E_{A,2}-E_{A,1}-2E_{P_{1,4}}$&$B_{(18,4,1)}$&&$\mathcal L_1$&1\\ &$-E_{P_{3,5}}-E_{P_{2,5}}-E_{P_{4,5}}|$&&&&\\ \hline	
2)&$|4H-2E_{B}-E_{B,4}-E_{B,5}-E_{B,6}-2E_{P_{3,4}}$&$B_{(18,4,0)}$&&$E_B$ &0\\ &$-2E_{P_{1,5}}-2E_{P_{1,6}}-E_{P_{2,4}}-E_{P_{3,5}}-E_{P_{3,6}}|$&&&&\\ \hline
3)&$|3H-2E_{A}-2E_{A,3}-E_{A,2}-E_{A,1}-2E_{B}$&$B_{(18,4,1)}$&&$\mathcal L_1$&1\\ &$-2E_{B,5}-E_{B,1}-E_{B,6}-2E_{P_{2,6}}-E_{P_{3,5}}|$&&&&\\ \hline
4)&$|3H-2E_{A}-2E_{A,3}-E_{A,2}-E_{A,1}-2E_{B}$&$B_{(18,4,1)}$&&$\mathcal L_6$&1\\ &$-2E_{B,1}-E_{B,6}-E_{B,5}-2E_{P_{2,5}}-E_{P_{3,6}}|$&&&&\\ \hline
5a)&$|2H-E_{P_{3,6}}-E_{P_{3,5}}-E_{P_{2,4}}-E_{P_{2,6}}|$&$B_{(18,4,1)}$&&$\mathcal L_4,\mathcal L_5$&1\\ \hline
5b)&$|2H-E_{P_{1,4}}-E_{P_{1,5}}-E_{P_{3,6}}-E_{P_{3,4}}|$&$B_{(18,4,0)}$&&$\mathcal L_5,\mathcal L_6$&0\\ \hline
6)&$|2H-E_{P_{3,6}}-E_{P_{4,6}}-E_{P_{4,5}}-E_{P_{2,5}}|$&$B_{(18,4,1)}$ &&$\mathcal L_2,\mathcal L_3$&1\\ \hline
7)&$|5H-2E_{B}-2E_{B,1}-E_{B,5}-E_{B,6}-3E_{P_{3,6}}$ &$B_{(18,4,1)}$&&$\mathcal L_1$&1\\ &$-2E_{P_{4,5}}-2E_{P_{1,4}}-2E_{P_{3,5}}-E_{P_{2,4}}-E_{P_{2,6}}|$&&&&\\ \hline
8)&$|4H-2E_{A}-2E_{A,2}-E_{A,1}-E_{A,3}-2E_{P_{3,6}}$&$B_{(18,4,0)}$ &&$\mathcal L_5$&0\\ &$-2E_{P_{1,5}}-2E_{P_{2,6}}-E_{P_{1,4}}-E_{P_{3,5}}|$&&&&\\ \hline
9)&$|3H-2E_{P_{2,5}}-E_{P_{1,4}}-E_{P_{4,6}}-E_{P_{2,6}}$&$B_{(18,4,1)}$& &$\mathcal L_3, \mathcal L_6$&1\\ &$-E_{P_{3,4}}-E_{P_{3,5}}|$&&&&\\ \hline
10)&$|3H-2E_{B}-2E_{B,6}-E_{B,1}-E_{B,5}-2E_{P_{1,4}}$&$B_{(18,4,1)}$&&$\mathcal L_3$&1\\ &$-E_{P_{3,5}}-E_{P_{4,5}}-E_{P_{3,6}}|$&&&&\\ \hline
11)&$|2H-E_{P_{3,6}}-E_{P_{3,4}}-E_{P_{4,5}}-E_{P_{2,5}}|$&$B_{(18,4,1)}$&&$\mathcal L_2,\mathcal L_6$&1\\ \hline
12)&$|2H-E_{P_{3,6}}-E_{P_{3,4}}-E_{P_{2,5}}-E_{P_{1,5}}|$&$B_{(18,4,0)}$&&$\mathcal L_1,\mathcal L_2,$&0\\ &&&&$\mathcal L_4, \mathcal L_6$&\\ \hline
13)&$|3H-2E_{A}-2E_{A,1}-E_{A,2}-E_{A,3}-2E_{P_{2,6}}$&$B_{(18,4,0)}$&&$\mathcal L_5$&0\\ &$-E_{P_{3,5}}-E_{P_{3,6}}-E_{P_{1,5}}|$&&&&\\ \hline
14a)&$|H-E_{(0:0:1)}|$&$B_{(13,9,1)}$&$c=b=0,a=f=-1,g=-1-e$&&1\\ \hline
14b)&$|H-E_{(g:1:0)}|$&$B_{(13,9,1)}$&$c=e=d=0,f=-1,a=0$&&1\\ \hline
15)&$|H-E_{(0:0:1)}|$&$B_{(13,9,1)}$&$a=0,b=-1,c=0,f=-1,g=-d-e$&&0\\ \hline
16)&$|H-E_{(1:0:1)}|$&$B_{(13,9,1)}$&$a=0,c=d=e=0,f=-1$&&1 \\ \hline
17)&$|H-E_{(0:0:1)}|$&$B_{(13,9,1)}$&$a=0,b=-1,d=1,f=-1-c,g=-1-e$&&1\\ \hline
18)&$|H-E_{(0:0:1)}|$&$B_{(13,9,1)}$&$c=e=0,g=-d,f=-1,a=-1$ &&1\\ \hline
19)&$|H-E_{(1:0:0)}|$&$B_{(13,9,1)}$&$a=0,b=-1,c=e=0,d=fg$&&0\\ \hline
20a)&$|H-E_{(0:0:1)}|$&$B_{(13,9,1)}$&$a=c=e=0,b=-1,d=fg$&&0\\ \hline 
20b)&$|H-E_{(0:0:1)}|$&$B_{(13,9,1)}$&$a=c=d=e=0,b=-1$&&1\\ \hline 
\end{tabular}\]
}\normalsize

\subsection{$r=19$}\label{subsec: existence r=19 tirvial MW} 
The sextic $B_{(19,3,1)}$ consists of six lines $l_1,\ldots, l_6$ with 3 triple points $A:=l_1\cap l_2 \cap l_3,B:=l_1\cap l_5\cap l_6,C:=l_3\cap l_4 \cap l_5$ and 6 double points: $P_{1,4},P_{2,5},P_{3,6},P_{2,4},P_{2,6},P_{4,6}$.
\scriptsize{
\[\begin{tabular}{|l|l|l|l|l|}
\hline
1)&$|3H-2E_A-2E_{A,2}-E_{A,1}-E_{A,3}-2E_{B}-2E_{B,5}-E_{B,1}$&$B_{(19,3,1)}$&&$\mathcal L_1$ \\&$-E_{B,6}-2E_{P_{3,6}}-E_{P_{2,5}}|$&&&\\ \hline	
2)&$|3H-2E_{A}-2E_{A,1}-E_{A,2}-E_{A,3}-2E_{B}-2E_{B,5}-E_{B,1}$ &$B_{(19,3,1)}$&& $\mathcal L_2$\\ &$-E_{B,6}-2E_{P_{3,6}}-E_{P_{2,5}}|$&&&\\ \hline
3)&$|2H-2E_{A}-E_{A,1}-E_{A,2}-E_{A,3}-E_{P_{1,4}}-E_{P_{4,6}}-E_{P_{3,6}}|$&$B_{(19,3,1)}$&&$\mathcal L_2,E_A$\\ \hline
4)&$|3H-2E_{B}-2E_{B,5}-E_{B,1}-E_{B,6}-2E_{C}-2E_{C,4}-E_{C,3}$&$B_{(19,3,1)}$&&$\mathcal L_1$\\ &$-E_{C,5}-2E_{P_{3,6}}-E_{P_{1,4}}|$&&&\\ \hline
5)&$|2H-2E_{A}-E_{A,1}-E_{A,2}-E_{A,3}-E_{P_{1,4}}-E_{P_{2,6}}-E_{P_{4,6}}|$&$B_{(19,3,1)}$&&$\mathcal L_3,E_A$ \\ \hline
6)&$|2H-2E_{A}-E_{A,1}-E_{A,2}-E_{A,3}-2E_{B}-E_{B,1}-E_{B,5}$&$B_{(19,3,1)}$&&$E_A, E_B$\\ &$-E_{B,6}-E_{P_{3,6}}-E_{P_{2,5}}|$&&&\\ \hline
7)&$|5H-4E_{C}-4E_{C,5}-2E_{C,4}-2E_{C,3}-2E_{B}-2E_{B,1}-E_{B,5}$&$B_{(19,3,1)}$&&$\mathcal L_1$\\
&$-E_{B,6}-3E_{P_{3,6}}-2E_{P_{1,4}}-E_{P_{2,6}}-E_{P_{2,4}}|$&&&\\ \hline
8)&$|3H-2E_{A}-2E_{A,3}-E_{A,1}-E_{A,2}-2E_{B}-2E_{B,6}-E_{B,1}$&$B_{(19,3,1)}$&&$\mathcal L_3,\mathcal L_6$\\ &$-E_{B,5}-2E_{P_{2,5}}-E_{P_{1,4}}|$&&&\\ \hline
9)&$|3H-2E_{B}-2E_{B,5}-E_{B,1}-E_{B,6}-2E_{C}-2E_{C,3}-E_{C,4}$&$B_{(19,3,1)}$&&$\mathcal L_3,\mathcal L_6$\\ &$-E_{C,5}-2E_{P_{2,4}}-E_{P_{2,6}}|$&&&\\ \hline
10a)&$|2H-E_{P_{1,4}}-E_{P_{2,5}}-E_{P_{2,6}}-E_{P_{4,6}}|$&$B_{(19,3,1)}$&&$\mathcal L_1,\mathcal L_5$ \\ \hline
10b)&$|2H-2E_B-2B_6-B_1-B_5-2C-C_3-C_4-C_5-E_{P_{1,4}}|$&$B_{(19,3,1)}$&&$\mathcal L_3,E_C$\\
\hline
11)&$|H-E_{(0:0:1)}|$&$B_{(13,9,1)}$&$b=-1,a=c=e=0,$&\\&&&$f=-1,g=-1$& \\ \hline
12)&$|H-E_{(0:0:1)}|$&$B_{(13,9,1)}$&$a=c=d=e=0,b=g=-1$&\\ \hline
13)&$|H-E_{(0:0:1)}|$&$B_{(13,9,1)}$&$a=c=d=e=0,b=-1,$&\\&&&$g=-1-f$&\\ \hline
\end{tabular}\]
}\normalsize

\subsection{$r=20$}
We consider the sextic $B_{(20,2,1)}$ made up of 6 lines $l_1,\ldots, l_6$ with four triple points: $A:=l_1\cap l_2 \cap l_3,B:=l_1\cap l_5\cap l_6,C:=l_3\cap l_4 \cap l_5, D:=l_2\cap l_4\cap l_6$, as in section \ref{section: non trivial MW r=20}.
\scriptsize{
\[\begin{tabular}{|l|l|l|l|l|}
\hline
1)& $|3H-2E_D-2E_{D,2}-E_{D,4}-E_{D,6}-2E_B-2E_{B,5}$& $B_{(20,2,1)}$&&$\mathcal L_6$\\&$-E_{B,1}-E_{B,6}-2E_{P_{1,4}}-E_{P_{2,5}}|$&&&\\ \hline	
2)& $|3H-2E_B-2E_{B,5}-E_{B,1}-E_{B,6}-2E_C-2E_{C,4}$& $B_{(20,2,1)}$&&$\mathcal L_1$\\ &$-E_{C,3}-E_{C,5}-E_{P_{1,4}}-2E_{P_{3,6}}|$&&&\\ \hline
3)&$|2H-2E_A-E_{A,1}-E_{A,2}-E_{A,3}-2E_{B}-E_{B,1}$& $B_{(20,2,1)}$&&$E_A, E_B$\\ &$-E_{B,5}-E_{B,6}-E_{P_{2,5}}-E_{P_{3,6}}|$&&&\\ \hline
4)&$|2H-2E_{C}-2E_{C,3}-E_{C,4}-E_{C,5}-2E_{D}-2E_{D,6}$ & $B_{(20,2,1)}$&&$\mathcal L_2,\mathcal L_5$\\ &$-E_{D,4}-E_{D,2}|$&&&\\ \hline
5)&$|5H-2E_{B}-2E_{B,1}-E_{B,5}-E_{B,6}-4E_{C}-4E_{C,5}-2E_{C,3}$&$B_{(20,2,1)}$&&$\mathcal L_1$\\
&$-2E_{C,4}-2E_{D}-2E_{D,2}-E_{D,4}-E_{D,6}-2E_{P_{1,4}}-3E_{P_{3,6}}|$&&& \\ \hline
6) &$|3H-2E_{C}-2E_{C,5}-E_{C,4}-E_{C,3}-2E_{D}-2E_{D,2}$ & $B_{(20,2,1)}$&&$\mathcal L_2,\mathcal L_5$\\ &$-E_{D,4}-E_{D,6}-E_{P_{1,4}}-2E_{P_{3,6}}|$&&&\\ \hline
7)& $|2H-2E_{C}-2E_{C,5}-E_{C,4}-E_{C,3}-E_{P_{1,4}}-E_{P_{3,6}}|$& $B_{(20,2,1)}$&&$\mathcal L_1,\mathcal L_6$\\ \hline
8)&$|H-E_{(0:0:1)}|$&$B_{(13,9,1)}$&$a=c=d=e=0,$&\\&&&$b=g=-1,f=-1$&\\ \hline
9)&$|H-E_{(1:0:0)}|$&$B_{(13,9,1)}$&$a=c=d=e=0,$&\\&&&$b=g=-1,f=-1$& \\ \hline
\end{tabular}\]
}\normalsize

\section{Van Geemen--Sarti involutions and isogenies}\label{section: van Geemen--Sarti involutions}
In \cite[Theorem 3.1]{GS3} it is proved that all the K3 surfaces of type $X_{(r,22-r,\delta)}$ admit at least a symplectic involution. It is clear that it is not necessarily a van Geemen--Sarti involution. The aim of this section is to classify all the van Geemen--Sarti involutions of the K3 surfaces $X_{(r,22-r,\delta)}$. Of course this is equivalent to classifying the sections of order 2 on the elliptic fibrations on such K3 surfaces. Once again we have to distinguish two cases: the one where the involution $\iota$ on $X_{(r,22-r,\delta)}$ is of type a) (cf.\  Example \ref{ex: involution on the basis} and Section \ref{section:existence case 1}) and the one where the involution $\iota$ on $X_{(r,22-r,\delta)}$ is of type b) (cf.\  Example \ref{ex: hyperelliptic involution} and Section  \ref{sec: existence case 2)}).\\ 
We also show that certain van Geemen--Sarti involutions on $X_{(r,22-r,\delta)}$ are Morrison--Nikulin involutions. The Morrison--Nikulin involutions are symplectic involution $\mu$ defined on a K3 surface $S$  such that the desingularization, $\widetilde{S/\mu}$, of $S/\mu$ is a Kummer surface $Km(A)$ (for a certain Abelian surface $A$) and $\T_{S}\simeq \T_A$. Equivalent conditions to be a Morrison--Nikulin involutions are: $\widetilde{S/\mu}$ is a Kummer surface $Km(A)$ and $\T_{\widetilde{S/\mu}}\simeq \T_A(2)$ or $\mu$ acts on the N\'eron--Severi group of $S$ switching two orthogonal copies of the lattice $E_8(-1)$. Morrison--Nikulin involutions induced by translation by section of order 2 (i.e.\ Morrison--Nikulin involutions which are also van Geemen--Sarti involutions) are described for the first time in \cite{vGS}. In \cite{Ko} Morrison--Nikulin involutions which are van Geemen--Sarti involutions acting on K3 surfaces of Picard number 17 with an elliptic fibration with a finite Mordell--Weil group are classified. In \cite{Sc} three 3-dimensional families of K3 surfaces 
with a Morrison--Nikulin involution are described. In each of these families the Morrison--Nikulin involution are van Geemen--Sarti involutions. We observe that for these families the Mordell--Weil group is not finite.
\subsection{van Geemen--Sarti involutions of type a)}
The elliptic fibrations on $X_{(r,22-r,\delta)}$ with non--symplectic involution $\iota$ of type a) are listed in Proposition \ref{prop: classification non finite MW}. Until now we did not describe the torsion part of the Mordell--Weil group, but now we are interested in the elliptic fibrations admitting a section of order 2, thus we have to consider exactly the torsion part of the Mordell--Weil group. 
\begin{proposition}\label{prop: van Geemen Sarti case a) and Morrison Nikulin} Among the elliptic fibrations of Table \eqref{table: possible fiber non trivial MW} the only ones admitting a section of order 2 are: $r=20$ reducible fibers $I_{4}+I_{16}$; $r=19$, reducible fibers $I_2+I_{16}$; $r=18$ reducible fibers either $I_4+I_{12}$ or $2I_8$.\\
In cases $r=20$ and $r=19$ the van Geemen--Sarti involution is also a Morrison--Nikulin involution, but this is not the case when $r=18$ holds.\end{proposition} 
\bprf A section of order 2 meets a fiber of type $I_{2k}$ either in the component $\Theta_0$ or in the component $\Theta_k$. Let us analyze elliptic fibrations with two reducible fibers $I_{2k_1}+I_{2k_2}$ with a section $t$ of order 2. The torsion section has to meet at least one fiber in a component which is not $\Theta_0$ (otherwise the section is not of order 2). If it meets one of the two reducible fibers (say the fiber $I_{2k_1}$) in the component $\Theta_0$, then the height formula implies that $k_2=8$. Thus, one of the two reducible fibers is $I_{16}$. Among the admissible elliptic fibrations of Table \eqref{table: possible fiber non trivial MW}, the only possibilities are $r=20$ with reducible fibers $I_4+I_{16}$ and $r=19$ with reducible fibers $I_2+I_{16}$. If the section of order 2 meets both the fibers $I_{2k_1}+I_{2k_2}$ in $\Theta_i$ with $i\neq 0$, then the height formula implies that $k_1+k_2=8$. By Table \eqref{table: possible fiber non trivial MW} this implies that $r=18$. There are no sections of order 2 on the elliptic fibration with a fiber of type $IV^*$ (one can use again the height formula, or consider the group law on the fibers of such a type).\\
If $t$ is a section of order 2 of an elliptic fibration $\mathcal{E}_r$ on $X_{(r,22-r,\delta)}$, then $t$ induces a 2-torsion section $t'$ on the elliptic fibration $\mathcal{E}'_r$ induced by $\mathcal{E}_r$ on the rational surface $X_{(r,22-r,\delta)}/\iota$. The elliptic fibrations $\mathcal{E}'_r$ are listed in Table \eqref{table: existence fibration non trivial MW}: the elliptic fibration with reducible fibers $I_{2k_1}+I_{2k_2}$ on $X_{(r,22-r,k)}$ induces on $X_{(r,22-r,\delta)}/\iota$ an elliptic fibration with reducible fibers $I_{k_1}+I_{k_2}$. If either both $k_1$ and $k_2$ are odd or one of $k_i$, $i=1,2$ is odd and the other is not 8, then the elliptic fibration with reducible fibers $I_{k_1}+I_{k_2}$ does not admit a section of order 2 (because of the height formula and the fact that a section of order 2 meets a fiber of type $I_{k}$ in the component $\Theta_0$ if $k$ is odd). Thus, on the rational elliptic fibrations with reducible fibers $I_1+I_7$ or $I_3+I_5$  there is no section of order 2, and hence on the elliptic K3 surface with $r=18$ and reducible fibers $I_2+I_{14}$, $I_{6}+I_{10}$ there is no section of order 2.\\
We proved that the only elliptic fibrations among the ones listed in Table \eqref{table: possible fiber non trivial MW} which can admit a section of order 2 are the ones listed in the statement. To prove that these elliptic fibrations admit in fact a section of order 2, it suffices to show that there exists a rational curve which intersects the reducible fibers as prescribed for a section of order 2 by the height formula. For example, if we consider the case $r=20$ and reducible fibers $I_4+I_{16}$, the section of order 2 has to meet the fiber of type $I_4$ in $\Theta_0$ and the fiber of type $I_{16}$ in $\Theta_8$. In the proof of Proposition  \ref{prop: existence case 1} we described the fiber of type $I_{16}$ of the fibration with fibers $I_4+I_{16}$ and we choose the curve $E_{P_{2,5}}$ as zero section. With respect to this choice the curve $E_{A,2}$ is a section of order 2. Analogously, in case $r=19$ the curve $E_{A,2}$ is a section of order 2. In case $r=18$, $\delta=0$ the presence of a section of order 2 is proved in \cite{O}.\\
If $t$ is a section of order 2 meeting a fiber of type $I_{16}$ in the component $\Theta_8$, the translation by this torsion section is a Morrison--Nikulin involution (cf. \cite{vGS}), thus in cases $r=20$ and $r=19$ the van Geemen--Sarti involution is also a Morrison--Nikulin involution.\\
On the surface $X_{(18,4,\delta)}$, $\delta=0,1$ there is no Morrison--Nikulin involution. Indeed  if the rank of the Picard number of a K3 surface $S$ is 18 then $S$ admits a Morrison--Nikulin involution if and only if the transcendental lattice $\T_S$ is isometric to $U\oplus T'$ for a certain lattice $T'$ of rank 2 (cf. \cite[Corollary 6.4]{Mo}). The length of a lattice isometric to $U\oplus T'$ is at most 2 (i.e. at most the maximal length of $T'$), but the length of the transcendental lattice of $X_{(18,4,\delta)}$ is 4. \erem

\begin{rem}{\rm If $r=18$ and $\delta=1$ holds then there exists at least one elliptic fibration with reducible fibers $I_4+I_{12}$ (resp.\ $2I_8$), which is the one described in the proof of Proposition \ref{prop: existence case 1}. This elliptic fibration does not admit a section of order 2. Indeed a $\Z$-basis of this K3 surface is made up of $\mathcal{L}_4$, $E_{P_{2,4}}$, $\mathcal{L}_2$, $E_{P_{2,6}}$,
$\mathcal{L}_6$, $E_{P_{3,6}}$, $\mathcal{L}_3$, $E_{A,3}$, $E_A$, $E_{A,1}$, $\mathcal{L}_1$, $E_{P_{1,4}}$, $\mathcal{L}_5$, $E_B$, $E_{B,5}$, $E_{P_{4,5}}$, $E_{B_1}$, $E_{B_6}$ and this choice of the basis corresponds to the components $\Theta_i$ $i=1,\ldots 12$ of the fiber $I_{12}$, the components $\Theta_j$, $j=0,1,2 $ of the fiber $I_4$, the zero section, and two sections of infinite order. If there was a section of order 2, then the set consisting of components of reducible fibers, of the zero section and of the sections of infinite order is a set of generators, but not a $\Z$-basis. Anyway, it is possible that there exists another elliptic fibration on $X_{(18,4,1)}$ with reducible fibers $I_4+I_{12}$ and admitting a section of order 2. Similarly, on the elliptic fibration with fibers $2I_8$ described in proof of Proposition \ref{prop: existence case 1} there is no section of order 2, but it is possible that there exists another elliptic fibration on $X_{(18,4,1)}$ with the same reducible fibers but admitting a section of order 2.\erem}\end{rem}
\begin{rem}\label{rem: symplectic involution non van geemen sarti}{\rm Let $\mathcal{E}_r$ be one of the elliptic fibrations of type a) (cf.\ Proposition \ref{prop: classification non finite MW}).  
For each section $P\in \MW(\mathcal{E}_r)$ such that the order of $P$ is not 2, there exists a symplectic involution $\alpha_P$ which preserves $\mathcal{E}_r$ but is not a van Geemen--Sarti involution of $\mathcal{E}_r$. Indeed let $\alpha_P:=\iota\beta_P$, where $\beta_P$ is the involution defined in Example \ref{ex: non symplectic involution Q-P}. The involutions $\iota$ and $\beta_P$ act respectively on the basis and on the fibers of the fibration, thus they commute and both are non--symplectic involutions. Hence $\alpha_P$ is a symplectic involution and it preserves $\mathcal{E}_r$. Clearly $\alpha_P$ is not a translation by a section of order 2 of the elliptic fibration $\mathcal{E}$. It is possible that $\alpha_P$ is a van Geemen--Sarti involution with respect to another elliptic fibration (cf.\ Remarks \ref{rem: vGT first} and \ref{rem: vGT is vGS}). \erem}\end{rem}
\begin{rem}\label{rem: vGT first}{\rm In \cite[Section 3.4]{vGT} a particular Morrison--Nikulin involution $i$ on the surface $X_{(19,3,1)}$ is analyzed in order to describe explicitly the 2-isogeny (defined over $\Q$)  between $X_{(19,3,1)}$ and the desingularization of $X_{(19,3,1)}/i$. In \cite{vGT} it is proved that the Morrison--Nikulin involution $i$ acts on the elliptic fibration with reducible fibers $I_{10}+I_8$ as the involution $\alpha_P$ described in Remark \ref{rem: symplectic involution non van geemen sarti} (we will show in Remark \ref{rem: vGT is vGS} that it is a van Geemen--Sarti involution with respect to another fibration). \erem}\end{rem}

\subsection{van Geemen--Sarti involutions case b)}
The classification of the van Geemen--Sarti involutions on the K3 surface $X_{(r,22-r,\delta)}$ in case b) coincides with the classification of all the sections of order 2 on elliptic fibrations on such a surface. This classification is contained in Propostion \ref{prop: classification finite MW}. Each van Geemen--Sarti involution induces an isogeny of order 2 between two K3 surfaces: describing the quotients of the surface $X_{(r,22-r,\delta)}$  by these van Geemen--Sarti involutions, we obtain a lot of isogenies between $X_{(r,22-r,\delta)}$ and other K3 surfaces. As an example we prove the following Proposition, then we summarize in tables  of Section \ref{subsec: quotient} similar results for all the van Geemen--Sarti involutions on the elliptic fibrations listed in Proposition \ref{prop: classification finite MW}. We observe that if the torsion part of the Mordell--Weil group is $(\Z/2\Z)^2$ then there are three sections of order 2, which are denoted by $t$, $u$ and $t+u$ (where $+$ is the group law in the Mordell--Weil group) and thus three distinct van Geemen--Sarti involutions, denoted by $\sigma_t$, $\sigma_u$, $\sigma_{t+u}$.

\begin{proposition}\label{prop: quotient elliptic fibration example} The K3 surface $X_{(19,3,1)}$ is 2-isogenous to a K3 surface $Y$ admitting an elliptic fibration with $\MW=(\Z/2\Z)^2$ and reducible fibers $2I_2^*+I_4+2I_2$. Let $t'$ and $u'$ be two sections of order 2 on $Y$, then $t'$ and $u'$ intersect precisely the following fiber components: $t'\cdot \Theta_1^1=t'\cdot \Theta_1^2=t'\cdot \Theta_2^3=t'\cdot\Theta_1^4=t'\cdot\Theta_1^5=1$, $u'\cdot \Theta_5^1=u'\cdot \Theta_5^2=u'\cdot \Theta_2^3=u'\cdot \Theta_0^4=u'\cdot \Theta_0^5=1$. The transcendental lattice of $Y$ is $\T_Y\simeq U(2)\oplus \langle 4\rangle$. The elliptic fibration on $Y$ has equation $$
w^2=x\left(x+\tau(\tau-1)^2\right)\left(x+\tau\left(\tau^2+2\tau(1-2d)+1\right)\right).
$$
The desingularization of the quotient $Y/\sigma_{t'}$ is $X_{(19,3,1)}$ and the 2-isogeny between $Y$ and $X_{(19,3,1)}$ is: $(x,w)\mapsto$ 
\small 
$$\left(\frac{\left(x^2+2\tau\left(\tau^2-2d\tau+1\right)x+
\tau^2\left(\tau-1\right)^2\left(\tau^2+2\tau(1-2d)+1\right)\right)}{4x}, \frac{w\left(x^2-\tau^2\left(\tau-1\right)^2\left(\tau^2+2\tau(1-2d)+1\right)\right)}{8x^2}\right).
$$
\normalsize
The desingularization of the quotient $Y/\sigma_{t'}$ is $X_{(19,3,1)}$.
\end{proposition}
\bprf It is well known (and easy to check directly) that the quotient of a fiber of type $I_{2k}^*$ by the translation by a section of order 2 is either $I_{k}^*$, if the section of order 2 meets one of the components $\Theta_{2k+3}$, $\Theta_{2k+4}$, or $I_{4k}^*$, if the section of order 2 meets the component $\Theta_1$. Similarly, the quotient of a fiber of type $I_{2k}$ by the translation by a section of order 2 is either $I_{k}$, if the section of order 2 meets the component $\Theta_{k}$, or $I_{4k}$, if the section of order 2 meets the component $\Theta_0$; the quotient of a fiber of type $I_1$ by the translation by a section of order 2 is a fiber of type $I_2$. The K3 surface $X_{(19,3,1)}$ admits an elliptic fibration with reducible fibers $2I_4^*+I_2+2I_1$ with a section of order 2 (see Section \ref{subsec: existence case b}, $r=19$ case 11)). Since the section $t$ of such an elliptic fibration meets the compontents $\Theta^1_7$, $\Theta_7^2$, $\Theta_0^3$ of the reducible fibers, the quotient elliptic fibration has fibers $2I_2^*+I_4+2I_2$. The elliptic fibration on $X_{(19,3,1)}$ has no sections of infinite order (see Proposition \ref{prop: classification finite MW}). Since the Mordell--Weil rank is invariant under isogeny, also the quotient elliptic fibration on $Y$ has no sections of infinite order.
By \cite{Sh} an elliptic fibration without sections of infinite order and with reducible fibers $2I_2^*+I_4+2I_2$ has Mordell--Weil group equal to $(\Z/2\Z)^2$. Let $t'$ and $u'$ be two generators of the Mordell--Weil group. Since the elliptic fibration on $Y$ is the quotient by the translation by a section of order 2, on $Y$ there exists a section of order 2, say $t'$, which gives the dual isogeny (i.e. such that the quotient by $\sigma_{t'}$ gives the original elliptic fibration). It follows that the section $t'$ meets the reducible fibers respectively in $\Theta_1^1$, $\Theta_1^2$, $\Theta_2^3$, $\Theta_1^4$, $\Theta_1^5$. By the height formula, it follows that the other two sections of order 2 meet the components in $\Theta^1_5$, $\Theta^2_5$, $\Theta^3_2$, $\Theta^4_0$, $\Theta^5_0$, and in $\Theta^1_6$, $\Theta^2_6$, $\Theta^3_0$, $\Theta^4_1$, $\Theta^5_1$ respectively. Thus one can assume that $u'$ is as in the statement.\\
The classes $b_1:=(\Theta_5^1+\Theta_6^1+\Theta_5^2+\Theta_6^2)/2$, $b_2:=(\Theta_1^2+\Theta_3^2+\Theta_5^2+\Theta_1^4)/2$, $b_3:=(\Theta_5^1+\Theta_6^1+\Theta_1^4)/2+ (\Theta_1^3+2\Theta_2^3+3\Theta_3^3)/4$ generate the discriminant group, $(\Z/2\Z)^2\times\Z/4\Z$, of $\NS(Y)$. The discriminant form (computed on these classes) has the property that $b_i^2\in\Z$, $i=1,2,3$ and $b_ib_j\in\frac{1}{2}\Z$, $i,j=1,2,3$. The discriminant form of $\T_Y$ is the opposite of the one of $\NS(Y)$, thus its generators satisfy the same conditions, moreover $rk(\T_Y)=l(\T_Y)$. This suffices to conclude that $\T_Y=L(2)$ for a certain even lattice $L$. The lattice $L$ has signature $(2,1)$ and has discriminant group $\Z/2\Z$ and thus satisfies \cite[Corollary 1.13.3]{Ni}. Hence $L$ is isometric to a unique lattice, which is $U\oplus\langle 2\rangle$. Thus $\T_Y\simeq U(2)\oplus \langle 4\rangle$. \\
Until now we did not use the explicit equation of the elliptic fibration on $X_{(19,3,1)}$, thus the previous arguments can be applied to any of the elliptic fibrations listed in Proposition \ref{prop: classification finite MW}. Now we observe that in Section \ref{subsec: existence r=19 tirvial MW} we proved that the following equation yields an elliptic fibration on $X_{(19,3,1)}$ with $2I_4^*+I_2+2I_1$ as fiber configuration: $w^2=x\left(x^2  -\tau\left(\tau^2-2d\tau+1\right)x+\tau^4(d-1)^2\right)$.
This is the equation of an elliptic curve defined over $\C[\tau]$ and the quotient by $\sigma_t$ is the quotient of this elliptic curve by the 2-torsion point $(0;0)$. By \cite[Pag. 79]{ST}, one immediately obtains the equation of the elliptic fibration on $Y$ and of the 2-isogeny.\erem

\begin{proposition} The van Geemen--Sarti involutions on $X_{(r,22-r,\delta)}$ in case b) are the ones associated to the sections of order 2 on the elliptic fibrations listed in Proposition \ref{prop: classification finite MW}. The quotient elliptic fibrations are listed in tables of Section \ref{subsec: quotient}. \\
Among these van Geemen--Sarti involutions those which are Morrison--Nikulin involutions are: $\sigma_t$ in cases $r=20$ cases 4), 6), 8), 9) and $r=19$ cases 8), 10b), 12), 13); $\sigma_u$ and $\sigma_{t+u}$ in case $r=20$ case 9).\\
The involution induced by $\iota$ on the quotient elliptic fibration acts as the identity on the N\'eron--Severi group in the following cases:\\
$r=20$ cases 3), 7) quotient by $\sigma_t$, $r=20$ case 8) quotient by $\sigma_u$ and $\sigma_{t+u}$;\\ 
$r=18$ cases 5b), 6), 14b), 15), 19) quotient by $\sigma_t$, $r=18$ cases 12), 16), 18) and 19) quotient by $\sigma_u$, 
$r=18$ cases 19), 20a) quotient by $\sigma_{t+u}$;\\ $r=16$ cases 6b), 9b), 11) quotient by $\sigma_t$.\end{proposition}
\bprf Proposition \ref{prop: classification finite MW} and Section \ref{sec: existence case 2)} give a complete classification of the elliptic fibrations of type b) on the surfaces $X_{(r,22-r,\delta)}$. Thus one obtains a complete classification of the van Geemen--Sarti involutions on these surfaces. The properties of the quotient elliptic fibrations can be proved as in Proposition \ref{prop: quotient elliptic fibration example}. In particular one computes the reducible fibers of the quotient from those of the original elliptic fibration and deduces the torsion part of the Mordell--Weil group by \cite{Sh} and by the fact that $\Z/2\Z\subset \MW$. Moreover, as in Proposition \ref{prop: quotient elliptic fibration example}, one proves that the Mordell--Weil rank is 0. In two cases ($r=18$ case 20a), quotient by $\sigma_t$ or $\sigma_u$ and $r=18$ case 20b) quotient by $\sigma_t$) the information given in \cite{Sh} and the inclusion $\Z/2\Z\subset \MW$ do not suffice to find the Mordell--Weil group: in both the cases it can be either $\Z/2\Z$ or $\Z/4\Z$. In case $r=18$ case 20a) an equation of the elliptic fibration on $X_{(18,4,1)}$ is $w^2=x\left(x+\tau\left(g\tau+1\right)\left(f+\tau\right)\right)\left(x+\tau^2\left(g+1\right)\left(f+1\right)\right)$ (cf.\ Section \ref{subsec: existence a 18}) and thus the equation of the elliptic fibration which is the quotient by $\sigma_t$ is \begin{equation}\label{eq:quotient non 4-torsion}w^2=x\left(x^2-2\tau\left(2fg\tau+\left(1+\tau\right)\left(g\tau+f\right)+2\tau\right)x+\tau^2\left(1-\tau\right)^2\left(g\tau-f\right)^2\right).\end{equation} In \cite[Table 1]{GS2} conditions to the equation of an  elliptic fibration to admit a 4-torsion section are given and thus we conclude that the quotient elliptic fibration \eqref{eq:quotient non 4-torsion} does not admit a 4-torsion section. Similarly, one can use the equations of the elliptic  fibrations  which are quotient by $\sigma_u$ of the case $r=18$, 20a) and quotient by $\sigma_t$ of the case $r=18$, 20b) to prove that in both cases the Mordell--Weil group is $\Z/2\Z$.\\ 
In the proof  of Proposition \ref{prop: van Geemen Sarti case a) and Morrison Nikulin} it is proved that there are no Morrison--Nikulin involutions on the surfaces $X_{(18,4,\delta)}$. With the same argument one shows that there are no Morrison--Nikulin involutions also on the surface $X_{(17,5,1)}$. Since the K3 surfaces with Picard number less than 17 do not admit Morrison--Nikulin involutions, we conclude that Morrison--Nikulin involutions can appear only in the cases $r=20$ or $r=19$. The discriminant group of the transcendental lattice of $X_{(r,22-r,\delta)}$ is $(\Z/2\Z)^{22-r}$ thus, if $\sigma_t$ is a Morrison--Nikulin involution then the discriminant group of the transcendental lattice of ${X_{(r,22-r,\delta)}/\sigma_t}$ is $(\Z/4\Z)^{22-r}$. This condition is satisfied only in cases $r=20$ cases 4), 6), 9) and $r=19$ cases 8), 10b), 12), 13) if one considers the quotient by $\sigma_t$, and in case $r=20$ case 8) if one considers the quotient by $\sigma_u$. 
Now one has to prove that in all these cases $\T_{X_{(r,22-r,\delta)}/\sigma_v}\simeq \T_{X_{(r,22-r,\delta)}}(2)$ where $v=t$ or $v=u$. One can make this computation directly (as in Proposition \ref{prop: quotient elliptic fibration example}) or refer to \cite[Cases 137, 177, 188, 286, Table 2]{SZ} if $r=20$.\\
The involution induced by $\iota$ on the quotient elliptic fibration is the hyperelliptic involution (i.e. of type b)). If it acts trivially on the N\'eron--Severi group, then the N\'eron--Severi group is 2-elementary. The unique quotient elliptic fibrations with this properties are the ones listed in the statement. Let us consider, for example, $r=16$ case 6b): the reducible fibers are $I_8^*+2I_2+6I_1$ and $\MW=\Z/2\Z$. The hyperelliptic involution fixes the zero section and the section $t'$ of order 2. Moreover, it fixes the bisection $b$ passing through the other points of order 2. Since the 2-torsion points are fixed, the curves $\Theta_0^1$, $\Theta_1^1$, $\Theta_{11}^1$, $\Theta_{12}^1$ are invariant. Thus $\Theta_2^1$ is invariant and the points $\Theta_0^1\cap \Theta_2^1$ and $\Theta_1^1\cap \Theta_2^1$  are fixed. But a non--symplectic involution does not admit isolated fixed points, thus $\Theta_2^1$ is a fixed curve ($\Theta_0^1$ is not fixed because it meets the zero section, which is fixed, and the fixed locus is smooth, analogously $\Theta_1^1$ meets $b$, thus it is not fixed). Similarly, one obtains that $\Theta_{2i}^1$, $i=1,2,3,4,5$ are fixed. Thus the hyperelliptic involution fixes 7 rational curves (2 sections and 5 components of the fiber of type $I_8^*$) and the curve $b$. This curve passes through two distinct points on each smooth fiber, on the fiber $I_8^*$ (it meets both the components $\Theta_1$ and $\Theta_{11}$) and on the fibers  of type $I_2$ (it meets both the components $\Theta_0$ and $\Theta_1$). It is tangent to the fibers of type $I_1$ and thus $b$ is a 2:1 cover of $\mathbb{P}^1$ (base of the fibration) branched along 6 points. Hence $b$ is a curve of genus 2 and the hyperelliptic involution fixes 7 rational curves and 1 curve of genus 2. Thus it acts as the identity on the N\'eron--Severi group if and only if the N\'eron Severi group is a 2-elementary lattice with rank 16 and length 2. These are the properties defining the N\'eron--Severi group of the quotient elliptic fibration, and this concludes the proof.\eprf

\begin{rem}{\rm The previous Proposition gives explicit isogenies among some K3 surfaces admitting a non--symplectic involution acting trivially on the N\'eron--Severi group, in particular between $X_{(16,6,1)}$ and $X_{(16,4,1)}$, $X_{(16,6,1)}$ and $X_{(16,2,1)}$, $X_{(18,4,0)}$ and $X_{(18,2,1)}$,  $X_{(18,4,1)}$ and $X_{(18,2,0)}$, $X_{(18,4,0)}$ and $X_{(18,0,0)}$. Moreover, there is an isogeny of order 2 between $X_{(18,4,0)}$ and $X_{(18,4,1)}$ given by the quotient of the elliptic fibration $r=18$ case 15) by $\sigma_t$.\erem}\end{rem}
\begin{rem}{\rm By the equation \eqref{eq: r=12} and the conditions $A=0$, $E=F=-1/2$, $H=-G$ one finds the equation of the elliptic fibration $r=16$ case 6b): it is $w^2=x(x^2+a(\tau)x+b(\tau))$ with $$\begin{array}{ll}a(\tau):=&(1-\tau)\left(8\tau\left(B+C\tau\right)+2(1-\tau)\left(D+\tau(1-3B-2C-2D)+\tau^2(B+D-1)\right)\right)\\ b(\tau):=&4\tau(B+C\tau)\left(1-\tau\right)^2\left(4G^2\tau(B+C\tau)+2G(1-\tau)\left(D+\tau\left(1-3B-2C-2D\right)+\right.\right.\\&\left.\left. \tau^2(B+D-1)\right)+(1-\tau)^2\left(B+C-1+\tau\right)\right).\end{array}$$ From this equation one finds an explicit equation for the isogeny described in \cite{CD:vGS involutions}, as in Proposition \ref{prop: quotient elliptic fibration example}.\erem}\end{rem}
\begin{rem}\label{rem: vGT is vGS}{\rm We showed in the previous Proposition that $\sigma_t$ is a Morrison--Nikulin involution in case $r=19$, 10b). One can easily identify the two copies of $E_8(-1)$ which are switched by this involution: one of them is $\mathcal{L}_1$, $E_{A,1}$, $E_A$, $E_{A,2}$, $E_{A,3}$, $\mathcal{L}_3$, $E_{P_{3,6}}$, $\mathcal{L}_6$. This shows that the involution $\sigma_t$ is the Morrison--Nikulin involution considered in \cite{vGT} (it suffices to identify $\mathcal{L}_1$ with $(x=0)$, $\mathcal{L}_2$ with $(x=-tz)$, $\mathcal{L}_3$ with $(z=0)$,  $\mathcal{L}_4$ with $(z=-1)$, $\mathcal{L}_5$ with $l_\infty$, $\mathcal{L}_6$ with $(x=-1)$, with the notation of \cite[Section 3.1]{vGT}).\erem}\end{rem}
\section{Appendix: Tables}\label{appendix}
\subsection{Elliptic fibrations case b): tables}\label{subsec: existence case b}
\scriptsize{
Case $r=20$ ($a=2$, $s=10$):
\begin{align*}
\begin{array}{|c|c|c|c|c|l|}
\hline
&\mbox{singular fibers}&m_1+m_2&n_1+2n_2&\MW&\mbox{components meeting}\\
&&&&&\mbox{ the torsion sections}\\
\hline
1)&2II^*& 2&0& \{1\}&\\
\hline
2)&II^*+I_6^*&0&2& \{1\}&\\
\hline
3)&2III^*+I_0^*&0&0&\Z/2\Z&t:\ \Theta^1_7,\ \Theta^2_7,\ \Theta_1^3\\
\hline
4)&III^*+I_6^*&1&3-(2m_1+3m_2)&\Z/2\Z&t:\ \Theta^1_7,\ \Theta^2_9,\ \Theta_0^3\\
\hline
5)&I_{14}^*&0&4&\{1\}&\\
\hline
6)&I_{12}^*&2&6-(2m_1+3m_2)&\Z/2\Z&t:\ \Theta^1_{15},\ \Theta^2_0,\ \Theta_0^3\\
\hline
7)&I_8^*+I_2^*&0&2&\Z/2\Z&t:\ \Theta^1_{11},\ \Theta^2_1\\
\hline
8)&2I_4^*&2&0&(\Z/2\Z)^2&t:\ \Theta^1_{7},\ \Theta^2_7,\ \Theta^3_0,\ \Theta^4_0\\
&&&&&u:\ \Theta^1_{8},\ \Theta^2_1,\ \Theta^3_1,\ \Theta^4_1\\
\hline
9)&3I_2^*&0&0&(\Z/2\Z)^2&t:\ \Theta^1_1,\ \Theta^2_5,\ \Theta^3_5\\
&&&&&u:\ \Theta^1_5,\ \Theta^2_6,\ \Theta^3_1\\
\hline
\end{array}\end{align*}
Case $r=19$ ($a=3$, $s=9$):
\begin{align*}
\begin{array}{|c|c|c|c|c|l|}
\hline
&\mbox{singular fibers}&m_1+m_2&n_1+2n_2&\MW&\mbox{components meeting}\\
&&&&&\mbox{ the torsion sections}\\
\hline
1)&II^*+III^*& 2&5-(2m_1+3m_2)& \{1\}&\\
\hline
2)&II^*+I_4^*&1&4-(2m_1+3m_2)& \{1\}&\\
\hline
3)&2III^*&3&0&\Z/2\Z&t:\ \Theta^1_7,\ \Theta^2_7,\ \Theta_1^3,\ \Theta_1^4,\ \Theta_0^5\\
\hline
4)&III^*+I_6^*&0&3&\{1\}&\\
\hline
5)&III^*+I_4^*&2&5-(2m_1+3m_2)&\Z/2\Z&t:\ \Theta_7^1,\ \Theta_7^2,\ \Theta_1^3,\ \Theta_0^4\\
\hline
6)&III^*+I_2^*+I_0^*&0&1&\Z/2\Z&t:\ \Theta_7^1, \Theta_5^2,\ \Theta_1^3\\
\hline
7)&I_{12}^*&1&6-(2m_1+3m_2)&\{1\}&\\
\hline
8)&I_{10}^*&3&8-(2m_1+3m_2)&\Z/2\Z&t:\ \Theta^1_{13},\ \Theta^2_1,\ \Theta^i_0,\ i=3,4\\
\hline
9)&I_8^*+I_0^*&1&4-(2m_1+3m_2)&\Z/2\Z&t:\ \Theta^1_{11},\ \Theta^2_1,\ \Theta^3_0\\
\hline
10a)&I_6^*+I_2^*&1&4-(2m_1+3m_2)&\Z/2\Z&t:\ \Theta^1_{9},\ \Theta^2_1,\ \Theta^3_1\\
\hline
10b)&I_6^*+I_2^*&1&4-(2m_1+3m_2)&\Z/2\Z&t:\ \Theta^1_{9},\ \Theta^2_5,\ \Theta^3_0\\
\hline
11)&2I_4^*&1&4-(2m_1+3m_2)&\Z/2\Z&t:\ \Theta^1_{7},\ \Theta^2_7,\ \Theta^3_0\\
\hline
12)&I_4^*+I_2^*&3&0&(\Z/2\Z)^2&t:\ \Theta^1_{7},\ \Theta^2_5,\ \Theta^3_1,\ \Theta^4_0,\ \Theta^5_0\\
&&&&&u:\ \Theta^1_1,\ \Theta^2_6,\ \Theta^3_1,\ \Theta^4_1,\ \Theta^5_1\\
\hline
13)&2I_2^*+I_0^*&1&0&(\Z/2\Z)^2&t:\ \Theta^1_{5},\ \Theta^2_5,\ \Theta^3_1,\ \Theta^4_0\\
&&&&&u:\ \Theta^1_1,\ \Theta^2_6,\ \Theta^3_3,\ \Theta^4_1\\
\hline
\end{array}\end{align*}
Case $r=18$ ($a=4$, $s=8$):
\begin{align*}
\begin{array}{|c|c|c|c|c|l|}
\hline
&\mbox{singular fibers}&m_1+m_2&n_1+2n_2&\MW&\mbox{components meeting}\\
&&&&&\mbox{ the torsion sections}\\
\hline
1)&II^*+I_2^*& 2&6-(2m_1+3m_2)& \{1\}&\\
\hline
2)&II^*+2I_0^*&0&2& \{1\}&\\
\hline
3)&2III^*&2&6-(2m_1+3m_2)&\{1\}&\\
\hline
4)&III^*+I_4^*&1&5-(2m_1+3m_2)&\{1\}&\\
\hline
5a)&III^*+I_2^*&3&7-(2m_1+3m_2)&\Z/2\Z&t:\ \Theta_7^1,\ \Theta_5^2,\ \Theta_1^3,\ \Theta_1^4,\ \Theta_0^5\\
\hline
5b)&III^*+I_2^*&3&7-(2m_1+3m_2)&\Z/2\Z&t:\ \Theta_7^1,\ \Theta_1^2,\ \Theta_1^3,\ \Theta_1^4,\ \Theta_1^5\\
\hline
6)&III^*+2I_0^*&1&3-(2m_1+3m_2)&\Z/2\Z&t:\ \Theta_7^1, \Theta_1^2,\ \Theta_1^3,\ \Theta_1^4\\
\hline
7)&I_{10}^*&2&8-(2m_1+3m_2)&\{1\}&\\
\hline
8)&I_{8}^*+I_0^*&0&4&\{1\}&\\
\hline
9)&I_8^*&4&10-(2m_1+3m_2)&\Z/2\Z&t:\ \Theta^1_{11},\ \Theta^2_1,\ \Theta^3_1,\ \Theta^4_0,\ \Theta^5_0\\
\hline
10)&I_6^*+I_2^*&0&4&\{1\}&\\
\hline
11)&I_6^*+I_0^*&2&6-(2m_1+3m_2)&\Z/2\Z&t:\ \Theta^1_{9},\ \Theta^2_1,\ \Theta^3_1,\ \Theta^4_0\\
\hline
12)&I_6^*&6&0&(\Z/2\Z)^2&t:\ \Theta^1_{9},\ \Theta^2_1,\ \Theta^3_1,\ \Theta^4_1,\ \Theta^5_0,\ \Theta^6_0,\ \Theta^7_0\\
&&&&&u:\ \Theta^1_{1},\ \Theta^2_1,\ \Theta^3_1,\ \Theta^4_1,\ \Theta^5_1,\ \Theta^6_1,\ \Theta^7_1\\
\hline
13)&2I_4^*&0&4&\{1\}&\\
\hline
14a)&I_4^*+I_2^*&2&6-(2m_1+3m_2)&\Z/2\Z&t:\ \Theta^1_{7},\ \Theta^2_5,\ \Theta^3_1,\ \Theta^4_0\\
\hline
14b)&I_4^*+I_2^*&2&6-(2m_1+3m_2)&\Z/2\Z&t:\ \Theta^1_{7},\ \Theta^2_1,\ \Theta^3_1,\ \Theta^4_1\\
\hline
15)&I_4^*+2I_0^*&0&2&\Z/2\Z&t:\ \Theta^1_{7},\ \Theta^2_1,\ \Theta^3_1\\
\hline
16)&I_4^*+I_0^*&4&0&(\Z/2\Z)^2&t:\ \Theta^1_{7},\ \Theta^2_1,\ \Theta^3_1,\ \Theta^4_1,\ \Theta^5_0,\ \Theta^6_0\\
&&&&&u:\ \Theta^1_{1},\ \Theta^2_3,\ \Theta^3_1,\ \Theta^4_1,\ \Theta^5_1,\ \Theta^6_1\\
\hline
17)&2I_2^*+I_0^*&0&2&\Z/2\Z&t:\ \Theta^1_{5},\ \Theta^2_5,\ \Theta^3_1\\
\hline
18)&I_2^*+2I_0^*&2&0&(\Z/2\Z)^2&t:\ \Theta^1_{5},\ \Theta^2_1,\ \Theta^3_1,\ \Theta^4_1,\ \Theta^5_0\\
&&&&&u:\ \Theta^1_{1},\ \Theta^2_3,\ \Theta^3_3,\ \Theta^4_1,\ \Theta^5_1\\
\hline
19)&4I_0^*&0&0&(\Z/2\Z)^2&t:\ \Theta^1_{1},\ \Theta^2_1,\ \Theta^3_1,\ \Theta^4_1\\
&&&&&u:\ \Theta^1_{3},\ \Theta^2_3,\ \Theta^3_3,\ \Theta^4_3\\
\hline
20a)&2I_2^*&4&0&(\Z/2\Z)^2&t:\ \Theta^1_{6},\ \Theta^2_6,\ \Theta^3_1,\ \Theta^4_1,\ \Theta^5_0,\ \Theta^6_0\\ 
&&&&&u:\ \Theta^1_{5},\ \Theta^2_5,\ \Theta^3_0,\ \Theta^4_0,\ \Theta^5_1,\ \Theta^6_1\\
\hline
20b)&2I_2^*&4&0&(\Z/2\Z)^2&t:\ \Theta^1_{5},\ \Theta^2_5,\ \Theta^3_1,\ \Theta^4_1,\ \Theta^5_0,\ \Theta^6_0\\ 
&&&&&u:\ \Theta^1_{1},\ \Theta^2_6,\ \Theta^3_0,\ \Theta^4_1,\ \Theta^5_1,\ \Theta^6_1\\
\hline
\end{array}\end{align*}
Case $r=17$ ($a=5$, $s=7$):
\begin{align*}
\begin{array}{|c|c|c|c|c|l|}
\hline
&\mbox{singular fibers}&m_1+m_2&n_1+2n_2&\MW&\mbox{components meeting}\\
&&&&&\mbox{ the torsion sections}\\
\hline
1)&II^*+I_0^*& 3&8-(2m_1+3m_2)& \{1\}&\\
\hline
2)&III^*+I_2^*&2&7-(2m_1+3m_2)& \{1\}&\\
\hline
3)&III^*+2I_0^*&0&3& \{1\}&\\
\hline
4)&III^*+I_0^*&4&9-(2m_1+3m_2)&\Z/2\Z&t:\ \Theta_7^1,\ \Theta_1^2,\ \Theta_1^3,\ \Theta_1^4,\ \Theta_1^5,\ \Theta_0^6\\
\hline
5)&I_8^*&3&10-(2m_1+3m_2)&\{1\}&\\
\hline
6)&I_6^*+I_0^*&1&6-(2m_1+3m_2)&\{1\}&\\
\hline
7)&I_6^*&5&12-(2m_1+3m_2)&\Z/2\Z&t:\ \Theta_9^1, \Theta_1^2,\ \Theta_1^3,\ \Theta_1^4,\ \Theta_0^5,\ \Theta_0^6\\
\hline
8)&I_4^*+I_2^*&1&6-(2m_1+3m_2)&\{1\}&\\
\hline
9)&I_4^*+I_0^*&3&8-(2m_1+3m_2)&\Z/2\Z&t:\ \Theta^1_7,\ \Theta^2_1,\ \Theta^3_1,\ \Theta^4_1,\ \Theta^5_0\\
\hline
10)&I_4^*&7&0&(\Z/2\Z)^2&t:\ \Theta^1_7,\ \Theta^2_1,\ \Theta^3_1,\ \Theta^4_1, \ \Theta^5_1,\ \Theta^6_0,\ \Theta^7_0,\ \Theta^8_0\\
&&&&&u:\ \Theta^1_1,\ \Theta^2_1,\ \Theta^3_0,\ \Theta^4_1, \ \Theta^5_1,\ \Theta^6_1,\ \Theta^7_1,\ \Theta^8_1\\
\hline
11)&I_2^*+2I_0^*&1&4-(2m_1+3m_2)&\Z/2\Z&t:\ \Theta^1_{5},\ \Theta^2_1,\ \Theta^3_1,\ \Theta^4_1\\
\hline
12)&I_2^*+I_0^*&5&0&(\Z/2\Z)^2&t:\ \Theta^1_{5},\ \Theta^2_1,\ \Theta^3_1,\ \Theta^4_1, \ \Theta^5_1,\ \Theta^6_0,\ \Theta^7_0\\
&&&&&u:\ \Theta^1_{1},\ \Theta^2_3,\ \Theta^3_0,\ \Theta^4_1, \ \Theta^5_1,\ \Theta^6_1,\ \Theta^7_1\\
\hline
13)&3I_0^*&3&0&(\Z/2\Z)^2&t:\ \Theta^1_{1},\ \Theta^2_1,\ \Theta^3_1,\ \Theta^4_1,\ \Theta^5_1,\ \Theta^6_0\\
&&&&&u:\ \Theta^1_3,\ \Theta^2_3,\ \Theta^3_3,\ \Theta^4_0,\ \Theta^5_1,\ \Theta^6_1\\
\hline
14a)&2I_2^*&3&8-(2m_1+3m_2)&\Z/2\Z&t:\ \Theta^1_{5},\ \Theta^2_5,\ \Theta^3_1,\ \Theta^4_1,\ \Theta^5_0\\
\hline
14b)&2I_2^*&3&8-(2m_1+3m_2)&\Z/2\Z&t:\ \Theta^1_{5},\ \Theta^2_1,\ \Theta^3_1,\ \Theta^4_1,\ \Theta^5_1\\
\hline
\end{array}\end{align*}
Case $r=16$ ($a=6$, $s=6$):
\begin{align*}
\begin{array}{|c|c|c|c|c|l|}
\hline
&\mbox{singular fibers}&m_1+m_2&n_1+2n_2&\MW&\mbox{components meeting}\\
&&&&&\mbox{ the torsion sections}\\
\hline
1)& II^*&6&14-(2m_1+3m_2)&\{1\}&\\
\hline
2)&III^*+I_0^*& 3&9-(2m_1+3m_2)& \{1\}&\\
\hline
3)&III^*&7&15-(2m_1+3m_2)& \Z/2\Z&t:\ \Theta^1_7,\ \Theta^2_1,\ \Theta^3_1,\ \Theta^4_1, \ \Theta^5_1,\ \Theta^6_1,\ \Theta^7_0,\ \Theta^8_0\\
\hline
4)&I_6^*&4&12-(2m_1+3m_2)& \{1\}&\\
\hline
5)&I_4^*+I_0^*&2&8-(2m_1+3m_2)&\{1\}&\\
\hline
6a)&I_4^*&6&14-(2m_1+3m_2)&\Z/2\Z&t:\ \Theta^1_7,\ \Theta^2_1,\ \Theta^3_1,\ \Theta^4_1, \ \Theta^5_1,\ \Theta^6_0,\ \Theta^7_0\\
\hline
6b)&I_4^*&6&14-(2m_1+3m_2)&\Z/2\Z&t:\ \Theta^1_1,\ \Theta^2_1,\ \Theta^3_1,\ \Theta^4_1, \ \Theta^5_1,\ \Theta^6_1,\ \Theta^7_1\\
\hline
7)&2I_2^*&2&8-(2m_1+3m_2)&\{1\}&\\
\hline
8)&I_2^*+2I_0^*&0&4&\{1\}&\\
\hline
9a)&I_2^*+I_0^*&4&10-(2m_1+3m_2)&\Z/2\Z&t:\ \Theta_5^1, \Theta_1^2,\ \Theta_1^3,\ \Theta_1^4,\ \Theta_1^5,\ \Theta_0^6\\
\hline
9b)&I_2^*+I_0^*&4&10-(2m_1+3m_2)&\Z/2\Z&t:\ \Theta_1^1, \Theta_1^2,\ \Theta_1^3,\ \Theta_1^4,\ \Theta_1^5,\ \Theta_1^6\\
\hline
10)&I_2^*&8&0&(\Z/2\Z)^2&t:\ \Theta_5^1, \Theta_1^i,\ i=2,\ldots,6,\ \Theta_0^7,\ \Theta_0^8,\ \Theta_0^9\\
&&&&&u:\ \Theta_1^1, \Theta_0^2,\ \Theta_0^3,\ \Theta_1^i,\ i=4,\ldots, 9\\
\hline
11)&3I_0^*&2&6-(2m_1+3m_2)&\Z/2\Z&t:\ \Theta^1_1,\ \Theta^2_1,\ \Theta^3_1,\ \Theta^4_1,\ \Theta^5_1,\\
\hline
12)&2I_0^*&6&0&(\Z/2\Z)^2&t:\ \Theta^1_1,\ \Theta^2_1,\ \Theta^3_1,\ \Theta^4_1,\ \Theta^5_1,\ \Theta_1^6,\ \Theta_0^7,\ \Theta_0^8,\\
&&&&&u:\ \Theta^1_3,\ \Theta^2_3,\ \Theta^3_0,\ \Theta^4_0,\ \Theta^5_1,\ \Theta_1^6,\ \Theta_1^7,\ \Theta_1^8,\\
\hline
\end{array}\end{align*}
Case $r=15$ ($a=7$, $s=5$):
\begin{align*}
\begin{array}{|c|c|c|c|c|l|}
\hline
&\mbox{singular fibers}&m_1+m_2&n_1+2n_2&\MW&\mbox{components meeting}\\
&&&&&\mbox{ the torsion sections}\\
\hline
1)& III^*&6&15-(2m_1+3m_2)&\{1\}&\\
\hline
2)&I_4^*&5&14-(2m_1+3m_2)& \{1\}&\\
\hline
3)&I_2^*+I_0^*&3&10-(2m_1+3m_2)&\{1\}&\\ 
\hline
4a)&I_2^*&7&16-(2m_1+3m_2)& \Z/2\Z&t:\ \Theta_5^1,\ \Theta_1^i,\ i=2,\ldots,6,\ \Theta_0^7,\ \Theta_0^8\\
\hline
4b)&I_2^*&7&16-(2m_1+3m_2)& \Z/2\Z&t:\ \Theta_1^i, i=1,\ldots,7,\ \Theta_0^8\\
\hline
5)&3I_0^*&1&6-(2m_1+3m_2)&\{1\}&\\
\hline
6)&2I_0^*&5&12-(2m_1+3m_2)&\Z/2\Z&t:\ \Theta^i_1,\ i=1,\ldots, 6,\ \Theta^7_0\\
\hline
7)&I_0^*&9&0&(\Z/2\Z)^2&t:\ \Theta^i_1,\ i=1,\ldots, 7, \Theta^8_0,\ \Theta^9_0,\ \Theta^{10}_0\\
&&&&&u:\ \Theta^1_3,\ \Theta^2_0,\ \Theta^3_0,\ \Theta^4_0, \Theta^i_1\ i=5,\ldots 10\\
\hline
\end{array}\end{align*}
Case $r=14$ ($a=8$, $s=4$):
\begin{align*}
\begin{array}{|c|c|c|c|c|l|}
\hline
&\mbox{singular fibers}&m_1+m_2&n_1+2n_2&\MW&\mbox{components meeting}\\
&&&&&\mbox{ the torsion sections}\\
\hline
1)&I_2^*&6&16-(2m_1+3m_2)&\{1\}&\\ 
\hline
2)&2I_0^*&4&12-(2m_1+3m_2)&\{1\}&\\
\hline
3)&I_0^*&8&18-(2m_1+3m_2)&\Z/2\Z&t:\ \Theta^i_1,\ i=1,\ldots, 7,\ \Theta^8_0,\ \Theta^9_0\\
\hline
4)&&12&0&(\Z/2\Z)^2&t:\ \Theta^i_1,\ i=1,\ldots, 8, \Theta^9_0, \Theta^{10}_0, \Theta^{11}_0, \Theta^{12}_0\\ 
&&&&&u:\ \Theta^i_0,\ i=1,\ldots, 4, \Theta^j_1,\ j=5,\ldots,12\\
\hline
\end{array}\end{align*}
Case $r=13$ ($a=9$, $s=3$):
\begin{align*}
\begin{array}{|c|c|c|c|c|l|}
\hline
&\mbox{singular fibers}&m_1+m_2&n_1+2n_2&\MW&\mbox{components meeting}\\
&&&&&\mbox{ the torsion sections}\\
\hline
1)&I_0^*&7&18-(2m_1+3m_2)&\{1\}&\\
\hline
2)&&11&24-(2m_1+3m_2)&\Z/2\Z&t:\ \Theta^i_1,\ i=1,\ldots, 8, \Theta^9_0, \Theta^{10}_0, \Theta^{11}_0\\
\hline
\end{array}\end{align*}

Case $r=12$ ($a=10$, $s=2$):
\begin{align*}
 m_1I_2+m_2III+n_1I_1+n_2II,\ m_1+m_2=10,\ n_1+2n_2=24-(2m_1+3m_2),\ \MW=\{1\}.
\end{align*}\\
}
\normalsize
\subsection{Quotient of elliptic fibrations case b): tables}\label{subsec: quotient}
\scriptsize
In the following tables we describe the elliptic fibration $\mathcal{E}/\sigma_t$. We refer to the examples given in Section \ref{sec: existence case 2)}, thus we assume $m_2=n_2=0$ as in Proposition \ref{prop: existence case b)}. The elliptic fibrations considered are listed in the first column, using the notation of the tables of Section \ref{subsec: existence case b}.
In the second column we indicate which van Geemen--Sarti involution we used to obtain the quotient surface: if the quotients by two distinct van Geemen--Sarti involutions have the same type of reducible fibers and the same torsion sections then we write them in a unique line, indicating both the van Geemen--Sarti involutions in the second column. In the third, fourth and fifth columns the quotient elliptic fibration is described. We denote by $t'$ the generator of the $\Z/2\Z$ in $\MW(\mathcal{E}/\sigma_t)$ such that $(\mathcal{E}/\sigma_t)/\sigma_{t'}$ is the original elliptic fibration on $X_{(r,22-r,\delta)}$. If $(\Z/2\Z)^2\subset \MW(\mathcal{E}/\sigma_t)$ we denote by $u'$ another generator of $\MW(\mathcal{E}/\sigma_t)$. In the fifth column the intersections between the torsion sections and the components of the reducible fibers are given: the reducible fibers (not the singular) are numbered according to the order given in the third column. In the sixth column we give the discriminant group of the N\'eron--Severi of the quotient surface. If the N\'eron--Severi group is 2-elementary, the quotient surface corresponds to a K3 surface of type $X_{(r,a,\delta)}$: in this case in the seventh column we give the invariants $(r,a,\delta)$ and, if $a= 22-r$, the number of the corresponding elliptic fibration with respect to the tables of Section 
\ref{subsec: existence case b}.\\
Quotients in case $r=20$:
\begin{align*}
\begin{array}{|c|c|c|c|l|c|c|}
\hline
&&\mbox{singular fibers}&\MW&&\NS^{\vee}/\NS&(r,a,\delta),\ case\\
\hline
3)&\sigma_t&2III^*+I_0^*&\Z/2\Z&t':\ \Theta_7^1,\ \Theta_7^2,\ \Theta_1^3&(\Z/2\Z)^2& (20,2,1),\ 3)\\
\hline
4)&\sigma_t&III^*+I_3^*+I_4+I_2&\Z/2\Z&t':\ \Theta_7^1,\ \Theta_1^2,\ \Theta_2^3&(\Z/4\Z)^2&\\
\hline
6)&\sigma_t&I_{6}^*+2I_4+2I_2&(\Z/2\Z)^2&t':\ \Theta_1^1,\ \Theta_2^2,\ \Theta_2^3,\ \Theta_1^4,\ \Theta_1^5&(\Z/4\Z)^2&\\
&&&&u':\Theta_9^7,\ \Theta_2^2,\ \Theta_0^3,\ \Theta_1^4,\ \Theta_0^5&&\\
\hline
7)&\sigma_t&2I_4^*+2I_2&(\Z/2\Z)^2&t': \Theta_1^1,\ \Theta_7^2,\ \Theta_1^3,\ \Theta_1^4&(\Z/2\Z)^2&(20,2,1),\ 8)\\
&&&&u':\Theta_7^1,\ \Theta_1^2,\ \Theta_1^3,\ \Theta_1^4&&\\
\hline
8)&\sigma_t&2I_2^*+2I_4&(\Z/2\Z)^2&t':\ \Theta_1^1,\ \Theta_1^2,\ \Theta_2^3,\ \Theta_2^4&(\Z/4\Z)^2&\\
&&&&u':\Theta_5^1,\ \Theta_5^2,\ \Theta_2^3,\ \Theta_0^4&&\\
\hline
8)&\sigma_u, \sigma_{t+u}&I_8^*+I_2^*+2I_1&\Z/2\Z&t':\ \Theta_{11}^1,\ \Theta_1^2&(\Z/2\Z)^2&(20,2,1),\ 7)\\
\hline
9)&\sigma_t,\ \sigma_u,\ \sigma_{t+u}&I_4^*+2I_1^*&\Z/2\Z&t':\ \Theta_7^1,\ \Theta_1^2,\ \Theta_1^3&(\Z/4\Z)^2&\\
\hline
\end{array}\end{align*}

Quotients in case $r=19$:
\begin{align*}
\begin{array}{|c|c|c|c|l|c|}
\hline
&&\mbox{singular fibers}&\MW&&\NS^{\vee}/\NS\\
\hline
3)&\sigma_t&2III^*+2I_1+I_4&\Z/2\Z&t':\ \Theta_7^1,\ \Theta_7^2,\ \Theta_2^3&\Z/4\Z\times (\Z/2\Z)^2\\
\hline
5)&\sigma_t&III^*+I_2^*+I_4+I_2+I_1&\Z/2\Z&t':\ \Theta_7^1,\ \Theta_1^2,\ \Theta_2^3,\ \Theta_1^4&\Z/4\Z\times (\Z/2\Z)^2\\
\hline
6)&\sigma_t&III^*+I_1^*+I_0^*+I_2&\Z/2\Z&t':\ \Theta_7^1,\ \Theta_1^2,\ \Theta_1^3,\ \Theta_1^4&\Z/4\Z\times (\Z/2\Z)^2\\
\hline
8)&\sigma_t&I_{5}^*+2I_4+2I_2+I_1&\Z/2\Z&t':\ \Theta_1^1,\ \Theta_2^2,\ \Theta_2^3,\ \Theta_1^4,\ \Theta_1^5&(\Z/4\Z)^3\\
\hline
9)&\sigma_t&I_4^*+I_0^*+I_4+2I_2&(\Z/2\Z)^2&t':\ \Theta_1^1,\ \Theta_1^2,\ \Theta_2^3,\ \Theta_1^4,\ \Theta_1^5&(\Z/4\Z)\times (\Z/2\Z)^2\\
&&&&u':\ \Theta_7^1,\ \Theta_3^2,\ \Theta_2^3,\ \Theta_0^4,\ \Theta_0^5&\\
\hline
10a)&\sigma_t&I_3^*+I_4^*+2I_2+I_1&\Z/2\Z&t':\ \Theta_1^1,\ \Theta_7^2,\ \Theta_1^3,\ \Theta_1^4&\Z/4\Z\times (\Z/2\Z)^2\\
\hline
10b)&\sigma_t&I_3^*+I_1^*+I_4+2I_2&\Z/2\Z&t':\ \Theta_1^1,\ \Theta_1^2,\ \Theta_2^3,\ \Theta_1^4,\ \Theta_1^5&(\Z/4\Z)^3\\
\hline
11)&\sigma_t&2I_2^*+I_4+2I_2&(\Z/2\Z)^2&t':\ \Theta_1^1,\ \Theta_1^2,\ \Theta_2^3,\ \Theta_1^4,\ \Theta_1^5&\Z/4\Z\\
&&&&u':\ \Theta_5^1,\ \Theta_5^2,\ \Theta_2^3,\ \Theta_0^4,\ \Theta_0^5&\\
\hline
12)&\sigma_{t}&I_2^*+I_1^*+2I_4+I_1&\Z/2\Z&t':\ \Theta_{1}^1,\ \Theta_1^2,\ \Theta_2^3,\ \Theta_2^4&(\Z/4\Z)^3\\
\hline
12)&\sigma_u&I_8^*+I_1^*+3I_1&\Z/2\Z&t':\ \Theta_{11}^1,\ \Theta_1^2&\Z/4\Z\\
\hline
12)&\sigma_{t+u}&I_2^*+I_4^*+2I_1+I_4&\Z/2\Z&t':\ \Theta_1^1,\ \Theta_7^2,\ \Theta_2^3&\Z/4\Z\times (\Z/2\Z)^2\\
\hline
13)&\sigma_{t}&2I_1^*+I_0^*+I_4&\Z/2\Z&t':\ \Theta_1^1,\ \Theta_1^2,\ \Theta_1^3,\ \Theta_2^4&(\Z/4\Z)^3\\
\hline
13)&\sigma_u,\ \sigma_{t+u}&I_1^*+I_4^*+I_0^*+I_1&\Z/2\Z&t':\ \Theta_1^1,\ \Theta_7^2,\ \Theta_1^3&\Z/4\Z\times (\Z/2\Z)^2\\
\hline
\end{array}\end{align*}

Quotients of case $r=18$ (in case 12) we indicate with $v'$ the 4-torsion section; in this case the section $t'$ is 2 times $v'$):
\begin{align*}
\begin{array}{|c|c|c|c|l|c|c|}
\hline
&&\mbox{singular fibers}&\MW&&\NS^{\vee}/\NS&(r,a,\delta),\ case\\
\hline
5a)&\sigma_t&III^*+I_1^*+2I_1+I_4+I_2&\Z/2\Z&t':\ \Theta_7^1,\ \Theta_1^2,\ \Theta_2^3,\ \Theta_1^4&(\Z/4\Z)^2&\\
\hline
5b)&\sigma_t&III^*+I_4^*+3I_1+I_2&\Z/2\Z&t':\ \Theta_7^1,\ \Theta_7^2,\ \Theta_1^3&(\Z/2\Z)^2&(18,2,1)\\
\hline
6)&\sigma_t&III^*+2I_0^*+I_1+I_2&\Z/2\Z&t':\ \Theta_7^1,\ \Theta_1^2,\ \Theta_1^3,\ \Theta_1^4&(\Z/2\Z)^4&(18,4,1),\ 6)\\
\hline
9)&\sigma_t&I_4^*+2I_1+2I_4+2I_2&\Z/2\Z&t':\ \Theta_1^1,\ \Theta_2^2,\ \Theta_2^3,\ \Theta_1^4,\ \Theta_1^5&(\Z/4\Z)^2\times (\Z/2\Z)^2&\\
\hline
11)&\sigma_t&I_3^*+I_0^*+I_1+I_4+2I_2&\Z/2\Z&t':\ \Theta_1^1,\ \Theta_1^2,\ \Theta_2^3,\ \Theta_1^4,\ \Theta_1^5&(\Z/4\Z)^2\times (\Z/2\Z)^2&\\
\hline
12)&\sigma_t,\ \sigma_{t+u}&I_3^*+3I_1+3I_4&\Z/4\Z&v':\ \Theta_6^1,\ \Theta_1^2,\ \Theta_1^3,\ \Theta_1^4&(\Z/4\Z)^2&\\
\hline
12)&\sigma_u&I_{12}^*+6I_1&\Z/2\Z&t':\ \Theta_{15}^1&&(18,0,0)\\
\hline
14a)&\sigma_t&I_2^*+I_1^*+I_1+I_4+2I_2&\Z/2\Z&t':\ \Theta_1^1,\ \Theta_1^2,\ \Theta_2^3,\ \Theta_1^4,\ \Theta_1^5&(\Z/4\Z)^2\times (\Z/2\Z)^2&\\
\hline
14b)&\sigma_t&I_2^*+I_4^*+2I_1+2I_2&\Z/2\Z&t':\ \Theta_1^1,\ \Theta_7^2,\ \Theta_1^3,\ \Theta_1^4&(\Z/2\Z)^4&(18,4,1),\ 14b)\\
\hline
15)&\sigma_t&I_2^*+2I_0^*+2I_2&(\Z/2\Z)^2&t':\ \Theta_1^1,\ \Theta_1^2,\ \Theta_1^3,\ \Theta_1^4,\ \Theta_1^5&(\Z/2\Z)^4&(18,4,1),\ 18)\\
&&&&u':\ \Theta_5^1,\ \Theta_3^2,\ \Theta_3^3,\ \Theta_1^3,\ \Theta_0^4&&\\
\hline
16)&\sigma_t,\ \sigma_{t+u}&I_2^*+I_0^*+2I_1+2I_4&\Z/2\Z&t':\ \Theta_1^1,\ \Theta_1^2,\ \Theta_2^3,\ \Theta_2^4&(\Z/4\Z)^4\times (\Z/2\Z)^2&\\
\hline
16)&\sigma_u&I_8^*+I_0^*+4I_1&\Z/2\Z&t':\ \Theta_{11}^1,\ \Theta_1^1&(\Z/2\Z)^2&(18,2,0)\\
\hline
17)&\sigma_t&2I_1^*+I_0^*+2I_2&\Z/2\Z&t':\ \Theta_1^1,\ \Theta_1^2,\ \Theta_1^3,\ \Theta_1^4,\ \Theta_1^5&(\Z/4\Z)^2\times (\Z/2\Z)^2&\\
\hline
18)&\sigma_t,\ \sigma_{t+u}&I_1^*+2I_0^*+I_1+I_4&\Z/2\Z&t':\ \Theta_1^1,\ \Theta_1^2,\ \Theta_1^3,\ \Theta_2^4&(\Z/4\Z)^2\times (\Z/2\Z)^2&\\
\hline
18)&\sigma_u&I_4^*+2I_0^*+2I_1&\Z/2\Z&t':\ \Theta_7^1,\ \Theta_1^2,\ \Theta_1^3&(\Z/2\Z)^4&(18,4,0),\ 15)\\
\hline
19)&\sigma_t,\ \sigma_u\,\sigma_{t+u}&4I_0^*&(\Z/2\Z)^2&t':\ \Theta_1^1,\ \Theta_1^2,\ \Theta_1^3,\ \Theta_1^4&(\Z/2\Z)^4&(18,4,0),\ 19)\\
&&&&u':\ \Theta_3^1,\ \Theta_3^2,\ \Theta_3^3,\ \Theta_3^4&&\\
\hline
20a)&\sigma_t,\ \sigma_u&2I_1^*+2I_1+2I_4&\Z/2\Z&t':\ \Theta_1^1,\ \Theta_1^2,\ \Theta_2^3,\ \Theta_2^4&(\Z/4\Z)^2\times (\Z/2\Z)^2&\\
\hline
20a)&\sigma_{t+u}&2I_4^*+4I_1&\Z/2\Z&t':\ \Theta_7^1,\ \Theta_7^2&(\Z/2\Z)^2&(18,2,0)\\
\hline
20b)&\sigma_{t}&2I_1^*+2I_4+2I_1&\Z/2\Z&t':\ \Theta_1^1,\ \Theta_1^2,\ \Theta_2^3,\ \Theta_2^4&(\Z/4\Z)^2\times (\Z/2\Z)^2&\\
\hline
20b)&\sigma_u,\ \sigma_{t+u}&I_1^*+I_4^*+3I_1+I_4&\Z/2\Z&t':\ \Theta_1^1,\ \Theta_7^2,\ \Theta_2^3&(\Z/4\Z)^2&\\
\hline
\end{array}\end{align*}

Quotients in case $r=17$:
\begin{align*}
\begin{array}{|c|c|c|c|l|c|}
\hline
&&\mbox{singular fibers}&\MW&&\NS^{\vee}/\NS\\
\hline
4)&\sigma_t&III^*+I_0^*+3I_1+I_4+I_2&\Z/2\Z&t':\ \Theta_7^1,\ \Theta_1^2,\ \Theta_2^3,\ \Theta_1^4&\Z/4\Z\times (\Z/2\Z)^2\\
\hline
7)&\sigma_t&I_3^*+3I_1+2I_4+2I_2&\Z/2\Z&t':\ \Theta_1^1,\ \Theta_2^2,\ \Theta_2^3,\ \Theta_1^4,\ \Theta_1^5&(\Z/4\Z)^3\\
\hline
9)&\sigma_t&I_2^*+I_0^*+2I_1+I_4+2I_2&\Z/2\Z&t':\ \Theta_1^1,\ \Theta_1^2,\ \Theta_2^3,\ \Theta_1^4,\ \Theta_1^5&\Z/4\Z\times (\Z/2\Z)^4\\
\hline
10)&\sigma_t&I_2^*+4I_1+3I_4&\Z/2\Z&t':\ \Theta_1^1,\ \Theta_2^2,\ \Theta_2^3,\ \Theta_2^4&(\Z/4\Z)^3\\
\hline
10)&\sigma_u&I_8^*+6I_1+I_4&\Z/2\Z&t':\ \Theta_{11}^1,\ \Theta_2^2&\Z/4\Z\\
\hline
11)&\sigma_t&I_1^*+2I_0^*+I_1+2I_2&\Z/2\Z&t':\ \Theta_1^1,\ \Theta_1^2,\ \Theta_1^3,\ \Theta_1^4,\ \Theta_1^5&\Z/4\Z\times (\Z/2\Z)^4\\
\hline
12)&\sigma_t,\ \sigma_{t+u}&I_1^*+I_0^*+3I_1+2I_4&\Z/2\Z&t':\ \Theta_1^1,\ \Theta_1^2,\ \Theta_2^3,\ \Theta_2^4&(\Z/4\Z)^3\\
\hline
12)&\sigma_u&I_4^*+I_0^*+I_4+4I_1&\Z/2\Z&t':\ \Theta_7^1,\ \Theta_1^2,\ \Theta_2^3&\Z/4\Z\times (\Z/2\Z)^2\\
\hline
13)&\sigma_t,\ \sigma_u,\ \sigma_{t+u}&3I_0^*+2I_1+I_4&\Z/2\Z&t':\ \Theta_1^1,\ \Theta_1^2,\ \Theta_1^3,\ \Theta_2^4&\Z/4\Z\times (\Z/2\Z)^4\\
\hline
14a)&\sigma_t&2I_1^*+2I_1+I_4+2I_2&\Z/2\Z&t':\ \Theta_1^1,\ \Theta_1^2,\ \Theta_2^3,\ \Theta_1^4,\ \Theta_1^5&(\Z/4\Z)^3\\
\hline
14b)&\sigma_t&I_1^*+I_4^*+3I_1+2I_2&\Z/2\Z&t':\ \Theta_1^1,\ \Theta_7^2,\ \Theta_1^3,\ \Theta_1^4&\Z/4\Z\times (\Z/2\Z)^2\\
\hline
\end{array}\end{align*}

Quotients in case $r=16$:
\begin{align*}
\begin{array}{|c|c|c|c|l|c|c|}
\hline
&&\mbox{singular fibers}&\MW&&\NS^{\vee}/\NS&(r,a,\delta),\ case\\
\hline
3)&\sigma_t&III^*+5I_1+2I_4+I_2&\Z/2\Z&t':\ \Theta_7^1,\ \Theta_2^2,\ \Theta_2^3,\ \Theta_1^4&(\Z/4\Z)^2&\\
\hline
6a)&\sigma_t&I_2^*+4I_1+2I_4+2I_2&\Z/2\Z&t':\ \Theta_1^1,\ \Theta_2^2,\ \Theta_2^3,\ \Theta_1^4,\ \Theta_1^5&(\Z/4\Z)^2\times (\Z/2\Z)^2&\\
\hline
6b)&\sigma_t&I_8^*+6I_1+2I_2&\Z/2\Z&t':\ \Theta_{11}^1,\ \Theta_1^2,\ \Theta_1^3&(\Z/2\Z)^2&(16,2,1)\\
\hline
9a)&\sigma_t&I_1^*+I_0^*+3I_1+I_4+2I_2&\Z/2\Z&t':\ \Theta_1^1,\ \Theta_1^2,\ \Theta_2^3,\ \Theta_1^4,\ \Theta_1^5&(\Z/4\Z)^2\times (\Z/2\Z)^2&\\
\hline
9b)&\sigma_t&I_4^*+I_0^*+4I_1+2I_2&\Z/2\Z&t':\ \Theta_7^1,\ \Theta_1^2,\ \Theta_1^3,\ \Theta_1^4&(\Z/2\Z)^4&(16,4,1)\\
\hline
10)&\sigma_t,\ \sigma_{t+u}&I_1^*+5I_1+3I_4&\Z/2\Z&t':\ \Theta_1^1,\ \Theta_2^2,\ \Theta_2^3,\ \Theta_2^4&(\Z/4\Z)^3&\\
\hline
10)&\sigma_u&I_4^*+2I_4+6I_1&\Z/2\Z&t':\ \Theta_7^1,\ \Theta_2^2,\ \Theta_2^3&(\Z/4\Z)^2&\\
\hline
11)&\sigma_t&3I_0^*+2I_1+2I_2&\Z/2\Z&t':\ \Theta_1^1,\ \Theta_1^2,\ \Theta_1^3,\ \Theta_1^4,\ \Theta_1^5&(\Z/2\Z)^6&(16,6,1),\ 11)\\
\hline
12)&\sigma_t,\ \sigma_u,\ \sigma_{t+u}&2I_0^*+4I_1+2I_4&\Z/2\Z&t':\ \Theta_1^1,\ \Theta_1^2,\ \Theta_2^3,\ \Theta_2^4&(\Z/4\Z)^2\times (\Z/2\Z)^2&\\
\hline
\end{array}\end{align*}

Quotients in case $r=15$:
\begin{align*}
\begin{array}{|c|c|c|c|l|c|}
\hline
&&\mbox{singular fibers}&\MW&&\NS^{\vee}/\NS\\
\hline
4a)&\sigma_t&I_1^*+5I_1+2I_4+2I_2&\Z/2\Z&t':\ \Theta_1^1,\ \Theta_2^2,\ \Theta_2^3,\ \Theta_1^4,\ \Theta_1^5&(\Z/4\Z)^3\\
\hline
4b)&\sigma_t&I_4^*+6I_1+I_4+2I_2&\Z/2\Z&t':\ \Theta_7^1,\ \Theta_2^2,\ \Theta_1^3,\ \Theta_1^4&\Z/4\Z\times (\Z/2\Z)^2\\
\hline
6)&\sigma_t&2I_0^*+4I_1+I_4+2I_2&\Z/2\Z&t':\ \Theta_1^1,\ \Theta_1^2,\ \Theta_2^3,\ \Theta_1^4,\ \Theta_1^2&\Z/4\Z\times (\Z/2\Z)^4\\
\hline
7)&\sigma_t,\ \sigma_u,\ \sigma_{t+u}&I_0^*+6I_1+3I_4&\Z/2\Z&t':\ \Theta_1^1,\ \Theta_2^2,\ \Theta_2^3,\ \Theta_2^4&(\Z/4\Z)^3\\
\hline
\end{array}\end{align*}

Quotients in case $r=14$:
\begin{align*}\label{table: quotient non trivial \MW 14}
\begin{array}{|c|c|c|c|l|c|}
\hline
&&\mbox{singular fibers}&\MW&&\NS^{\vee}/\NS\\
\hline
3)&\sigma_t&I_0^*+6I_1+2I_4+2I_2&\Z/2\Z&t':\ \Theta_1^1,\ \Theta_2^2,\ \Theta_2^3,\ \Theta_1^4,\ \Theta_1^5&(\Z/4\Z)^2\times (\Z/2\Z)^2\\
\hline
4)&\sigma_t,\ \sigma_u,\ \sigma_{t+u}&8I_1+4I_4&\Z/2\Z&t':\ \Theta_2^1,\ \Theta_2^2,\ \Theta_2^3,\ \Theta_2^4&(\Z/4\Z)^3\\
\hline
\end{array}\end{align*}

Quotient in case $r=13$:
\begin{align*}
\begin{array}{|c|c|c|c|l|c|}
\hline
&&\mbox{singular fibers}&\MW&&\NS^{\vee}/\NS\\
\hline
2)&\sigma_t&8I_1+3I_4+2I_2&\Z/2\Z&t':\ \Theta_2^1,\ \Theta_2^2,\ \Theta_2^3,\ \Theta_1^4,\ \Theta_1^5&(\Z/4\Z)^3\\
\hline
\end{array}\end{align*}
\normalsize

\bibliographystyle{amsplain}

\end{document}